\documentclass[12pt]{report}
\pdfoutput=1 
\usepackage{graphicx}
\usepackage{latexsym}
\usepackage{amsmath}
\usepackage{amsthm}
\usepackage{amssymb}
\usepackage{bm}
\usepackage{subfig}

\usepackage{algpseudocode}
\usepackage{algorithm}

\newtheorem{theorem}{Theorem}[section]
\newtheorem{lemma}[theorem]{Lemma}

\theoremstyle{definition}
\newtheorem{definition}[theorem]{Definition}

\newtheorem{p}{Problem}[section]

\theoremstyle{remark}
\newtheorem{remark}[theorem]{Remark}

\begin{document}

\marginparsep 0pt
\textwidth 15.0 truecm

%\setcounter{page}{0}
%\begin{document}

\setlength{\baselineskip}{0.780cm}

\pagestyle{empty}
\begin{center}
\Large{ \bf Parametric FEM for Shape Optimization applied to Golgi Stack}
\end{center}

\vspace{20mm}
\begin{center}
CHEN, Xinshi
\end{center}

\vspace{20mm}
\begin{center}
A Thesis Submitted in Partial Fulfilment \\
of the Requirements for the Degree of \\
 Master of Philosophy \\
 in \\
 Mathematics
\end{center}

\vspace{20mm}
\begin{center}
The Chinese University of Hong Kong\\
July 2017
\end{center}

\newpage
\begin{titlepage}\centering
\vspace*{\fill}
\begin{center}
\underline{Thesis Assessment Committee}\\
$ $\\
Professor CHAN Hon Fu Raymond (Chair)\\
Professor CHUNG Tsz Shun Eric (Thesis Supervisor)\\
Professor LUI Lok Ming (Committee Member)\\
Professor KIM Hyea Hyun (External Examiner)	
\end{center}
\vspace*{\fill}
\end{titlepage}

\newpage
\setcounter{page}{1}
\pagestyle{myheadings}
\markright{Parametric FEM for Shape Optimization applied to Golgi Stacks}

\noindent
{\Huge {\bf Abstract}}
\vspace{1.2cm}

\noindent

The thesis is about an application of the shape optimization to the morphological evolution of Golgi stack. Golgi stack consists of multiple layers of cisternae. It is an organelle in the biological cells. Inspired by the Helfrich Model \cite{Helfrich}, which is a model for vesicles typically applied to biological cells, a new model specially designed for Golgi stack is developed and then implemented using FEM in this thesis.

In the Golgi model, each cisternae of the Golgi stack is viewed as a closed vesicle without topological changes, and our model is adaptable to both single-vesicle case and multiple-vesicle case.  The main idea of the math model is to minimize the elastic energy(bending energy) of the vesicles, with some constraints designed regarding to the biological properties of Golgi stack. With these constraints attached to the math model, we could extend this model to an obstacle-type problem. Hence, in the thesis, not only the simulations of Golgi stack are shown, some interesting examples without biological meanings are also demonstrated. Also, as multiple cisternaes are considered as a whole, this is also a model handling multiple objects.

 A set of numerical examples is shown to compare with the observed shape of Golgi stack, so we can lay down some possible explanations to the morphological performance of trans-Golgi cisternae.

%The second part of the thesis is on photoacoustic tomography. We apply the Mixed GMsFEM for numerical approximation of the wave propagation in this inverse problem. This paper is organized as follows. First, some background materials about the Neumann-based algorithm of this inverse problem will be presented. Second, the reconstruction method with the implementation details will be presented. Third, the explanation of the Mixed GMsFEM is included. Finally, various numerical examples will be shown as the demonstration of the method.
\newpage

\begin{center}
{\large {\bf ACKNOWLEDGMENTS}}
\end{center}
\vspace{5mm}

I would like to express my sincere gratitude to Prof. Eric CHUNG, who took me as his M.Phil. student two years ago. It is with this chance he gave me that I could participate this awesome graduate program in Department of Mathematics in CUHK. His constant guidance on my thesis research through the past two years is of great importance to me. Besides, I wish to express my appreciation to Prof. Raymond CHAN, who provided me with many opportunities of presenting my work in some meetings and symposium. Here I extend my gratitude to Prof. LUI Lok Ming and Prof. KIM Hyea Hyun, for spending their precious time on the examination of my thesis research.

Besides, I want to thank all the faculties, staffs and colleagues that I met in CUHK Math Department, for their various kinds of supports to my two years' study. Especially, I want to thank my group-mates including Dota Chi Yeung LAM, Tony Siu Wun CHEUNG, Nina Yue QIAN, Ivan Tak Shing AU YEUNG, Simon Sai Mang PUN, John Ming Fai LAM, Tommy Chor Hung LI, and quasi-group-mate John Yufei ZHANG, who give me countless help and bring me lots of fun everyday.

Last but not least, I would like to delicate this thesis to my parents, Mr. and Mrs. CHEN, and my boyfriend, Rui, for their unconditional company and support.

\newpage
{\large {\tableofcontents}}

%%% notation
\chapter*{Abbreviations and Notations}
\begin{table}[h]
\label{abbrev}
\begin{tabular}{ll}
FEM & Finite Element Method \\
SDG & Staggered Discontinuous Galerkin  \\
%GMsFEM & Generalized Multiscale Finite Element Method  \\
  &  
\end{tabular}
\end{table}

\begin{table}[h]
\label{notate}
\begin{tabular}{ll}
$\mathbb{R}$ & the set of real numbers \\
$\Gamma$ & a hypersurface of dimension $k-1$ embedded in $\mathbb{R}^k$  \\
${\boldsymbol{\nu}}$ & outer unit normal vector field  \\
$\|\cdot\|$  & Euclidean norm \\
$\mathbb{P}_2(K)$ & The space of polynomials of degree $\leq$ 2 defined over the set $K$\\
$C^1$ & continuously differentiable functions\\
$D(\phi)$ & $D(\phi)_{ij}=[\nabla_\Gamma \phi_j]_i +[\nabla_\Gamma \phi_i]_j$\\
$\nabla_\Gamma \phi$ & tangential gradient of $\phi$\\
$\mathcal{D}_{\Gamma}(\boldsymbol{\phi})$ & tangential Jacobian matrix of $\mathbf{w}$\\
$\text{div}_\Gamma \boldsymbol{\phi}$ & tangential divergence of $\mathbf{w}$\\
$\Delta_\Gamma \phi$ & tangential Laplace of $\phi$\\
\end{tabular}
\end{table}

%%% start Golgi project

\chapter{Introduction}

The thesis project is a mathematical modeling on the shape evolution of Golgi stack. This project was proposed by Prof. Kang from Life Science Department, and the final goal of this project is to mimic the growing process of Golgi stack so as to give a potential explanation to some special properties of the observed shapes of Golgi stack mathematically. The following image is shown to give a first impression of the Golgi stack \cite{golgi_photo}.
\begin{figure}[h]
  \centering
    \includegraphics[width=0.7\textwidth]{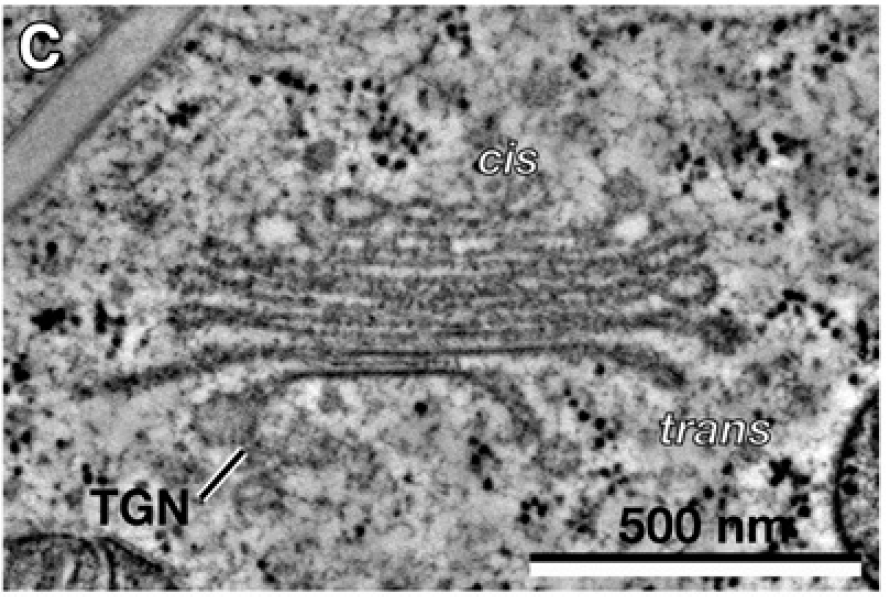}
      \caption{An ET slice image \cite{golgi_photo}. Scale bar: 500nm.}
\end{figure}

A Golgi stack consists of multiple layers. Each layer is a Golgi cisternae. The cisternae in different stages has different morphological performances. Those stages are called cis-, med- and trans-. For the cis-Golgi, the cisternal assembly is in process, but in the trans-Golgi, the assembly is finished. Our model mainly concerns with trans-Golgi cisternae, for which we can view each layer(cisternae) as a vesicle with closed surface without topological changes. 

Our math model is constructed by forming a series of geometric evolution equations. Suppose $\Gamma$ is the surface of a vesicle, which represents one Golgi cisternae. There are four shape functionals considered in our model defined on $\Gamma$. First, the Willmore energy $W(\Gamma)=\int_\Gamma h^2$, where $h$ is the mean curvature. It is believed that the shapes of the biomembranes are closely relevant to the elastic energy. Because Willmore energy is equivalent to the elastic energy under certain condition, which is explained in Section 1.1.1, we regard the Willmore energy as the dominant energy in our model of Golgi stack. Second, the area functional $A(\Gamma)=\int_\Gamma 1$. The model include this functional to serve for the surface area constraint. Because the number of molecules is believed to be fixed in our model, the surface area of the shape is conserved. Third, a heaviside functional $H(\Gamma)=\int_\Gamma 1_{\bf B}$, where $1_{\bf B}$ is an indicator function. To mimic the barriers (the intercisternal elements for Golgi stack) above and belong each Golgi layer, we form this new functional by integrating the heaviside function. Note that the set ${\bf B}$ represents the region of the barriers/obstacles. Fourth, a distance functional $D(\Gamma)$. When multiple vesicles are considered as a whole, the intersection of them should be avoided. Hence, this distance functional is included to handle the relations among those objects. We couple these functionals to form the model. Then we implement it by parametric finite element method, using Matlab.

The followings are the main results of the Golgi simulation. First, for a single layer of trans-cisternae, we mimic the evolution process of the cisternae when the protein vesicles bud from the marginal part of the cisternae. The numerical result explains that it is the barriers placed above and below the cisternae that inhibit the vertical expansion of the cisterae. Second, it is revealed from our result that the swelling of the marginal regions of the trans-cisternae thins its central part to decrease the elastic energy. Third, we also mimic the shape evolution of Golgi with multiple cisternae stacked upon each other. All these simulations are done by mathematical computation, which may give a potential explanation to the mechanism of the Golgi stack.

What's more, the model can also be extended to more general cases without biological meanings. In Section 5.1, many examples of the applications of this model are demonstrated.

\section{Overall Introduction and Motivation}
%% bio background of the model
In this section, I describe the construction of the math models with explanations on the corresponding biological properties of Golgi stack. From now on, we use the geometric surface $\Gamma\subset \mathbb{R}^k$ to describe the surface of a single Golgi cisternae, which is assumed to be a closed vesicle without topological changes. For the case including multiple layers of Golgi cisternaes in one model, we use the family $\{\Gamma_i\}_{i=1}^M$, where each layer of Golgi cisternae is represented by a surface $\Gamma_i$ for some $i$. The followings are four important functionals related to the thesis.

\subsection{Willmore Energy}

 First, based on a well-accepted fact suggested in some previous works \cite{Helfrich,Helfrich2,Redblood,elastic,bending} that the biomembranes are closely related to elastic energy, we consider the elastic energy as the dominated energy to be minimized in our model. We include this energy because $\Gamma$ is also a surface of biomembranes. The elastic energy, also named bending energy, to the lowest order, take the form
     $$E(\Gamma) = \frac{K_b}{2}\int_\Gamma h^2\,dS + \frac{K_G}{2}\int_\Gamma G\,dS,$$
     where the integral $\int_\Gamma\,dS$ is taken over the surface $\Gamma$, $K_b$ is the bending rigidity with respect to mean curvature and $K_G$ is the bending rigidity with respect to Gaussian curvature. $h$ and $G$ are the mean and Gaussian curvature respectively defined as $h:=\frac{1}{2}(C_1+C_2),\ G:=C_1C_2$. $C_1$ and $C_2$ are the principle curvatures. Assume that $\Gamma$ is a closed surface. Gauss-Bonnet Theorem \cite[Ch.\ 8]{bonnet} tells,
     $$\int_\Gamma G\,dS = 2\pi \chi(\Gamma),$$
     where $\chi(\Gamma)$ is the Euler characteristic of $\Gamma$, which is topological invariant. Hence, we only need to consider the Willmore energy \cite{Willmore}
     $$W(\Gamma):=\int_\Gamma h^2\, dS,$$ when $\Gamma$ is a closed surface without topological changes.    
     In summary, to find the optimal shape of Golgi cisternae by minimizing the Willmore energy is the first main point in our model.

\subsection{Area Constraint}   
Second, in many cell models, the surface area of a cell is set to be fixed \cite{Helfrich,Helfrich2}. Agreed by Prof. Kang (Life Science Department), we assume the surface area of each trans-Golgi cisternae $\Gamma$ is conserved. This is based on the fact that the cisternal assembly is completed for tran-cisternae, so the number of molecules of the membrane surface is assumed to be invariant. To enforce this constraint, we consider the functional $A(\Gamma)$ defined by
 $$A(\Gamma):=\int_\Gamma 1\,d\Gamma,$$
 which indicates the area of the surface $\Gamma$. Utilizing this functional, we impose the area constraint into the model. 
 
 \subsection{Barrier Functional}
Third, we want to include the inter-cisternal elements of Golgi stack, which is biologically relevant to the membrane stacking, into our model. The inter-cisternal elements may limit the expansion of Golgi cisternae in the lateral direction. It serves as the obstacles/barriers placed above and below each Golgi cisternae. To model these restraining factors, we construct a shape functional $H(\Gamma)$ by the integration of an indicator function $1_{\bf B}:\Gamma\rightarrow \{0,1\}$, where ${\bf B}$ is a subset of $\mathbb{R}^n$ indicating the region of the obstacles. We define the functional $H(\Gamma)$ as
$$H(\Gamma):=\int_{\Gamma}1_{\bf B}({\bf x})\,d\Gamma,$$
where ${\bf x}$ is the identity on $\Gamma$. We want to include this functional to our model so that the surface $\Gamma$ gets hard to cross the region indicated by the set ${\bf B}$. Detailed explanation of the usage of this functional is given in Section 3.2.

\subsection{Distance Functional - for multiple vesicles case}
Fourth, Golgi stack consists of multiple layers of Golgi cisternaes. When multiple vesicles $\{\Gamma_i\}_{i=1}^n$ are placed in one model, the interaction among the vesicles should be considered. In other words, the vesicles $\Gamma_i$ should not cross each other, and even the repulsion between the vesicles should be taken into consideration, because the lipid bilayer structure of the bio-membranes could cause repulsion when the vesicles approach each other. In this case, we consider the distance among the cisternaes. We also use it to construct a functional $D(\Gamma)$. The detailed formulas and explanations are demonstrated in Section 3.3. In summary, $D(\Gamma)$ is built to control the multiple vesicles case.

\section{Outline}
 \begin{itemize}
\item Chapter 2: The preliminary definitions and lemmas related to the thesis are stated. Most of them are in the field of Differential Geometry. Though we did not work on the theories of Geometry, the theoretical results worked by the predecessors are important for us to construct the numerical algorithm.
\item Chapter 3: The detailed constructions of three models are illustrated. Model 2 and Model 3 are main contributions of this thesis. The three models are introduced in Section 3.1, 3.2 and 3.3 respectively. In each of the section, the motivation of building the model, the idea of the model and the detailed problem setting of the model are stated.
\item Chapter 4: In the first part of this chapter, the detailed explanation of the time discretization and the space discretization (the construction of the mesh) is given. In the second part, the linearization of some nonlinear functions are illustrated, and then the fully discretized weak formulas for the discrete problems are written for each model introduced in Chapter 3. Finally, the full algorithm is given.
\item Chapter 5: The first part of this section gives some numerical examples without biological meanings, only to demonstrate the models and to see the conservation of $A(\Gamma)$ and the decrease of the energy. The second part gives some simulations of Golgi stack, including single cisternae case and multi-cisternae case. Some possible explanations to the morphological properties of the Golgi are given, according to the comparison of the numerical results and the observed Golgi.
\end{itemize}

\section{Previous Work and Our Contributions}

The mathematical study of the shape of biomembranes is introduced by the Helfrich model in 1970s \cite{Helfrich,Helfrich2}, which is a model aiming to study the equilibrium shape of vesicles dominated by elastic energy (or called bending energy). After that, further works on this topic have been done \cite{elastic,bending,Barrett,Dziuk2008,Nochetto2010,Qiang05,Mayer2002,Droske2004,Bobenko2005}, including theoretical analysis and numerical implementation. Many of them applied FEM \cite{Dziuk2008,Nochetto2010}. Besides, other methods were also studied, for instance, finite difference method \cite{Mayer2002}, level set method \cite{Droske2004} and discrete Willmore flow method \cite{Bobenko2005}. Previous works give us inspirations. For example, the Lagrange Multiplier Method is commonly used for the area constraint. Also, previous studies provided excellent formulas for the shape derivative of Willmore energy.

In the thesis, I introduce three models in Section 3.1, 3.2 and 3.3 respectively, and then these models are implemented by FEM with the algorithm stated in Chapter 4. In Section 3.1, it brings out a model dominated by the Willmore energy with conservation of the surface area. This is not a new model developed and solved by us. The works mentioned above \cite{Barrett,Dziuk2008,Nochetto2010,Qiang05} have studied this model. Some of them only includes the area constraint \cite{Barrett,Dziuk2008}. Some also take the volume into account \cite{Nochetto2010,Qiang05}. If one includes both the area and volume constraints, the numerical results could explain the concave shape of the blood cells in a numerical way \cite{Helfrich,Nochetto2010}. 

Besides, Model 2 (Section 3.2) and Model 3 (Section 3.3) are presented in Chapter 3. These two models are constructed by considering the properties of the Golgi stack and these are the main contribution of our work. Model 2 is an extension of Model 1 by adding some obstacles into the problem. This is inspired by the existence of some biological elements which locate above and below each Golgi cisternae and may confine the vertical extension of the cisternae. Besides, as those confining elements could be moving, we also demonstrate the examples of the moving obstacles. Model 3 is designed for the Golgi stack of multiple cisternaes as a whole. These cisternaes could not cross each other, and even the repulsion of their surface should also be considered. Hence, our method could handle the problem with multiple objects. In summary, Model 2 and 3 are newly formed by us, and we also implement the models by FEM. The numerical results are shown in Chapter 5.

Lastly, by comparing our numerical simulation of the Golgi trans-cisternae to the observed Golgi, some possible explanations are drawn on the morphological performance of the trans-cisternae.

\chapter{Mathematical Background}

The theoretical background of this research is introduced in this chapter. It contains many theories in the field of Differential Geometry. Without the theoretical results of Geometry worked by the predecessors, the numerical methods will be hard to implement. In this chapter, $\Gamma$ denotes a hyperspace embedded in $\mathbb{R}^k$, which is piecewise smooth. The mathematical concepts and theories defined on $\Gamma$, which are relevant to this thesis research, are presented. 

\section{Tangential Calculus}
\begin{definition} {\bf (Tangential Gradient)}\\
The tangential gradient of a function $\phi \in C^1(\Gamma; \mathbb{R})$ is defined as
$$\nabla_{\Gamma}\phi := \nabla\tilde{\phi}|_{\Gamma}- \frac{\partial \tilde{\phi}}{\partial n}n,$$
where $\tilde{\phi}$ is a smooth extension to $\phi$ such that $\tilde{\phi}|_{\Gamma} = \phi$.
\begin{remark}
The value of $\nabla_\Gamma\phi$ is independent on the extension function $\tilde{\phi}$ chosen.	
\end{remark}

\end{definition}

With the above definition, we can write down the corresponding tangential Jacobian matrix for a vector function $\boldsymbol{\phi}\in C^1(\Gamma;\mathbb{R}^d)$:
$$[\mathcal{D}_\Gamma (\boldsymbol{\phi})]_{ij} = [\nabla_\Gamma\boldsymbol{\phi}_i]_j.$$

\begin{definition} {\bf (Tangential Divergence)}\\
The tangential divergence of a function $\boldsymbol{\phi}\in C^1(\Gamma;\mathbb{R}^d)$ is defined as
$$\text{div}_\Gamma (\boldsymbol{\phi}):= \text{tr}(\mathcal{D}_\Gamma (\boldsymbol{\phi})),$$
where $\text{tr}(\cdot)$ represents the trace of the matrix.
	
\end{definition}

\begin{definition} {\bf (Tangential Laplace-Beltrami operator)}\\
The tangential Laplace of a function $\phi\in C^2(\Gamma; \mathbb{R})$ is defined as
$$\Delta_\Gamma \phi = \text{div}_\Gamma\,\nabla_\Gamma \phi$$ 
	
\end{definition}

\section{Shape Differential Calculus}
Consider a domain $\Gamma\subset{\mathbb{R}^{k}}$ and a functional defined on $\Gamma$ in the following form
$$J(\Gamma) = \int_\Gamma \phi\, d\Gamma.$$
Consider a family of transformations of $\Gamma$, $\{\Gamma_t\}_{t\in [0,T]}$. Denote the transformation by $\mathcal{T}_t$ such that $\Gamma_t = \mathcal{T}_t(\Gamma)$ and $\Gamma_0 = \Gamma$. We assume that $\mathcal{T}_t$ is a diffeomorphism from $\Gamma$ to $\Gamma_t$ (see \cite[Ch.\ 5]{Shapes} for the definition of diffeomorphism). Denote by $\Gamma_{all}$ the domain containing $\overline{\Gamma_t}$ for all $t\in[0,T]$. Let $\boldsymbol{v}$ be the velocity field associated with the transformation $\mathcal{T}_t$, then we can arrive at the following definition.

\begin{definition} {\bf (shape derivative)} Suppose $J$ is {\it shape differentiable} at $\Gamma$ (see \cite[Ch.\ 5]{Shapes} for the definition). The shape derivative $dJ$ at $\Gamma$ according to the direction $\boldsymbol{v}$ is defined as
$$dJ(\Gamma; \boldsymbol{v}) = \lim_{t\rightarrow 0^+} \dfrac{J(\Gamma_t)-J(\Gamma)}{t},$$
where the transformation $\Gamma_t = \mathcal{T}_t(\Gamma)$ is associated with $\boldsymbol{v}$.
	
\end{definition}

Here are two examples of the formulas of shape derivatives of the functional $W(\Gamma)=\int_\Gamma h^2\,d\Gamma$, where $h$ is the mean curvature. And the detailed derivations of these formulas are provided in \cite{Nochetto2010} and \cite{Dziuk2008}. Note that ${\bf h}:=h\boldsymbol{\nu}$.
$$\begin{array}{ll}
	dW(\Gamma;\boldsymbol{v}) = \int_\Gamma \nabla_\Gamma{\boldsymbol{v}}\colon \nabla_\Gamma{\bf h} &- \int_\Gamma \nabla_\Gamma\boldsymbol{v} (\nabla_\Gamma \text{Id}+\nabla_\Gamma\text{Id}^T)\colon\nabla_\Gamma {\bf h}\\
\ &+\frac{1}{2}\int_\Gamma\text{div}_\Gamma{\bf h}\ \text{div}_\Gamma\boldsymbol{v}\end{array}
$$

$$	\begin{array}{ll}
dW(\Gamma;\boldsymbol{v}) = &-\int_\Gamma \nabla_\Gamma\boldsymbol{v}\colon \nabla_\Gamma{\bf h} + \int_\Gamma D(\boldsymbol{v})\nabla_\Gamma \text{Id}\colon\nabla_\Gamma {\bf h}\\
\ &-\int_\Gamma \text{div}_\Gamma{\bf h}\ \text{div}_\Gamma \boldsymbol{v}-\frac{1}{2}\int_\Gamma|h|^2\text{div}_\Gamma\boldsymbol{v}\end{array},
$$
where $D(\phi)_{ij}=[\nabla_\Gamma \phi_j]_i +[\nabla_\Gamma \phi_i]_j$.

Recall the functional
$$J(\Gamma) = \int_\Gamma \phi\,d\Gamma.$$
Now suppose the function value of $\phi$ only depends on ${\bf x}$. In other words, the function value can be represented by $\phi({\bf x})$. According to \cite[Ch.\ 2]{introshape}, the following lemma is applicable.
\begin{lemma}
\label{lemma:dJ}
Suppose $\phi \in W_1^2(\mathbb{R}^k)$ is independent on the geometry $\Gamma \subset{\mathbb{R}^k}$ and $\Gamma$ is of class $C^2$. Then in the direction $\boldsymbol{v}$,
$$dJ(\Gamma;\boldsymbol{v}) = \int_\Gamma \nabla\phi \cdot \boldsymbol{v}+\phi(\text{\rm div}_\Gamma \boldsymbol{v})\,d\Gamma.$$ 	
\end{lemma}

For example, this lemma can be applied to the functional $A(\Gamma)=\int_\Gamma 1\,d\Gamma$. Here $\phi = 1$ is clearly independent of $\Gamma$. One can simply apply {\bf Lemma 2.1.6} to obtain the formula
$$dA(\Gamma;\boldsymbol{v}) = \int_\Gamma \text{div}_\Gamma \boldsymbol{v}\,d\Gamma.$$ 
Since $\int_\Gamma \text{div}_\Gamma \boldsymbol{v}\,d\Gamma = \int_\Gamma \boldsymbol{v}\cdot {\bf h} +\int_{\partial \Gamma} \boldsymbol{v}\cdot \boldsymbol{\nu}_s$ (see \cite{ellipPDE}), where $\boldsymbol{\nu}_s$ represents the conormal vector, and the surface $\Gamma$ we consider do not have boundary. Hence, we arrive at the following equation that we widely use in our work:
$$dA(\Gamma;\boldsymbol{v}) = \int_\Gamma \boldsymbol{v}\cdot {\bf h}\,d\Gamma.$$

% shape derivatives

% material derivative

% derive dA

% lemma for calculating shape derivatives

% semi-implicit Euler Method

\chapter{Model and Problem Setting}

\section{Willmore energy with Area Constraint}
With the functionals introduced above, we can define a functional $J(\Gamma)$ dependent on the specific problem we are working on, and this functional $J(\Gamma)$ functions as the dominant energy that we want to minimize on the surface $\Gamma$. To track the motion of $\Gamma$ dominated by the energy $J(\Gamma)$, a typical way is to define a geometric evolution equation using the shape derivative $dJ$. Hence, the main idea of our numerical method to solve the shape evolution $\Gamma(t)$ is to find the velocity $\bf v$, which satisfies the equation
  \begin{equation}\label{eq:evolution}\langle {\bf v},{\bf \phi}\rangle = -dJ(\Gamma; {\bf \phi}),\quad \forall {\phi}\in Hil(\Gamma), \end{equation}
  where $(Hil(\Gamma), \langle\cdot,\cdot\rangle)$ is a Hilbert space of functions defined on $\Gamma$.

Consequently, we can define different models by constructing different energies $J(\Gamma)$. In the following sections, I formula three sets of problems corresponding to three models, which can be numerically solved by the discrete schemes described in Chapter 4. In the following three sub-sections, each model is illustrated with detailed equations.

\subsection{The Model 1 and the Functional $J_1(\Gamma)$}
The model 1 is a basic model. It considers only one vesicle $\Gamma$. The main idea of Model 1 is to minimize the Willmore energy under the condition that the surface area of $\Gamma$ is fixed. Hence, it does with the functionals $W(\Gamma)=\frac{1}{2}\int_\Gamma h^2\,d\Gamma$ and $A(\Gamma)=\int_\Gamma 1\,d\Gamma$. They are briefly discussed in Section 2.1.1 and Section 2.1.2. The biological reason of the conservation of the area is also given in Section 2.1.2. 
\begin{equation}\begin{array}{cc} \displaystyle
	\text{minimize} & W(\Gamma)=\frac{1}{2}\int_\Gamma h^2\,d\Gamma,\\
	\text{subject to} & A(\Gamma) = A(\Gamma_0),
	\end{array}
\end{equation}
where $\Gamma_0$ is a given initial shape. The confinement $A(\Gamma)=A(\Gamma_0)$ is imposed to the problem by using a multiplier $\lambda\in \mathbb{R}$. Then $J_1(\Gamma)$ is formulated as
\begin{equation}J_1(\lambda,\Gamma) := W(\Gamma)+\lambda (A(\Gamma)-A(\Gamma_0))=\frac{1}{2}\int_\Gamma h^2\,d\Gamma+\lambda\left(\int_\Gamma \,d\Gamma-\int_{\Gamma_0}\,d\Gamma_0\right).\end{equation}
We aim to find the optimal $\lambda$ and $\Gamma$. 

\subsection{Problem Setting of Model 1}
First, we define
$$G_T = \bigcup_{t\in [0,T]} \Gamma(t)\times \{t\}$$
and ${\bf u}: G_T \rightarrow \mathbb{R}^n$ by 
\begin{equation}\label{eq:u}{\bf u}(x,t)= x\end{equation}
for all $x \in \Gamma(t)$ and $t\in [0,T]$. Hence, we can regard ${\bf u}(\cdot,t)= \text{Id}_{\Gamma(t)}$. The goal is now to minimize the functional $J_1(\lambda, \Gamma)$. Hence, it is a typical way to form an evolution equation in the form of \eqref{eq:evolution}. Consequently, we have the following problem setting.

\begin{p} {\bf (Willmore Flow with Area Constraint, Weak Form)}\\
For a given initial shape $\Gamma_0 = \Gamma(0) \subset \mathbb{R}^n$, find the multiplier $\lambda:[0,T]\rightarrow \mathbb{R}$ and the function ${\bf u}: G_T\rightarrow \mathbb{R}^n$ according to the family of surface $\{\Gamma(t)\}_{t\in[0,T]}$, such that on the time interval $t\in[0,T]$,
   \begin{equation}\label{eq:problem1}\int_{\Gamma(t)}{\dot{\bf u}}\cdot {\bf \phi} = -dJ_1(\lambda(t),\Gamma(t);\phi)= -dW(\Gamma(t);\phi) - \lambda(t) dA(\Gamma(t);\phi),\end{equation}
   \begin{equation}A(\Gamma(t))=A(\Gamma_0),\end{equation}
for all test function $\phi$. The function space of $\phi$ will be later discussed.
\end{p}
The shape derivatives of the above functionals can be computed using the following equations
\begin{equation}\label{dW1}\begin{array}{ll}
dW(\Gamma;\phi) = \int_\Gamma \nabla_\Gamma{\bf \phi}\colon \nabla_\Gamma{\bf h} &- \int_\Gamma \nabla_\Gamma{\bf \phi} (\nabla_\Gamma \text{Id}+\nabla_\Gamma\text{Id}^T)\colon\nabla_\Gamma {\bf h}\\
\ &+\frac{1}{2}\int_\Gamma\text{div}_\Gamma{\bf h}\ \text{div}_\Gamma{\bf \phi},\end{array}
\end{equation}
or 
\begin{equation}\label{dW2}\begin{array}{ll}
dW(\Gamma;\phi) = &-\int_\Gamma \nabla_\Gamma{\bf \phi}\colon \nabla_\Gamma{\bf h} + \int_\Gamma D(\phi)\nabla_\Gamma \text{Id}\colon\nabla_\Gamma {\bf h}\\
\ &-\int_\Gamma \text{div}_\Gamma{\bf h}\ \text{div}_\Gamma {\bf \phi}-\frac{1}{2}\int_\Gamma|h|^2\text{div}_\Gamma{\bf \phi},\end{array}
\end{equation}
and
\begin{equation}\label{dA}
	dA(\Gamma;{\bf \phi})=\int_\Gamma {\bf h}\cdot {\bf \phi}.
\end{equation}
The details of the first two equations can be found in \cite{Nochetto2010} and \cite{Dziuk2008}. The third equation is simpler, so its derivation is given in the Preliminary Section.

Some of the notations used in the above equations are explained as follows:
\begin{itemize}
\item ${\bf h}:=h\boldsymbol{\nu}$ is the vector form of mean curvature $h$ on the direction of the outer unit normal vector $\boldsymbol{\nu}$. 
\item $\text{Id} = {\bf u}(\cdot,t)$ is the identity on $\Gamma(t)$.
\item $D$ is a symmetric tensor defined by $D(\phi)_{ij}=(\nabla_\Gamma)_i\phi^j +(\nabla_\Gamma)_j\phi^i$.
\end{itemize}

The equation \eqref{dW1} and \eqref{dA} are implemented for the Model 1.

This Model 1 is not a model firstly produced and solved by us. Instead, many works \cite{Barrett,Dziuk2008,Nochetto2010,Qiang05} have been done on the study of this model, of which many also consider the volume constraint \cite{Nochetto2010,Qiang05}. In summary, Model 1 is a fundamental model, which has been studied for many years. However, it catches an important property of the bio-membranes, that the membranes is closely relevant to the elastic energy. Also, it imposes the idea that the number of molecules of the membranes is fixed such that the surface area of the membrane is conserved. If one includes one more constraint, the volume constraint, to Model 1, the numerical results can explain the concave shape of the blood cells mathematically \cite{Helfrich,Nochetto2010}. However, since our work is triggered by the shape evolution of Golgi stack, instead of considering the volume constraints, we consider other special properties of the Golgi stack and build our own models. The coming two sections present two models of Golgi stacks. Model 2 is designed for a singer layer of Golgi cisternae and Model 3 is designed for the Golgi stack of multiple cisternaes as a whole.

\section{Willmore Energy with Area Constraint and Obstacles}
\subsection{The Model 2 and the Functional $J_2(\Gamma)$}
Model 2 is an extension to Model 1 by considering the existence of some obstacles/barriers in the problem. This is motivated by a property of Golgi cisternae, which may be confined by the intercisternal elements above and below each Golgi cisternae layer. (See Section 2.1.3). In Model 2, the obstacles are considered and estimated by the functional $$H(\Gamma)=\int_\Gamma 1_{\bf B}(x)\,d\Gamma$$ defined in Section 2.1.3, where ${\bf B}$ represents the region of the obstacles/barriers. One can easily observe that $H(\Gamma)$ will take nonzero values only if $\Gamma\cap {\bf B}\neq \emptyset$. Hence, if we add this functional $H(\Gamma)$ to the energy that we want to minimize, it will be hard for the shape $\Gamma$ to touch and cross the region ${\bf B}$ so as to avoid the increase of the total amount of the energy.

To conclude, in Model 2, we want to do the same optimization as in Model 1, but also to include some obstacles indicated by the set ${\bf B}\subset\mathbb{R}^n$. Normally, this model can be explained by the following problem:
\begin{equation}\begin{array}{cc} \displaystyle
	\text{minimize} & W(\Gamma)+\alpha H(\Gamma),\\
	\text{subject to} & A(\Gamma) = A(\Gamma_0).
	\end{array}
\end{equation}
The resulted functional $J_2$ can then be formulated as
\begin{equation}
\begin{array}{ll}
J_2(\lambda,\Gamma) & := W(\Gamma)+\alpha H(\Gamma)+\lambda (A(\Gamma)-A(\Gamma_0))\\
\ &=\frac{1}{2}\int_\Gamma h^2\,d\Gamma+\alpha\int_\Gamma 1_{\bf B}(x)\,d\Gamma+\lambda\left(\int_\Gamma \,d\Gamma-\int_{\Gamma_0}\,d\Gamma_0\right),\end{array}\end{equation}
which is the augmented energy to be minimized in Model 2.

\begin{remark}
The constant $\alpha\in \mathbb{R}$ is a weight of the functional $H(\Gamma)$ to control the impact of $H$. The dominated energy of our model should be the Willmore energy $W(\Gamma)$. $H(\Gamma)$ is only a constraint. Hence, we don't want the functional $H(\Gamma)$ to  dominate the whole energy $J_2(\Gamma)$.	
\end{remark}

\subsection{Problem Setting of Model 2}

Use the same notation as those in Section 3.1.2, we can form the following weak problem to minimize $J_2(\lambda, \Gamma)$.

\begin{p} {\bf (Willmore Flow with Area Constraint and Obstacles, Weak Form)}\\
Suppose that $\alpha\in \mathbb{R}$ is fixed as a weight coefficient. Now given $\Gamma_0 = \Gamma(0) \subset \mathbb{R}^n$, find $\lambda:[0,T]\rightarrow \mathbb{R}$ and the function ${\bf u}: G_T\rightarrow \mathbb{R}^n$ such that on the time interval $t\in[0,T]$,
   \begin{equation}
   \begin{array}{lcl}\int_{\Gamma(t)}{\dot{\bf u}}\cdot {\bf \phi} &=& -dJ_2(\lambda(t),\Gamma(t);\phi)\\
   \ &=&-dW(\Gamma(t);\phi) -\alpha dH(\Gamma(t);\phi) - \lambda(t) dA(\Gamma(t);\phi),\end{array}\end{equation}
   \begin{equation}A(\Gamma(t))=A(\Gamma_0),\end{equation}
for all test function $\phi$.
\end{p}

The calculation of $dH$ is discussed as follows.

Since the indicator function $1_{\bf B}:\Gamma \rightarrow\{0,1\}$ is discontinuous, the implementation of it using FEM is not applicable. In the finite element method, we choose a smooth version of the indicator. More precisely, assume that ${\bf B}$ represents a very regular shape. Then, $1_{\bf B}$ can be written as a composition of the heaviside function $$heav(x):=1_{[0,\infty]}.$$
\begin{remark}
For example, if ${\bf B}=[a,b]\times ([m,\infty]\cup[-\infty,-m])\subset \mathbb{R}^2$, then it can be written as a composition of $heav(x)$ by
    $$1_{\bf B} = heav(x-a)heav(-x+b)\cdot \left(heav(y-m)+heav(-y+m)\right).$$
    \end{remark}
    Though $heav(x)$ is still a discontinuous step function, we use the smooth approximation $$heav(x)\approx\dfrac{1}{1+e^{-2kx}}$$
    of it. Now a smooth version of $1_{\bf B}$ is obtained.
    
    With the smoothness of $1_{\bf B}$, we can now derive the formula for $dH$. By lemma \ref{lemma:dJ}, the formula of $dH(\Gamma)$ is obtained:
    \begin{equation}\label{dH}dH(\Gamma;\phi) = \int_\Gamma \nabla1_{\bf B} \cdot \phi +
    \int_\Gamma 1_{\bf B}\,\text{div}_\Gamma \phi.\end{equation}
    
    \begin{remark}
    	Usually, in the math models, ${\bf B}$ is some fixed obstacles and independent of time. However, inspired by the hypothesis mentioned by Prof. Kang, that the intercisternal elements (such as Golgi matrix, which works as the constraining factor in our math model on Golgi stacks) maintain the same distance with the membrane since some of its components are embedded in the membrane. Hence, it could be more realistic to keep the distance between the membrane and the barriers when we mimic the growing process of Golgi cisternae. Nevertheless, since the shape of the cisternae $\Gamma(t)$ evolves with time, the position of the barriers also need changesx to keep the distance. Based on this, we also make some numerical experiments for the moving obstacles ${\bf B}(t)$. These examples can be found in Chapter 5. 
    \end{remark}

\section{Multiple Vesicles Case}
\subsection{The Distance Functional}
As succinctly introduced in Section 2.1.4, the Golgi stack consist of multiple layers of Golgi cisternaes, denoted by $\Gamma_1, \Gamma_2, \cdots, \Gamma_n$ in our math model. The family $\{\Gamma_i\}_{i=1}^n$ is aimed to be modeled on the whole, at the same time the vesicles $\Gamma_i$ should not cross or even should repulse from each other. We applied the following function with Euclidean norm $\|\cdot\|$ to measure the distance between the vesicles:
 \begin{equation}
 d(x) = \min_{\substack{y\in V_j\\j\neq i}}\|x-y\|^2,\ \ \forall x\in V_i,	
 \end{equation}
where $V_i$ denotes the $i$-th vesicle. Equivalently, $d(x)$ can be defined as
\begin{equation}
\forall x\in V_i,\ d(x) = \|x-y(x)\|^2,\ y(x) = \arg\min_{\substack{y\in V_j\\j\neq i}}\|x-y\|^2.	
\end{equation}
Using this measurement of distance, we define the following functional
\begin{equation}
D(\Gamma) = \int_\Gamma	\frac{1}{d(x)}\,d\Gamma,
\end{equation}
where $x$ is the identity on $\Gamma$. Note that $\Gamma$ only represents the surface of one single vesicle $V_i$.
\subsection{The Model 3 and the Functional $J_3(\Gamma)$}
Well-prepared with the above functional $D(\Gamma)$, we can now extend the single vesicle case - Model 2 to the multiple case - Model 3. Let $\Gamma = \partial V_i$ for some $i$. For each $i$, we do the following problem:
\begin{equation}\begin{array}{cc} \displaystyle
	\text{minimize} & W(\Gamma)+\alpha H(\Gamma)+\beta D(\Gamma),\\
	\text{subject to} & A(\Gamma) = A(\Gamma_0),
	\end{array}
\end{equation}
where $\alpha, \beta\in \mathbb{R}$ are the weight coefficients for $H(\Gamma)$ and $D(\Gamma)$ respectively. To translate this optimization problem: in fact, we are minimizing the Willmore energy $W(\Gamma)$, under the condition that, firstly, $\Gamma$ is hard to tough the region $B$ with a weight $\alpha$ and secondly, $\Gamma$ is hard to get very closed to other vesicles with a weight $\beta$.

The associated functional $J_3$ corresponding to this problem is given by
\begin{equation}
\begin{array}{ll}
J_3(\lambda,\Gamma) & := W(\Gamma)+\alpha H(\Gamma)+\beta D(\Gamma)+\lambda (A(\Gamma)-A(\Gamma_0))\\
\ &=\frac{1}{2}\int_\Gamma h^2 +\alpha\int_\Gamma 1_{\bf B}(x) +\beta\int_\Gamma \dfrac{1}{d(x)} +\lambda\left(\int_\Gamma 1 -\int_{\Gamma_0} 1 \right),\end{array}\end{equation}
which is the augmented energy to be minimized in Model 3.
\begin{remark}
The Model 3 is also newly formed by us. It is motivated by the component of the multiple layers of Golgi. It could be extended to other problems which include multiple objects.
\end{remark}

\subsection{Problem Setting of Model 3}
Similar to what we do for the above two models, we can now form the following weak problem to minimize the functional $J_3(\lambda, \Gamma)$.
\begin{p} {\bf (Willmore Flow with Area Constraint and Obstacles - applied to Multiple Vesicles Case, Weak Form)}\\
Suppose that $\alpha,\beta\in \mathbb{R}$ are fixed as a weight coefficient for the functional $H$ and $D$ respectively. Given an initial shape $\Gamma_0 = \Gamma(0) \subset \mathbb{R}^n$, find the multiplier $\lambda:[0,T]\rightarrow \mathbb{R}$ and the function ${\bf u}: G_T\rightarrow \mathbb{R}^n$ according to the family of surfaces $\{\Gamma(t)\}_{t\in[0,T]}$, such that on the time interval $t\in[0,T]$,
   \begin{equation}
   \begin{array}{lcl}\int_{\Gamma(t)}{\dot{\bf u}}\cdot {\bf \phi} &=& -dJ_3(\lambda(t),\Gamma(t);\phi)\\
   \ &=&-dW(\Gamma(t);\phi) -\alpha dH(\Gamma(t);\phi) -\beta dD(\Gamma(t);\phi)- \lambda(t) dA(\Gamma(t);\phi),\end{array}\end{equation}
   \begin{equation}A(\Gamma(t))=A(\Gamma_0),\end{equation}
for all test function $\phi$.
\end{p}

The formulas of the shape derivatives of $W(\Gamma)$, $H(\Gamma)$ and $A(\Gamma)$ are clearly explained in the previous sections, so I only demonstrate the calculation of $dD(\Gamma;\phi)$ in this section. Since our model is implemented by FEM at the end, I estimate the shape derivative of $D(\Gamma)$ in a tricky way. First, for any point $x^*\in \Gamma_i$, consider an open ball $S = B(x^*;r)$ centered at $x^*$ with radius $r$ and define a function $d^*:S\rightarrow \mathbb{R}^+$ by
$$d^*(x)=\|x-y^*\|^2,\ \text{where}\ y^* = \arg\min_{\substack{y\in V_j\\j\neq i}}\|x^*-y\|.$$
Since $1/d^*$ is independent on the geometry, we can apply Lemma \ref{lemma:dJ} and then obtain the shape derivative of $D^*(S\cap \Gamma):=\int_{s\cap\Gamma}f^*(x)\,d\Gamma$.
\begin{equation}\label{dD}
dD^*(S\cap\Gamma;\phi) = -2\int_{S\cap\Gamma}\dfrac{1}{d^*(x)^2}(x-y^*)\cdot \phi +\int_{S\cap \Gamma}\dfrac{1}{d^*(x)}\text{dir}_\Gamma\phi .	
\end{equation}
Applying this formula, then we can approximate the shape derivative of $D(\Gamma)$ by piecewise implementation in FEM.

\chapter{Numerical Schemes}

\section{Time Discretization and Equation Split}
\subsection{Time Discretization}
The model is implemented on the time-interval $[0,T]$, though the final time $T$ can be chosen dependent on the stopping criteria set in the algorithm. Let $$\bigcup_{i=0}^{N-1}[t_i,t_{i+1}] = [0,T],$$ where $t_0 = 0,\ t_N=T$, be a partition of the interval $[0,T]$. 

\begin{remark} {\bf About Time Adaptivity:} If one wants to obtain a more effective algorithm to reach the optimization shape of $\Gamma$ faster, it is more reasonable to make the time step $$\tau_n := t_{n+1}-t_n$$ adaptive to the mesh size, because this FEM is using a moving-mesh. For me, I simply choose a comparatively small time step $\tau = \tau_n$, which is fixed, for convenience. However, there is actually a goodness of using a small time step $\tau$, because we apply the linear approximations (see Section 4.3.1) on some shape functionals in our method and a small time step is beneficial to the linear approximations.
\end{remark}

Now the solution we want to find is the family $\{\Gamma(t_i)\}_{i=1}^N$ when $\Gamma(t_0)$ is given. Recall the function ${\bf u}: G_T \rightarrow \mathbb{R}^n$ defined by equation \eqref{eq:u}. From now on, we denote the numerical solution by $\Gamma^n$ at each time $t_n$, which is viewed as the image of ${\bf U}^n(\cdot)$
 \begin{equation}\label{Gamman}
 \Gamma^n:= \{{\bf U}^n(x);\ \forall x\in \Gamma^{n-1} \}
 \end{equation}
 with given approximation $\Gamma_0=\Gamma^0\approx\Gamma(t_0)$.
\subsection{Split}
From equation \eqref{eq:problem1} and \eqref{dW1}, one can see that, if we solve ${\bf u}$ directly, the order of the differential equation is high. To reduce this order, many previous work chose to split the formula \cite{Nochetto2010}. In our work, we solve a pair of unknown $({\bf V}^{n+1},{\bf H}^{n+1})$ first, and then update ${\bf U}^{n+1}$. ${\bf V}^{n+1}$ is defined as
\begin{equation}\label{def:V}{\bf V}^{n+1}(x) = \frac{1}{\tau}\left({\bf U}^{n+1}(x)-x\right),\ \forall x\in \Gamma^n.\end{equation}
By applying the equation $-\Delta_\Gamma x = {\bf h}$ \cite[Page 390]{ellipPDE} and the above equation \eqref{def:V}, the following identity for ${\bf H}^{n+1}$ is obtained
\begin{equation}
-\Delta_{\Gamma^n} {\bf V}^{n+1} = \dfrac{1}{\tau}\left({\bf H}^{n+1}+\Delta_{\Gamma^n}x\right),
\end{equation}
and hence the following weak formula
\begin{equation}\label{eq:HX}
\tau\int_{\Gamma^n} \nabla_{\Gamma^n}{\bf V}^{n+1}\colon \nabla_{\Gamma^n}{\bf \phi} = \int_{\Gamma^n} {\bf H}^{n+1}\cdot {\bf \phi}-\int_{\Gamma^n}\nabla_{\Gamma^n}x\colon\nabla_{\Gamma^n}{\bf \phi},
\end{equation}
for all test function $\phi$. The function space of the test function will be discussed in the coming section.

In summary, with equation \eqref{eq:HX}, we can solve the problems stated in Chapter 3 by solving the pair of unknowns $({\bf V}^{n+1}, {\bf H}^{n+1})$ first and then update ${\bf U}^{n+1}$ by equation \eqref{def:V}. 
%%% other benefit??

\section{Finite Elements}
The following are some notations and definitions used in this section.
\begin{itemize}
\item $\mathbb{R}^k$ is the space containing $\Gamma$.
\item $\mathbb{R}^{k-1}$ is the parametrization space.
\end{itemize}
{\it Polyhedral Approximation} $(\overline{\Gamma}_h, \overline{\mathcal{T}}_h, \overline{K}, V_{\overline{K}})$
\begin{itemize}
\item $\overline{\Gamma}_h = \bigcup_{\overline{K}\in \overline{\mathcal{T}}_h} \overline{K}$ is a polyhedral approximation of $\Gamma$, where $\overline{\mathcal{T}}_h$ is the triangulation of $\overline{\Gamma}_h$ and the vertices of $\overline{\Gamma}_h$ lie on $\Gamma$.
\item $\overline{K}$ is a $(k-1)$-simplex in $\mathbb{R}^k$ with its $k$ vertices $\{v_i\}_{i=1}^k$.
\item $V_{\overline{K}}=\{v = \sum_{i=1}^k c_i v_i;\ \  \sum_{i=1}^k c_i = 1, c_i\in \{0,\frac{1}{2},1\}\}$ is the vertex set that we use in FEM, which includes the vertices of $\overline{K}$ and the mid-points of each edge of $\overline{K}$. 
\end{itemize}
\begin{figure}[h!]
\centering
\includegraphics[width=0.8\textwidth]{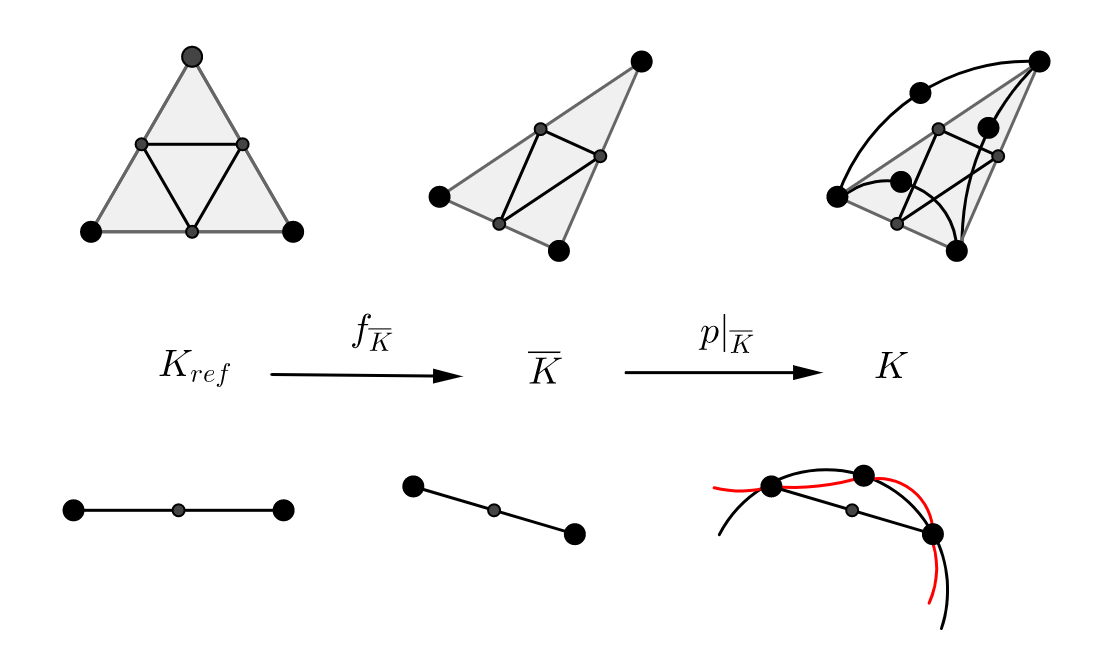}
\caption{Example of the elements in $\mathbb{R}^3$ and $\mathbb{R}^2$; the red curve represent a piece of $\Gamma$.}
\end{figure}

{\it Polynomial Approximation} $(\Gamma_h, K, V_K)$ 
\begin{itemize}	
\item $\Gamma_h$ is the image of a function $p$ defined on $\overline{\Gamma}_h$ such that $p|_{\overline{K}}$ is a polynomial of degree $\leq 2$. 
\item Denote $K = p(\overline{K})$. Then we have $\Gamma_h = \bigcup K$.
\item $V_K = \{p(v);\ v\in V_{\overline{K}}\}$. 
\item[] The function values of $p(v)$ indicate the position of the vertices $V_K$. The values of $p(v)$ is obtained by following the rules:
 \begin{enumerate}
 \item If $v\in \cup_{i=1}^k\{v_i\}$ (i.e. $v$ is a vertex of the $(k-1)$-simplex $\overline{K}$), then $p(v) = v$. Since $v_i$ are on $\Gamma$, then $p(v_i)$ are also on $\Gamma$.
 \item If $v\in V_k \setminus\cup_{i=1}^k\{v_i\}$ (i,e, $v$ is a mid-point of an edge of $\overline{K}$), then $p(v)$ is an orthogonal projection onto $\Gamma$. Hence, $p(v)$ also lies on $\Gamma$.
 \end{enumerate}
 \begin{figure}[h!]
 \centering
 \subfloat[The big black dots are vertices in $V_K^n$.]{
 \includegraphics[width=0.3\textwidth]{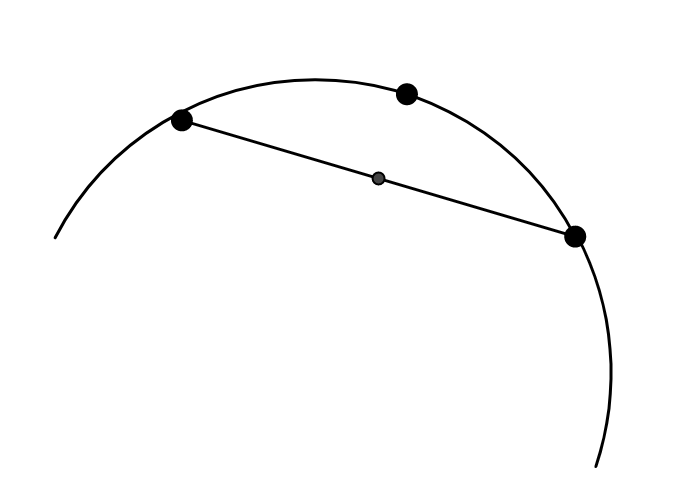}
 }
 \subfloat[The blue dots are the next-step position of the vertices.]{
 \includegraphics[width=0.3\textwidth]{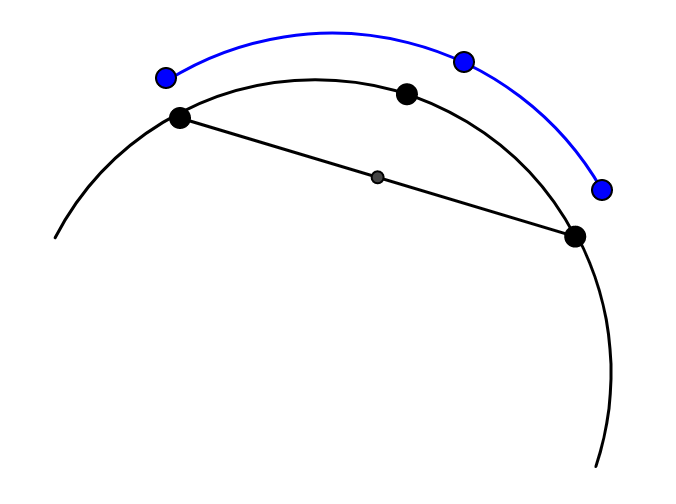}
 }
 \subfloat[The 2nd blue dot is adjusted to be the midpoint orthogonal projection. Then the three red dots form the $V_K^{n+1}$.]{
 \includegraphics[width=0.3\textwidth]{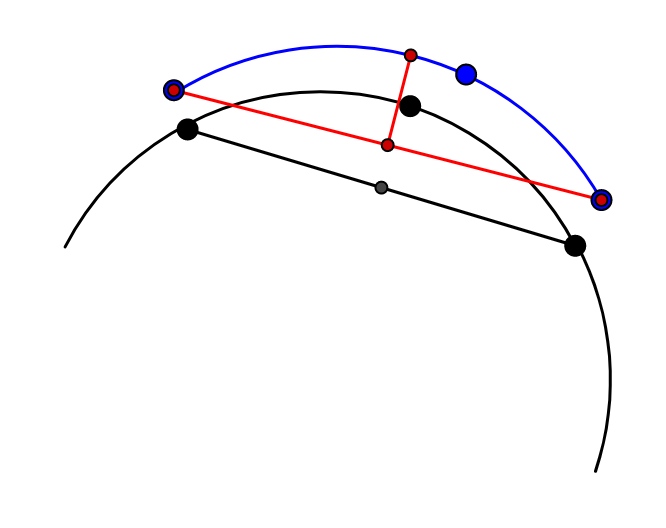}
 }
 \caption{Adjustment of the vertices}
 \end{figure}

To make sure that the set $V_K^{n+1}$ is the orthogonal projection of $V_{\overline{K}}$ on to $\Gamma_h^{n+1}$, we adjust the position of the vertices whenever we obtain a new set $V_K^{n+1}$ from $V_K^{n}$ by solving the numerical problem.
 
 By the above construction, one can conclude that all the points in $V_K$ are on $\Gamma$. Besides, the existence and uniqueness of the polynomial $p|_K$ satisfying the above two rules are proved by \cite[Theorem 2.2.1]{ciarFEM}.
\end{itemize}

{\it Reference Element}
\begin{itemize}
\item $K_{ref}$, a $(k-1)$-simplex in $\mathbb{R}^{k-1}$. We define the standard reference $K_{ref}$ as the convex envelope with vertices $\cup_{i=1}^{k-1}\{{\bf e}_i\}\cup\{{\bf 0}\}$.
 \begin{itemize}\item In $\mathbb{R}^1$, $K_{ref}=[0,1]$.
    \item In $\mathbb{R}^2$, $K_{ref} = \{a_0(0,0)+a_1(0,1)+a_2(1,0):\ \forall a_0, a_1,a_2\ge 0\ \text{s.t.}\ a_0+ a_1+a_2 = 1\}$.\end{itemize}
\item $V_{ref}= \{v = \sum_{i=1}^{k-1} c_i {\bf e}_i;\ \  \sum_{i=1}^k c_i = 1, c_i\in \{0,\frac{1}{2},1\}\}$ represents the vertex set of $K_{ref}$.
\item[] For each $(k-1)$-simplex $\overline{K}$ in $\mathbb{R}^k$ mentioned above, there exists a bijective mapping $f_{\overline{K}}:K_{ref}\rightarrow \overline{K}$ such that $f_{\overline{K}}$ maps $V_{ref}$ to $V_{\overline{K}}$.
\end{itemize}
{\it Finite Element Space}
\begin{itemize}
\item[] The finite element space defined over the set $\Gamma_h$ is defined as \begin{equation}\label{FEMspace}\mathbb{F}(\Gamma_h) := \{u_h\in C^0(\Gamma_h);\ u_h|_K \circ p_{\overline{K}}\circ f_{\overline{K}} \in \mathbb{P}_2(K_{ref}),\ \forall \overline{K}\in \overline{\mathcal{T}_h}\}.\end{equation}
\item[] At each time step $t_{n+1}$, $\Gamma_h^n$ is known. Hence, the finite element function space that we consider at the time step $t_{n+1}$ to find $\Gamma_h^{n+1}$ is $\mathbb{F}(\Gamma_h^n)$.
\end{itemize}

\section{Discrete Problems}
With the time and space discretization discussed in the previous sections, we can now rewrite our problems in the discrete forms, which can be implemented.

In Chapter 3, we formulate the three problems, Problem 3.1.1, Problem 3.2.1 and Problem 3.3.1. In each of them, we aim to solve ${\bf u}(\cdot, t) = \text{Id}_{\Gamma(t)}$ for $t\in[0,T]$. After time discretization, we aim to solve $\Gamma^n :=\{{\bf U}^n(x);\ \forall x\in \Gamma^{n-1}\}$ (equation \eqref{Gamman}) at each time $t_n$. Instead of solving ${\bf U}^{n+1}$ directly from ${\bf U}^{n}$, as discussed in Section 4.1.2, we split the formula by adding the equation \eqref{eq:HX}, so we can solve $({\bf V}^{n+1}, {\bf H}^{n+1})$ first and then update ${\bf U}^{n+1}$ by equation \eqref{def:V}. The discretization (Section 4.2) gives us the space $\mathbb{F}(\Gamma_h^n)$ to solve $({\bf V}^{n+1}, {\bf H}^{n+1})$. Based on these ideas, we can formulate the problems in Chapter 3 into the following discrete forms.

\begin{p}{\bf (Discrete Form of Problem 3.1.1.)}\\
Suppose that $\alpha,\beta\in \mathbb{R}$ are fixed as a weight coefficient for the functional $H$ and $D$ respectively. Given an initial shape $\Gamma_h^0\approx \Gamma(0)\subset \mathbb{R}^k$, find ${\bf V}^{n+1}, {\bf H}^{n+1}\in \mathbb{F}(\Gamma_h^n)$ and $\lambda^{n+1}\in \mathbb{R}$ such that for each $n$, $\forall \Phi\in \mathbb{F}(\Gamma_h^n)$,
\begin{equation}\int_{\Gamma_h^n} {\bf V}^{n+1}\cdot {\Phi}	 = -dW^{n+1}(\Gamma_h^n;\Phi)- \lambda^{n+1}dA^{n+1}(\Gamma_h^n;\Phi),\end{equation}
\begin{equation}
A(\Gamma_h^{n+1}) = A(\Gamma_h^{0}),	
\end{equation}
\begin{equation}
	\tau\int_{\Gamma^n} \nabla_{\Gamma^n}{\bf V}^{n+1}\colon \nabla_{\Gamma^n}{\Phi} = \int_{\Gamma^n} {\bf H}^{n+1}\cdot {\Phi}-\int_{\Gamma^n}\nabla_{\Gamma^n}x\colon\nabla_{\Gamma^n}{\Phi},
\end{equation}
where
\begin{equation}
\begin{array}{ll}
	dW^{n+1}(\Gamma_h^n;\Phi) &= \int_{\Gamma_h^n} \nabla_{\Gamma_h^n}{\Phi}\colon \nabla_{\Gamma_h^n}{\bf H}^{n+1}\\ &- \int_{\Gamma_h^n} \nabla_{\Gamma_h^n}{\Phi} (\nabla_{\Gamma_h^n} \text{Id}+\nabla_{\Gamma_h^n}\text{Id}^T)\colon\nabla_\Gamma {\bf H}^{n+1}\\
\ &+\frac{1}{2}\int_{\Gamma_h^n}\text{div}_{\Gamma_h^n}{\bf H}^{n+1}\ \text{div}_{\Gamma_h^n}{\Phi}\end{array}
\end{equation}
and
\begin{equation}
		dA^{n+1}(\Gamma_h^n;{\Phi})=\int_{\Gamma_h^n} {\bf H}^{n+1}\cdot{\Phi}.
\end{equation}
At each time step, $\Gamma_h^{n+1}$ is updated by 
\begin{equation}
\Gamma_h^{n+1} = \Gamma_h^n +\tau {\bf V}^{n+1}(\Gamma_h^n).	
\end{equation}

\end{p}

\begin{p}
	{\bf (Discrete Form of Problem 3.3.1.)}\\
Given an initial shape $\Gamma_h^0\approx \Gamma(0)\subset \mathbb{R}^k$, find ${\bf V}^{n+1}, {\bf H}^{n+1}\in \mathbb{F}(\Gamma_h^n)$ and $\lambda^{n+1}\in \mathbb{R}$ such that for each $n$, $\forall \Phi\in \mathbb{F}(\Gamma_h^n)$,
\begin{equation}
\begin{array}{ll}
\int_{\Gamma_h^n} {\bf V}^{n+1}\cdot {\Phi}	 = &-dW^{n+1}(\Gamma_h^n;\Phi)- \alpha dH^{n+1}(\Gamma_h^n; \Phi) - \beta dD^{n+1}(\Gamma_h^n;\Phi)\\ &-\lambda^{n+1}dA^{n+1}(\Gamma_h^n;\Phi),
\end{array}
\end{equation}
\begin{equation}
A(\Gamma_h^{n+1}) = A(\Gamma_h^{0}),	
\end{equation}
\begin{equation}
	\tau\int_{\Gamma^n} \nabla_{\Gamma^n}{\bf V}^{n+1}\colon \nabla_{\Gamma^n}{\Phi} = \int_{\Gamma^n} {\bf H}^{n+1}\cdot {\Phi}-\int_{\Gamma^n}\nabla_{\Gamma^n}x\colon\nabla_{\Gamma^n}{\Phi},
\end{equation}
where the formulations of $dH^{n+1}$ and $dD^{n+1}$ are linearized in the next section. They are specially treated and formulated in an implicit form.
\end{p}
Note that the discrete form of {\bf Problem 3.2.1.} is the same as {\bf Problem 4.3.2.} by simply taking $\beta$ to be zero.

\subsection{Linearization of $dH$ and $dD$ in FEM}
\subsubsection{Linearization of $dH$}
 Recall the formula for $dH$ \eqref{dH}:
 $$dH(\Gamma;\phi) = \int_\Gamma \nabla1_{\bf B}({\bf x}) \cdot \phi +
    \int_\Gamma 1_{\bf B}({\bf x})\,\text{dir}_\Gamma \phi.$$
 The function $1_{\bf B}$ with smoothness (by applying exponential functions) is nonlinear. So is $\nabla1_{\bf B}$. In FEM, we use their linear approximation. They are simply made in the standard way:
 $$
 1_{\bf B}({\bf x})\approx \mathcal{L}_{1_{\bf B}}({\bf x})= 1_{\bf B}({\bf x_0})+[\mathcal{D}{1_{\bf B}}({\bf x_0})][({\bf x}-{\bf x_0})],
 $$
$$
 	\nabla1_{\bf B}({\bf x})\approx \mathcal{L}_{\nabla1_{\bf B}}({\bf x})= \nabla1_{\bf B}({\bf x_0})+[\mathcal{D}{(\nabla1_{\bf B})}({\bf x_0})][({\bf x}-{\bf x_0})].
$$
Here $\mathcal{D}$ is the Fréchet derivative operator. It is well-known that the approximation is good for ${\bf x}$ if it is close enough to ${\bf x_0}$, so it is very natural that when ${\bf x}\in \Gamma^{n+1}$ we take the points ${\bf x_0}$ from $\Gamma^n$. When $\tau$ is small, the approximation could be good. Hence, we now replace ${\bf x_0}$ by $x\in \Gamma_h^n$ and ${\bf x}$ by ${\bf U}_h^{n+1}(x)\in \Gamma_h^{n+1}$. More specifically, we write the approximation as: $\forall x\in \Gamma_h^n$,
\begin{equation}
	 	\mathcal{L}_{1_{\bf B}}({\bf U}_h^{n+1}({x}))= 1_{\bf B}({x})+[\mathcal{D}{1_{\bf B}}({x})][({\bf U}_h^{n+1}(x)-{x})],
\end{equation}
\begin{equation}
	 \mathcal{L}_{\nabla1_{\bf B}}({\bf U}_h^{n+1}({x}))= \nabla1_{\bf B}({x})+[\mathcal{D}{(\nabla1_{\bf B})}({x})][({\bf U}_h^{n+1}(x)-{x})].
\end{equation}
 Clearly, in the discretized form, we can replace $\left({\bf U}_h^{n+1}(x)-{x}\right)$ by $\left(\tau {\bf V}^{n+1}\right)$ and then obtain the following semi-implicit formula for $dH^{n+1}(\Gamma_h^n;\Phi)$:
 \begin{equation}
 \begin{array}{ll}
 	dH^{n+1}(\Gamma_h^n;\Phi)=&\int_{\Gamma_h^n}\nabla 1_{\bf B}({\bf x})\cdot \Phi +\tau\int_{\Gamma_h^n}[\mathcal{D}{(\nabla1_{\bf B})}({\bf x})]{\bf V}^{n+1}\cdot \Phi\\
 	 &+\int_{\Gamma_h^n}1_{\bf B}({\bf x})\,\text{dir}_{\Gamma_h^n} \Phi +\tau\int_{\Gamma_h^n}[\mathcal{D}{1_{\bf B}}({\bf x})]{\bf V}^{n+1}\,\text{dir}_{\Gamma_h^n} \Phi.
 	\end{array}
 \end{equation}
\subsubsection{Linearization of $dD$}
 Recall the formula \eqref{dD} for $dD$ on the neighborhood $S\cap \Gamma$ of a point $x^*\in \Gamma$:
 $$dD^*(S\cap\Gamma;\phi) = -2\int_{S\cap\Gamma}\dfrac{1}{d^*(x)^2}(x-y^*)\cdot \phi +\int_{S\cap \Gamma}\dfrac{1}{d^*(x)}\text{dir}_\Gamma\phi. $$	
Similarly, we first linearize the functions $\dfrac{1}{d^*(x)^2}(x-y^*)$ and $\dfrac{1}{d^*(x)}$ in the standard way:
$$\dfrac{1}{d^*({\bf x})}\approx \mathcal{L}_{\frac{1}{d^*}}({\bf x})= \dfrac{1}{d^*({\bf x_0})} - \dfrac{2}{d^*({\bf x_0})^2}[{\bf x_0}-{\bf y^*}]\cdot [{\bf x}-{\bf x_0}],$$
$$\dfrac{1}{d^*({\bf x})^2}({\bf x}-{\bf y^*})\approx \mathcal{L}_{\nabla\frac{1}{d^*}}({\bf x})= \dfrac{1}{d^*({\bf x_0})^2}[{\bf x_0}-{\bf y^*}]+\left[\mathcal{H}\left(\frac{1}{d^*({\bf x_0})}\right)\right] [{\bf x}-{\bf x_0}].$$
Here $\mathcal{H}\left(\dfrac{1}{d^*({\bf x_0})}\right)$ is the Hessian matrix of the function $\dfrac{1}{d^*}$ at the point ${\bf x_0}$. Similarly, we take $x_0$ to be the points ${\bf x}$ on $\Gamma_h^n$ and then obtain the following equations for all ${\bf x}\in \Gamma_h^n$:
\begin{equation}
	\mathcal{L}_{\frac{1}{d^*}}({\bf U}_h^{n+1}({\bf x})) = \dfrac{1}{d^*({\bf x})} - \dfrac{2}{d^*({\bf x})^2}[{\bf x}-{\bf y^*}]\cdot [{\bf U}_h^{n+1}({\bf x})-{\bf x}],
\end{equation}
\begin{equation}
	\mathcal{L}_{\nabla\frac{1}{d^*}}({\bf U}_h^{n+1}({\bf x})) = \dfrac{1}{d^*({\bf x})^2}[{\bf x}-{\bf y^*}]+\left[\mathcal{H}\left(\frac{1}{d^*({\bf x})}\right)\right] [{\bf U}_h^{n+1}({\bf x})-{\bf x}].
\end{equation}
With this approximation, we construct the formula for $dD^{n+1}(\Gamma_h^n; \Phi)$:
\begin{equation}\begin{array}{ll}
	dD^{n+1}(\Gamma_h^n; \Phi) =& -2\sum_{K\subset \Gamma_h^n}\int_K\left(\dfrac{1}{d_K({\bf x})^2}[{\bf x}-{\bf y_K}]\cdot \Phi +\tau\left[\mathcal{H}\left(\dfrac{1}{d_K({\bf x})}\right){\bf V}^{n+1}\right]\cdot \Phi
	\right)\\
	\ & +\sum_{K\subset \Gamma_h^n}\int_K\left(\dfrac{1}{d_K({\bf x})}\text{dir}_\Gamma\Phi - \dfrac{2\tau}{d_K({\bf x})^2}[{\bf x}-{\bf y}_K]\cdot [{\bf V}^{n+1}]\,\text{dir}_{\Gamma}\Phi
	\right)\end{array}
\end{equation}
Note that $d_K({\bf x})$ is exactly the distance function $d^*({\bf x})$ by taking the point ${\bf y^*}$ dependent on each element $K$, denoted by ${\bf y}_k$. Understandably, ${\bf y}_k$ is taken in the multiple vesicles case where ${\bf y}_k$ is a point on one of the vesicles other than $\Gamma_h^n$ such that it is closest to the element $K$.
\section{Algorithm}
The problem is fully discretized as discussed in the last chapter. Since the {\bf Problem 4.3.2.} is the full problem when the other two problems can be obtained by taking either $\alpha$ or $\beta$ to be zero, we illustrate the full algorithm for this problem in this chapter. We now have the {\bf Problem 4.3.2.} discretized and linearized, which make it standard to be solved by FEM, except that the area constraint $A(\Gamma_h^{n+1}) = A(\Gamma_h^0)$ need to be reconsidered. Hence, we start our discussion from the part of 'area constraint'.
\subsection{Area Constraint}
The method to treat the area constraint is based on the method presented in \cite{Nochetto2010} to compute the Lagrange multiplier $\lambda^{n+1}$, but with some nontrivial difference.

First, recall the equations for {\bf Problem 4.3.2.} here for convenience:
\begin{equation}\label{dWdHdD}\begin{array}{ll}
\int_{\Gamma_h^n} {\bf V}^{n+1}\cdot {\Phi}	 = &-dW^{n+1}(\Gamma_h^n;\Phi)- \alpha dH^{n+1}(\Gamma_h^n; \Phi) - \beta dD^{n+1}(\Gamma_h^n;\Phi)\\ &-\lambda^{n+1}dA^{n+1}(\Gamma_h^n;\Phi),
\end{array}
\end{equation}
\begin{equation}\label{vhx}	\tau\int_{\Gamma^n} \nabla_{\Gamma^n}{\bf V}^{n+1}\colon \nabla_{\Gamma^n}{\Phi} = \int_{\Gamma^n} {\bf H}^{n+1}\cdot {\Phi}-\int_{\Gamma^n}\nabla_{\Gamma^n}{\bf x}\colon\nabla_{\Gamma^n}{\Phi},
\end{equation}
\begin{equation}\label{area} A(\Gamma_h^{n+1}) = A(\Gamma_h^{0}).\end{equation}

In their method, they make the shape derivative $dA$ in the above equation to be explicit \cite{Nochetto2010}. That is to say, it becomes $dA^{n}(\Gamma_h^n;\Phi)=\int{\bf H}^n\cdot \Phi$ instead of $dA^{n+1}(\Gamma_h^n;\Phi)=\int{\bf H}^{n+1}\cdot \Phi$ as stated in our method. However, when I try to make it explicit as they said and implement the method in Matlab, the program always breaks down. Even if I take $\alpha$ and $\beta$ to be zero, which makes our model very similar to their problem, the numerical results still break down. However, when I use the implicit formula $dA^{n+1}$ as in equation \eqref{dWdHdD}, the numerical results reveal to be stable and convergent. Hence, we apply the ideas of computing Lagrange Multiplier as stated in \cite{Nochetto2010}, but use $dA^{n+1}$. The idea for solving our problem exactly stated as in equations (\ref{dWdHdD} -~\ref{area}) is explained below.

Rewrite the pair of unknown as $({\bf V}^{n+1}, {\bf H}^{n+1}) = ({\bf V}^{n+1}_1, {\bf H}^{n+1}_1) + \lambda^{n+1}({\bf V}^{n+1}_2, {\bf H}^{n+1}_2)$. We solve $({\bf V}^{n+1}_1, {\bf H}^{n+1}_1)$ and $({\bf V}^{n+1}_2, {\bf H}^{n+1}_2)$ separately corresponding to the {\bf Problem 4.4.1} and {\bf Problem 4.4.2}, and then find the Lagrange Multiplier $\lambda^{n+1}$ such that the resulted $\Gamma_h^{n+1}$ updated by ${\bf V}^{n+1}={\bf V}^{n+1}_1 +\lambda^{n+1}{\bf V}^{n+1}_2$ satisfies the area constraint.

\begin{p} Find $({\bf V}^{n+1}_1, {\bf H}^{n+1}_1)$ such that $\forall \Phi\in \mathbb{F}(\Gamma_h^n)$,
	$$\int_{\Gamma_h^n} {\bf V}_1^{n+1}\cdot {\Phi}	 = -dW^{n+1}(\Gamma_h^n;\Phi)- \alpha dH^{n+1}(\Gamma_h^n; \Phi) - \beta dD^{n+1}(\Gamma_h^n;\Phi),$$
	$$\tau\int_{\Gamma^n} \nabla_{\Gamma^n}{\bf V}_1^{n+1}\colon \nabla_{\Gamma^n}{\Phi} = \int_{\Gamma^n} {\bf H}_1^{n+1}\cdot {\Phi}-\int_{\Gamma^n}\nabla_{\Gamma^n}{\bf x}\colon\nabla_{\Gamma^n}{\Phi}.
	$$
	\end{p}

\begin{p} Find $({\bf V}^{n+1}_2, {\bf H}^{n+1}_2)$ such that $\forall \Phi\in \mathbb{F}(\Gamma_h^n)$,
	$$\int_{\Gamma_h^n} {\bf V}_2^{n+1}\cdot {\Phi}	 = -dA^{n+1}(\Gamma_h^n;\Phi),$$
	$$\tau\int_{\Gamma_h^n} \nabla_{\Gamma_h^n}{\bf V}_2^{n+1}\colon \nabla_{\Gamma_h^n}{\Phi} = \int_{\Gamma_h^n} {\bf H}_2^{n+1}\cdot {\Phi}-\int_{\Gamma_h^n}\nabla_{\Gamma_h^n}{\bf x}\colon\nabla_{\Gamma_h^n}{\Phi}.
	$$
\end{p}
The above two problems are standardly formulated to be solved by FEM. Now the values of $({\bf V}^{n+1}_1, {\bf H}^{n+1}_1)$ and $({\bf V}^{n+1}_2, {\bf H}^{n+1}_2)$ are known and we only need to find the suitable $\lambda^{n+1}$. More specifically, by equation \label{area}, we solve the Lagrange Multiplier $\lambda^{n+1}\in \mathbb{R}$ by Newton's Method as a root of the function
$$f_N(\lambda) = A\left(\Gamma_h^n + \tau {\bf V}_\lambda(\Gamma_h^n)\right)-A(\Gamma_h^0),$$
where $${\bf V}_\lambda(\cdot)={\bf V}^{n+1}_1(\cdot) +\lambda{\bf V}^{n+1}_2(\cdot).$$
The detailed information about the derivation of the differential of $f_N$ and the initial guess for this Newton's method is given by \cite{Nochetto2010}. The iterative equation is
$$\lambda_i^{n+1} = \lambda_{i-1}^{n+1} - [\mathcal{D}f_N(\lambda_{i-1}^{n+1})]^{-1}f_N(\lambda_{i-1}^{n+1}),$$
where $\mathcal{D}f_N(\lambda_{i-1}) = \tau \int_{\Gamma^n(\lambda)}\text{dir}_{\Gamma^n(\lambda)}{\bf V}_2^{n+1}$. The initial guess is
$$\lambda_0^{n+1} = -\left(\int_{\Gamma_h^n} \text{dir}_{\Gamma_h^n}\,{\bf V}_2^{n+1}\right)^{-1}\int_{\Gamma_h^n} \text{dir}_{\Gamma_h^n}\,{\bf V}_1^{n+1}.$$

\subsection{Full Algorithm}
\begin{algorithm}[H]
    \caption{Full Algorithm for {\bf Problem 4.3.2}}
    \label{fullalgo}
    \begin{algorithmic}[1] % The number tells where the line numbering should start
        \Procedure{FEM}{$\Gamma_h^0,T$} %\Comment{The g.c.d. of a and b}
            \State Given initial discretized shapes $\{\Gamma_h^{i,0}\}_{i=1}^M$
            \State Given essential datas: the obstacles ${\bf B}$; weights $\alpha$ and $\beta$
            \State Set time-step size $\tau$, tolerance $\epsilon$ and maximum number of iterations $N$
            \State $T=0$
            \State $n=0$ 
            \While{$1$} \Comment{Breaking rule is in the loop}
                \For{$m=1:M$}\Comment{Move the shapes one by one}
                    \State $k \gets \text{TOEP}(mod(n,M),m)$ 
                    \State Solve {\bf Problem 4.4.1} and {\bf 4.4.2} on $\mathbb{F}(\Gamma_h^{k,n})$
                    \State Use {\bf Newton's Method} to solve $\lambda^{k,n+1}$
                    \State Update $\Gamma_h^{k,n+1}$ by the results from line 9 and 10
                    \State Adjust the position of the vertices on $\Gamma_h^{k,n+1}$
                \EndFor
                \State $T \gets T+\tau$
                \State $n \gets n+1$ \Comment{Index of iteration}
                \If{$n > N$}
                    \State Break;
                \ElsIf{$|J_3(\lambda^{i,j+1},\Gamma_h^{i,j+1})-J_3(\lambda^{i,j},\Gamma_h^{i,j})|<\epsilon$ for all $i=1,2,\cdots,M$ and for $j=n,n-1,n-2$.}
                    \State Break;
                \EndIf
            \EndWhile            
        \EndProcedure
    \end{algorithmic}
\end{algorithm}

% no final time T

\newpage
\chapter{Numerical Examples}
In this chapter, we demonstrate various numerical examples to show the various applications of the models (from Model 1 to Model 3). In Section 5.1, the experiments are not aimed to mimic something. We just try different initial shapes to make as many interesting experiments as possible. Also, we plot some graphs to see the decrease of energy. In Section 5.2, the experiments are aimed to mimic the shape properties of the Golgi stacks, so the initial shapes are set goal-oriented. In this section, we provide additional information about the biological motivation of forming those special experiments. Also, we combine those numerical results with the observed image of Golgi cisternae and draw some conclusions based on that.  
\newpage
\section{Examples}
In this section, we present various examples, from the most basic one. The first example is a minimization of the functional $A(\Gamma)$. It is simply the length of $\Gamma$ in $\mathbb{R}^2$. This model is simple but contained in every model. Figure \ref{fig:A} shows that the ellipse shrinks to a point.

\subsubsection*{Example 1}
\begin{figure}[h!]
\centering
\subfloat[T=0]{
\includegraphics[width=60mm]{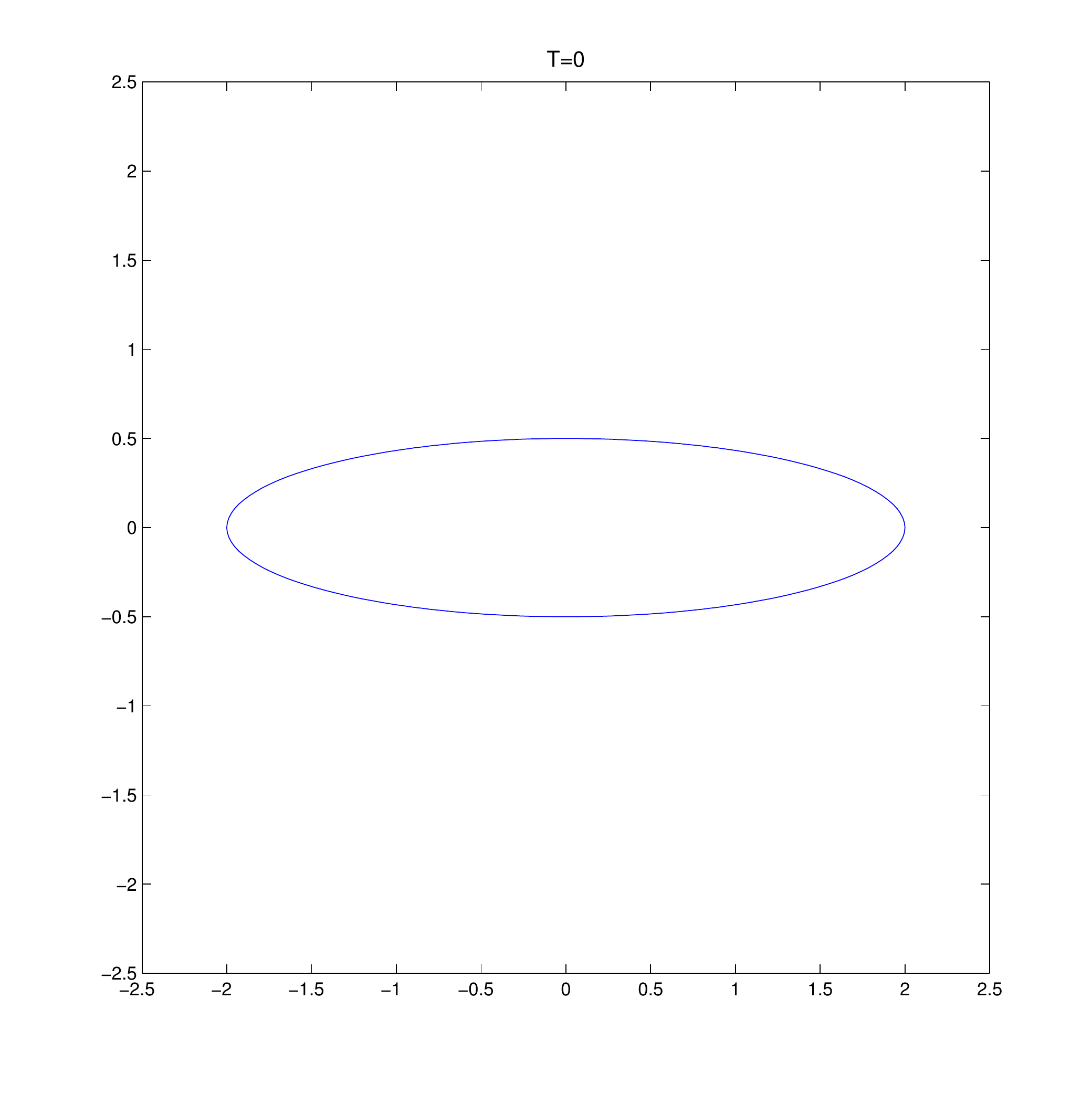}
}
\subfloat[T=0.29]{
\includegraphics[width=60mm]{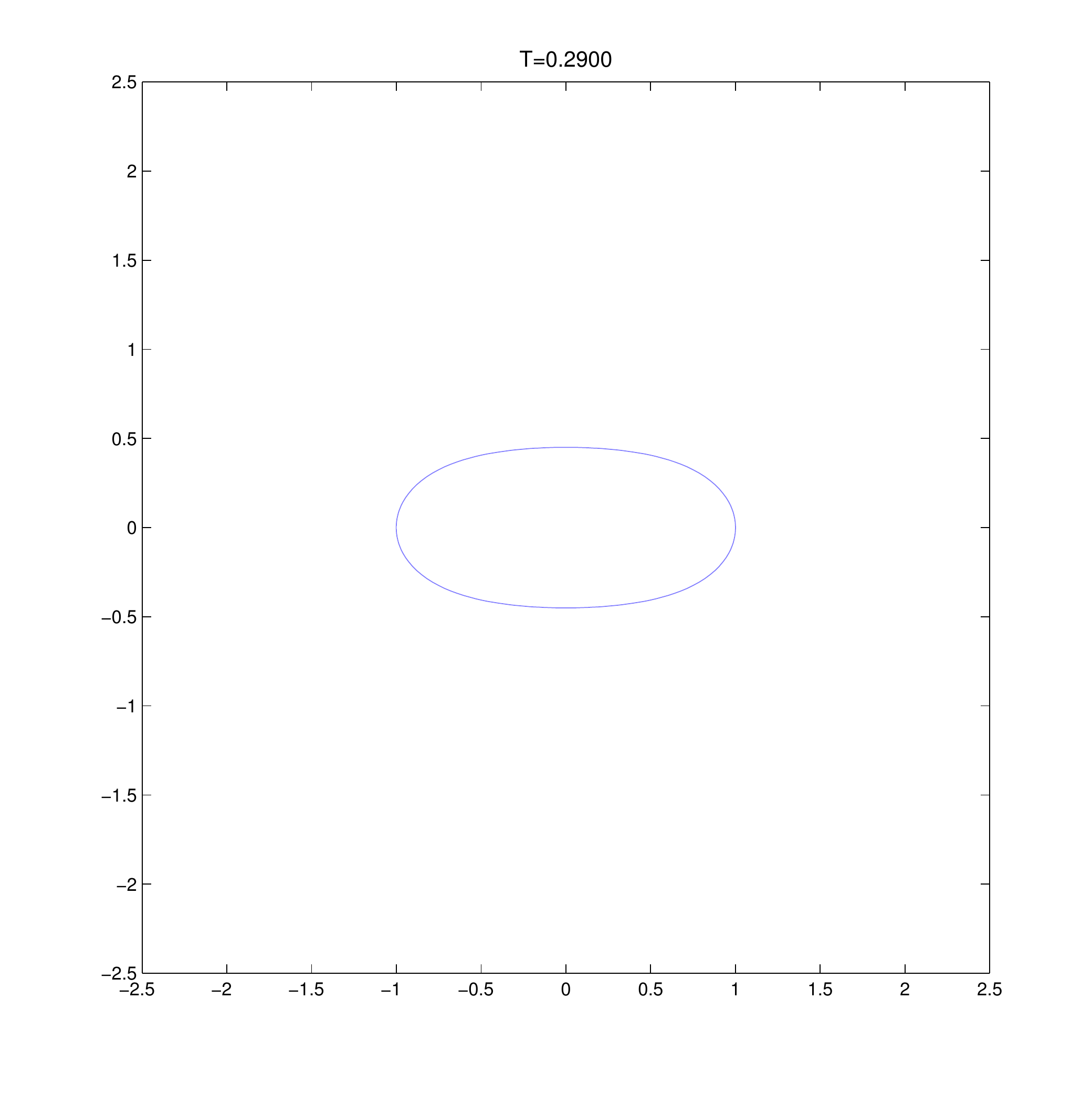}
}
\hspace{0mm}
\subfloat[T=0.48]{
\includegraphics[width=60mm]{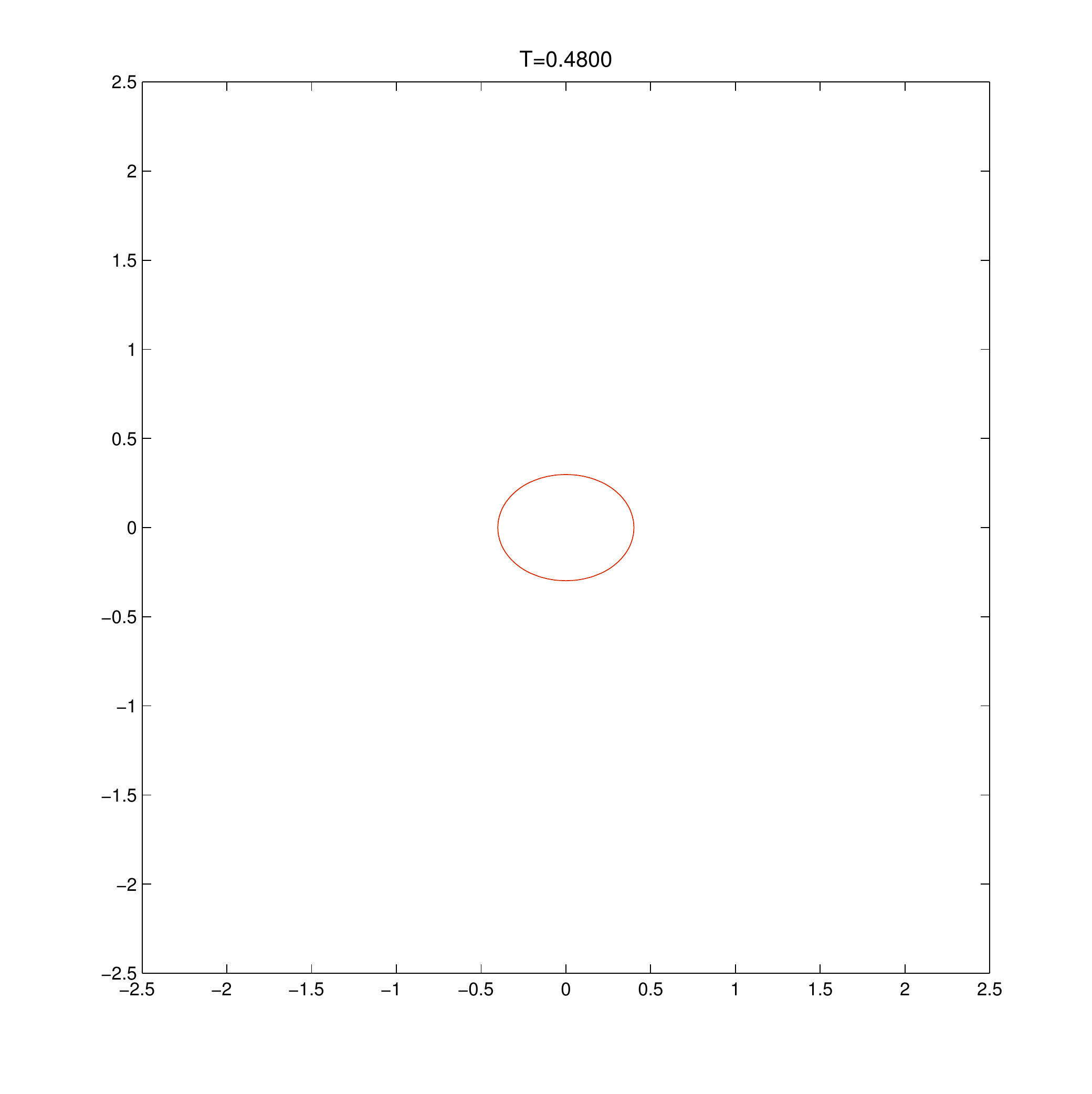}
}
\subfloat[T=0.56]{
\includegraphics[width=60mm]{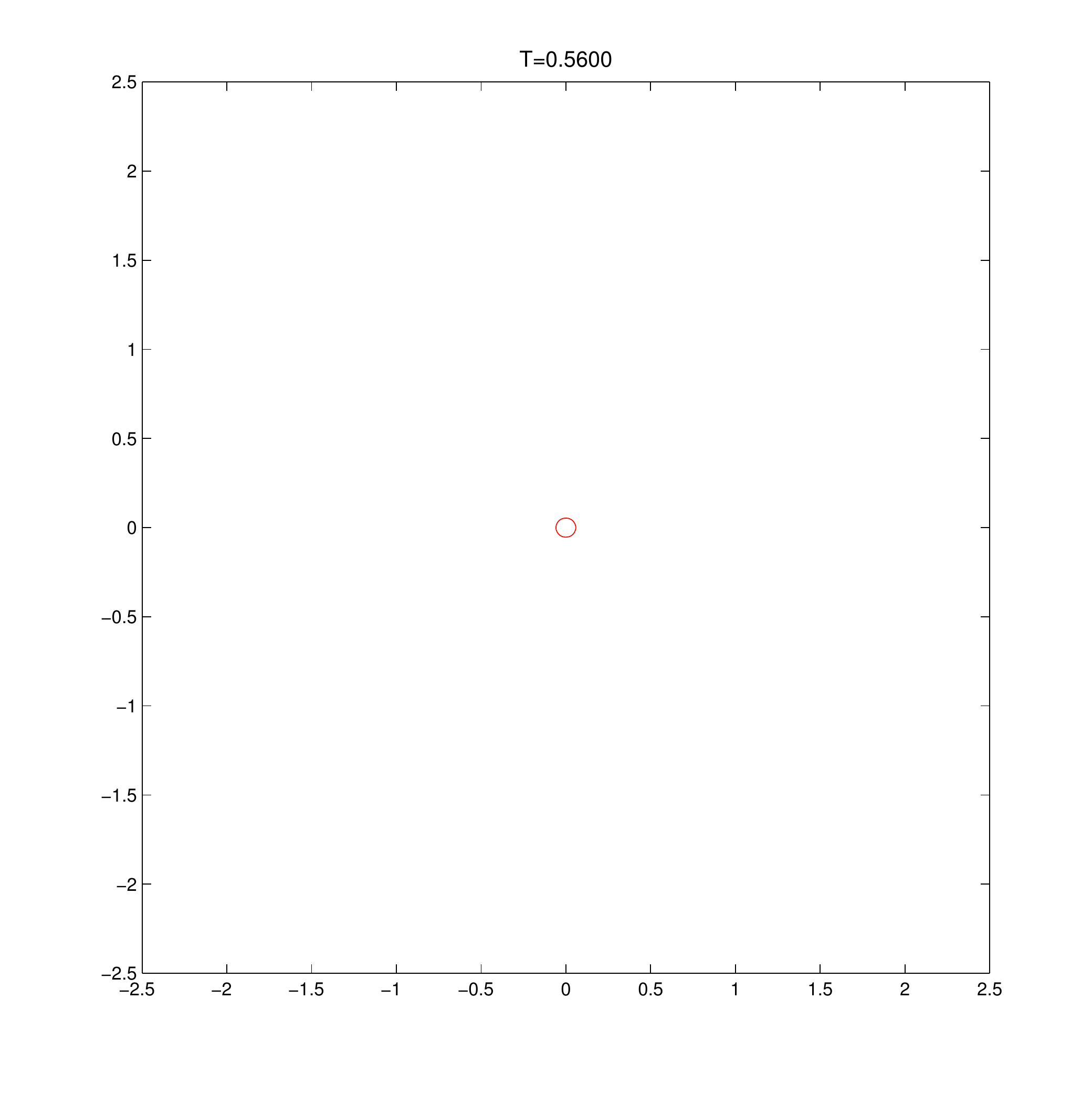}
}
\caption{Minimization of $A(\Gamma)$ (To be continued)}
\end{figure}

\begin{figure}[h!]\ContinuedFloat
\centering
\subfloat[From T=0 to T=0.58]{
\includegraphics[width=60mm]{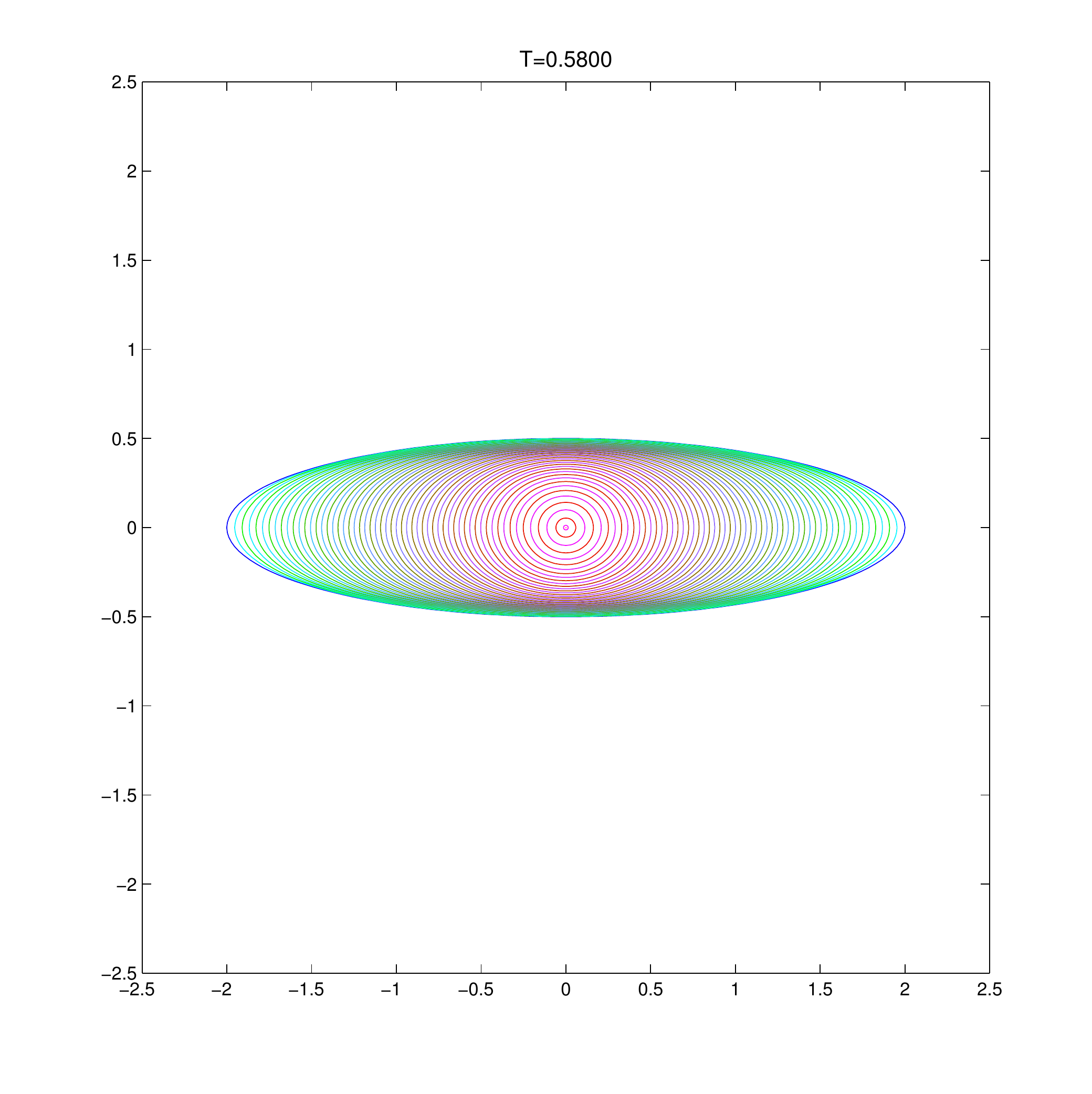}
}
\subfloat[Decrease of the Length]{
\includegraphics[width=60mm]{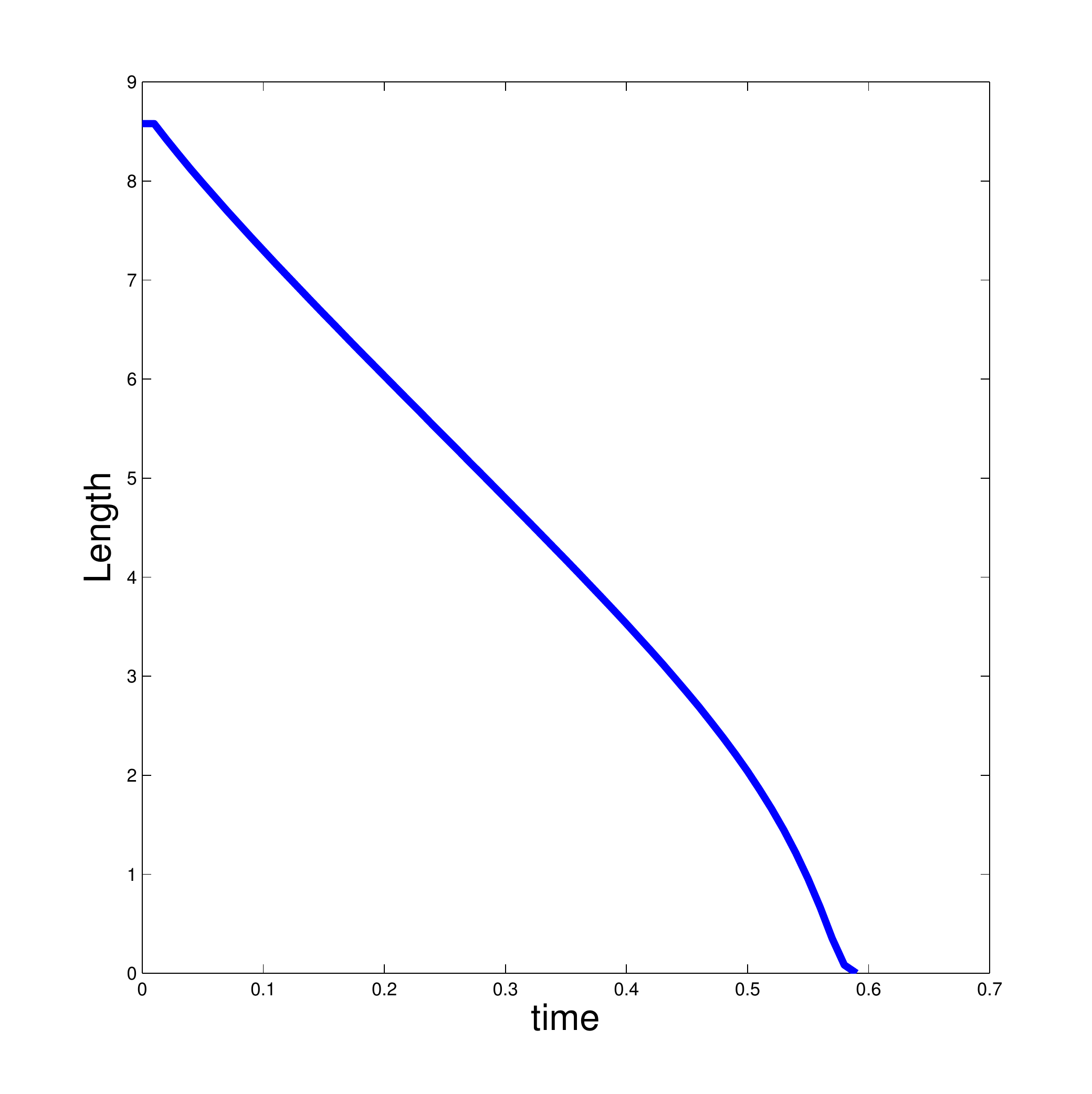}
}
\caption{(continue) Minimization of $A(\Gamma)$}
\label{fig:A}
\end{figure}

\subsubsection*{Example 2 － Model 1}
The second example is the minimization of Willmore Energy $W(\Gamma)$ under the condition that the area (length in 2D) $A(\Gamma)$ is fixed.
\begin{figure}[h!]
\centering
\subfloat[Initial shape]{
\includegraphics[width=65mm]{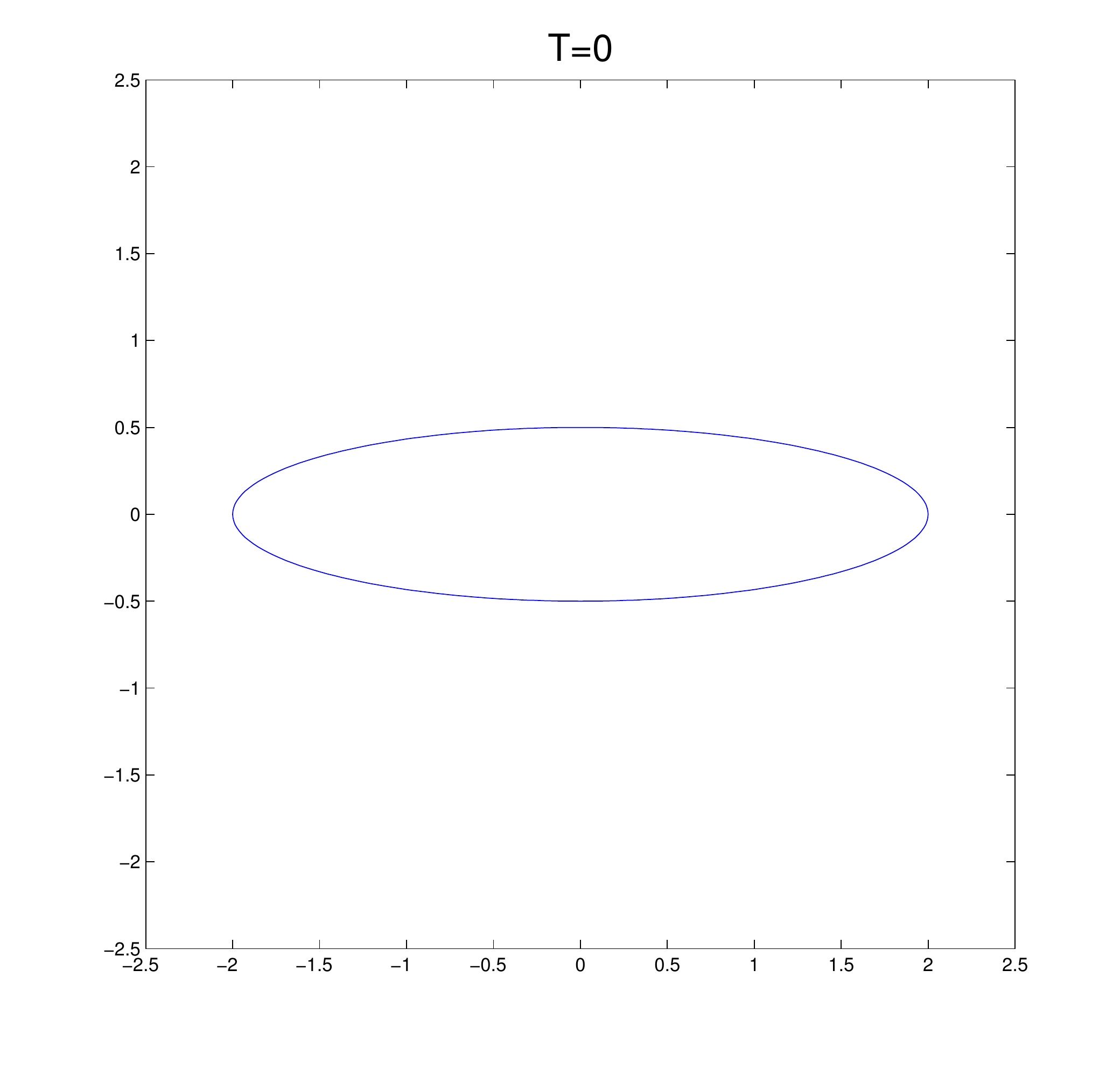}
}
\subfloat[]{
\includegraphics[width=65mm]{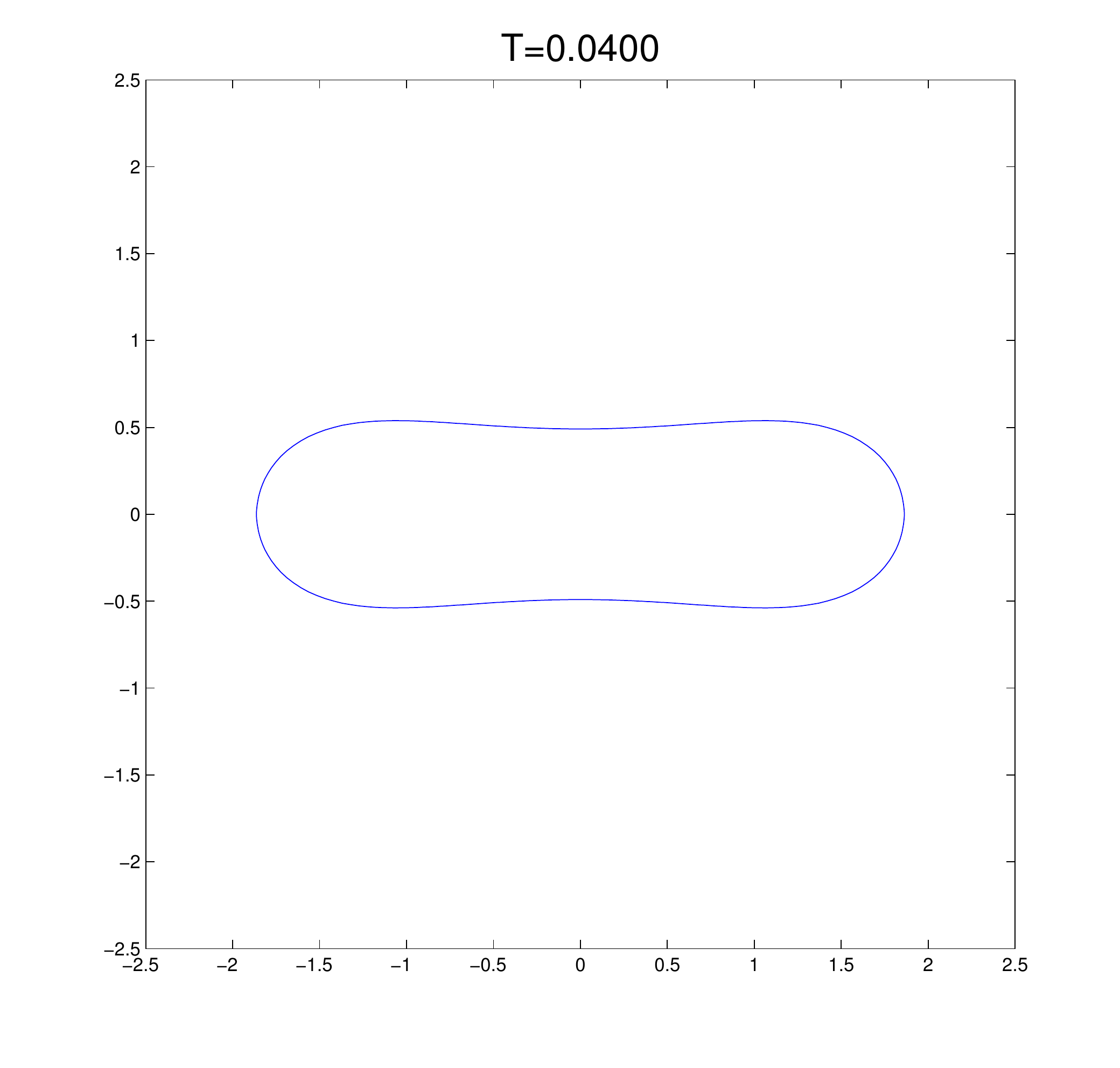}
}
\caption{Example 2 (to be continued)}	
\end{figure}
\begin{figure}[h!]\ContinuedFloat
\centering
\subfloat[]{
\includegraphics[width=65mm]{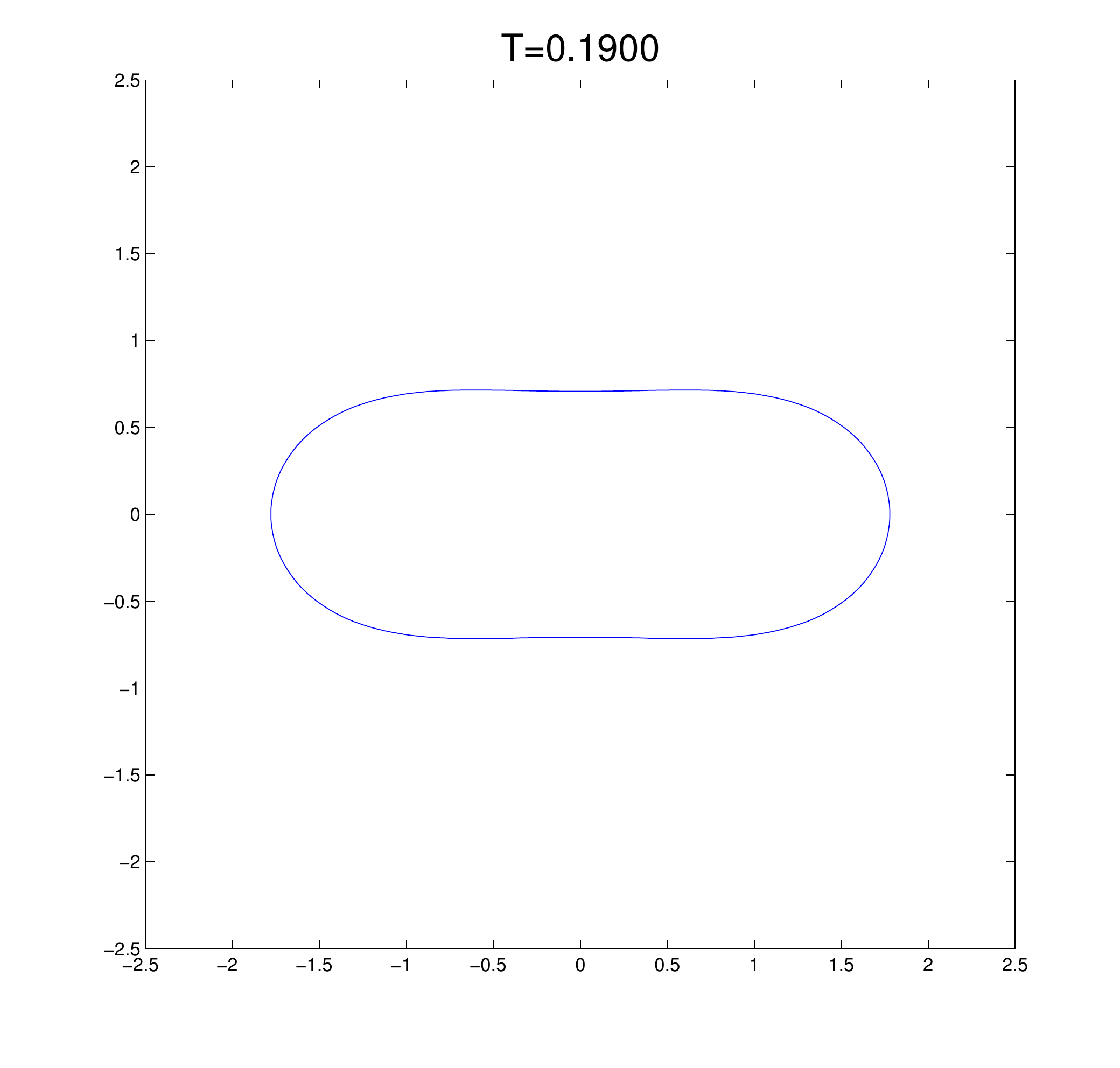}
}
\subfloat[]{
\includegraphics[width=65mm]{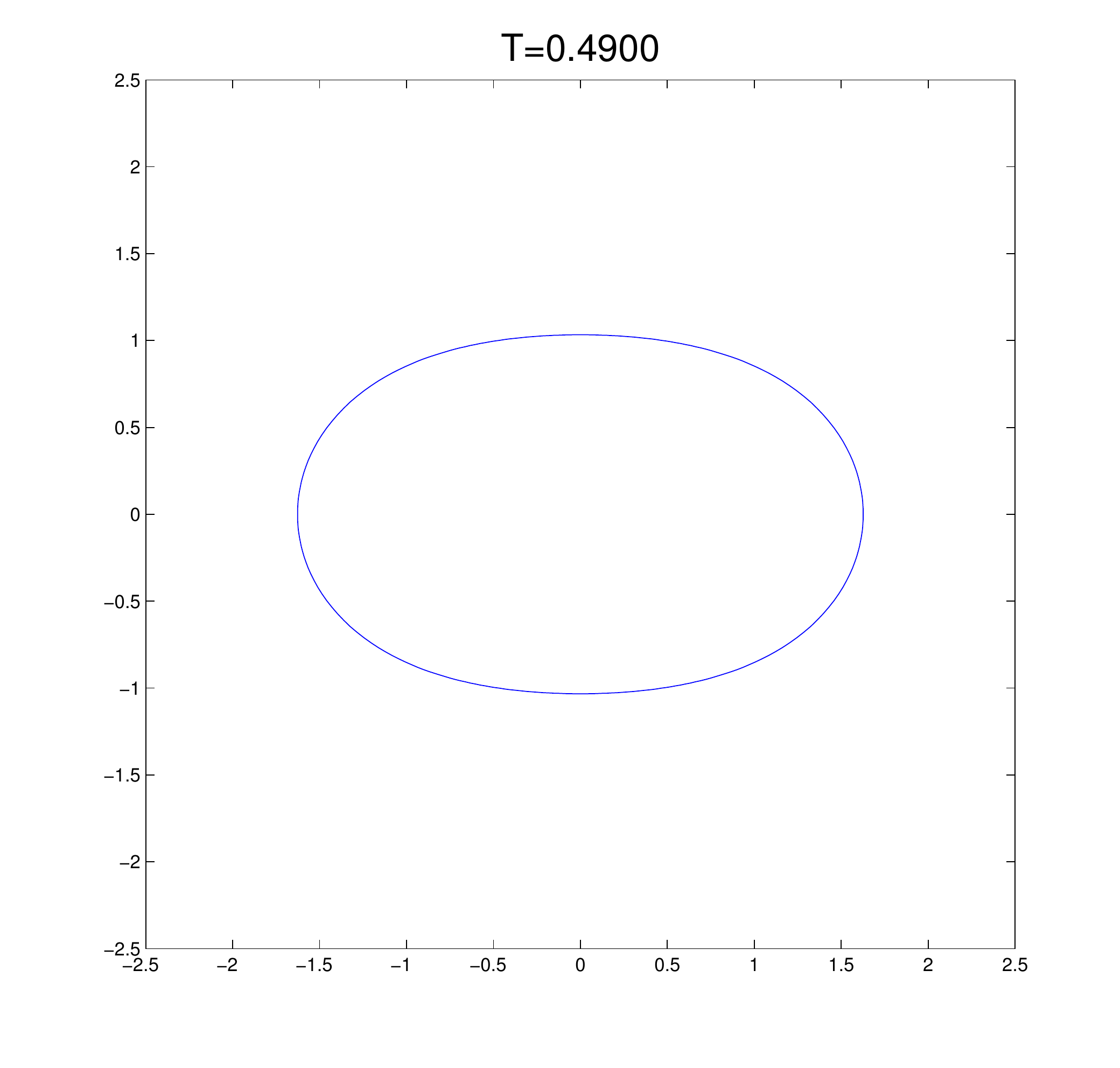}
}
\hspace{0mm}
\subfloat[]{
\includegraphics[width=65mm]{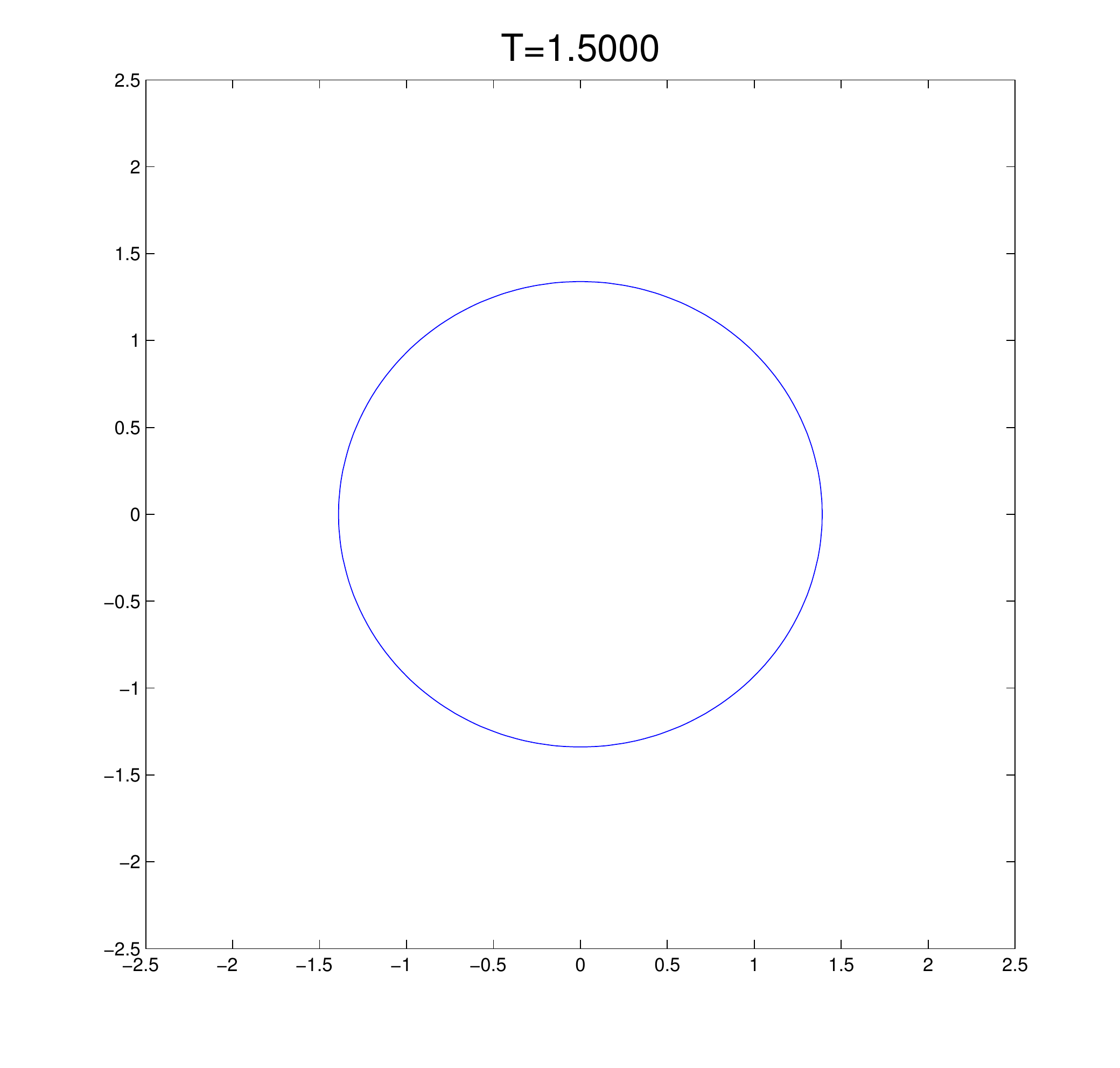}
}
\subfloat[From T=0 to T=1.50]{
\includegraphics[width=65mm]{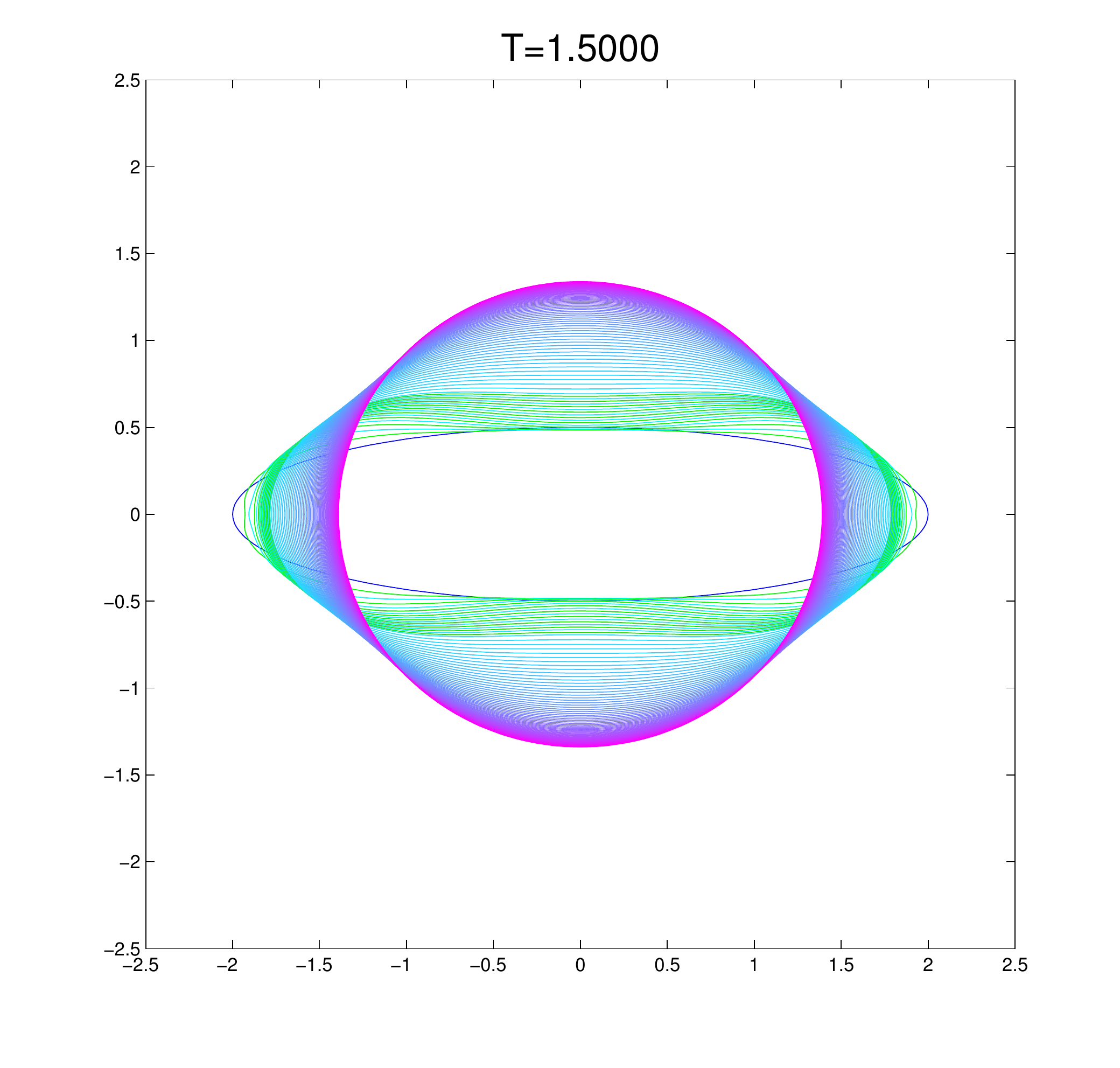}
}
\hspace{0mm}
\subfloat[Conservation of the Length $A(\Gamma)$]{
\includegraphics[width=65mm]{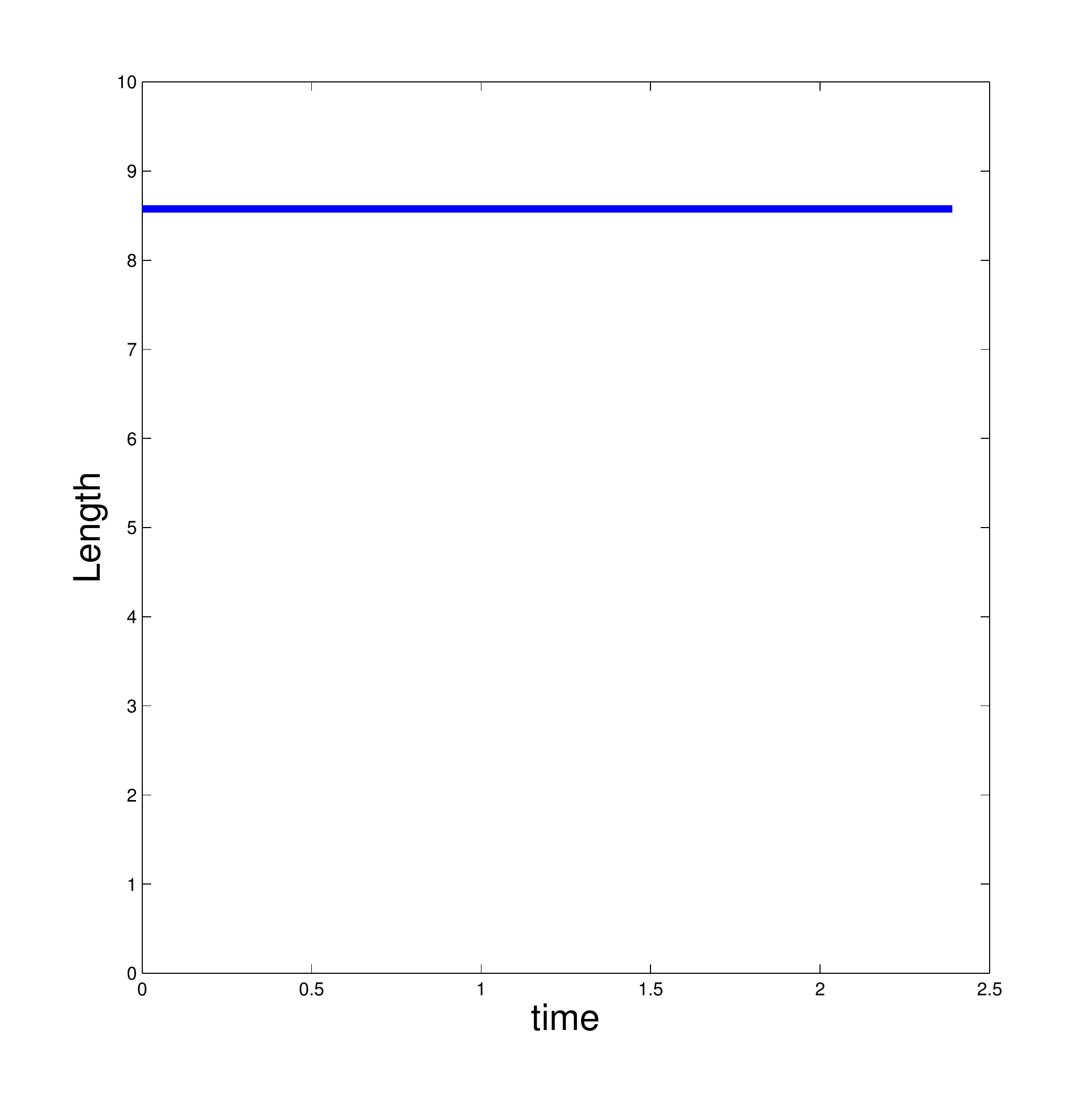}
}
\subfloat[Decrease of the Willmore Energy $W(\Gamma)$]{
\includegraphics[width=65mm]{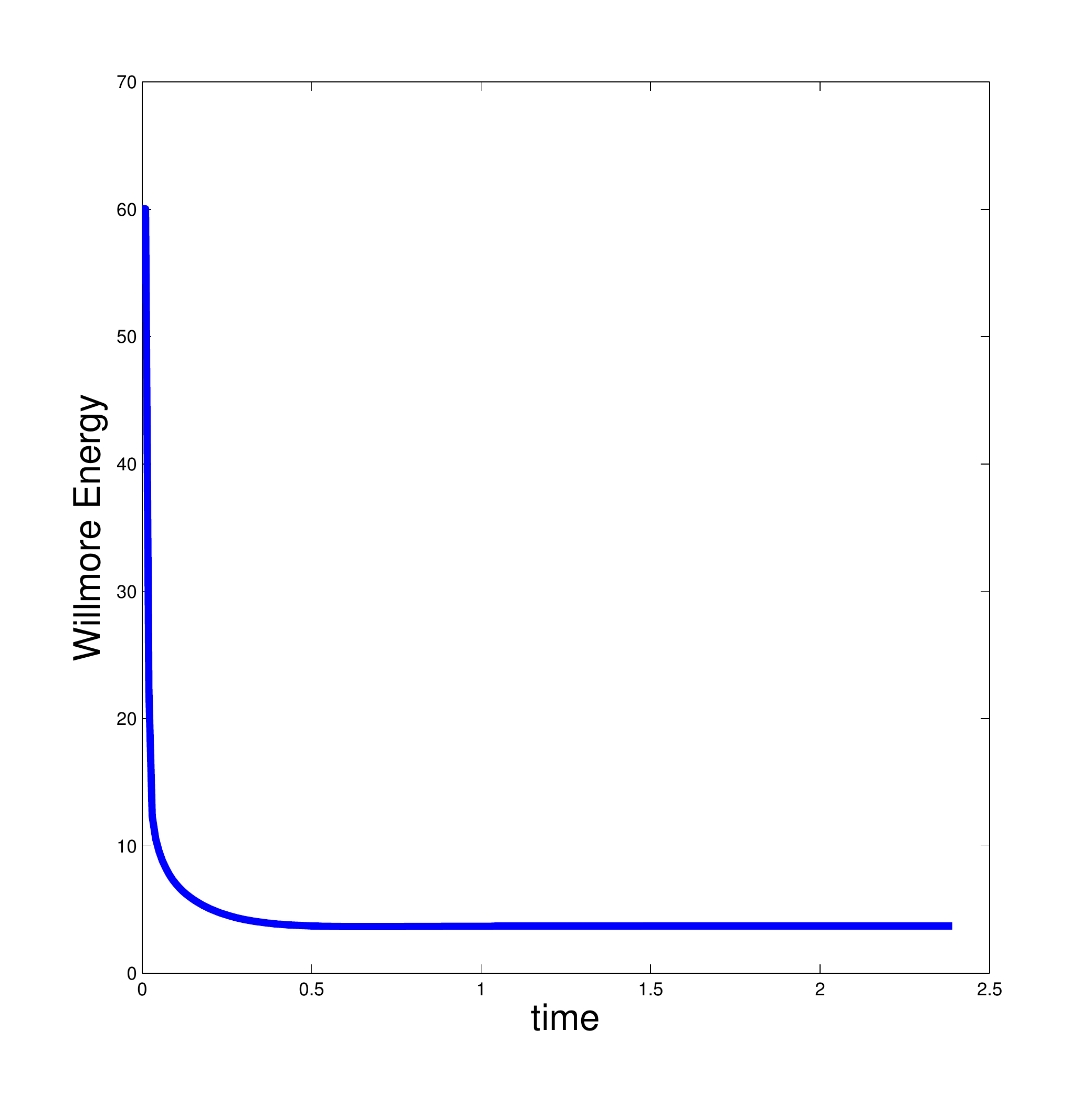}
}
\caption{Example 2}	
\end{figure}

\clearpage
\subsubsection*{Example 3 - Model 1}
This is similar to example 2 but it has a more interesting initial shape - the C shape.
\begin{figure}[h!]
\centering
\subfloat[Initial C shape]{
\includegraphics[width=0.5\textwidth]{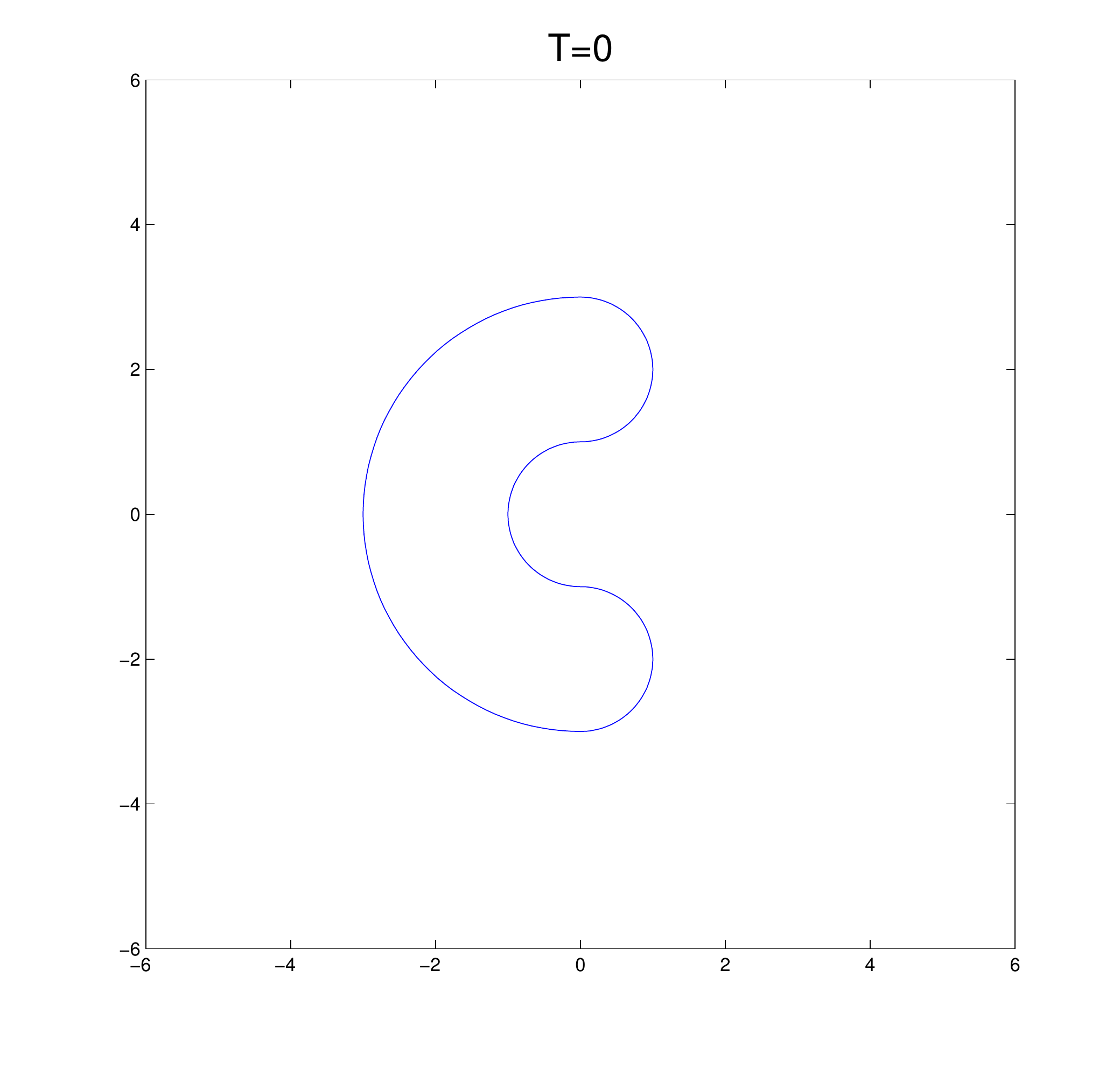}
}
\subfloat[]{
\includegraphics[width=0.5\textwidth]{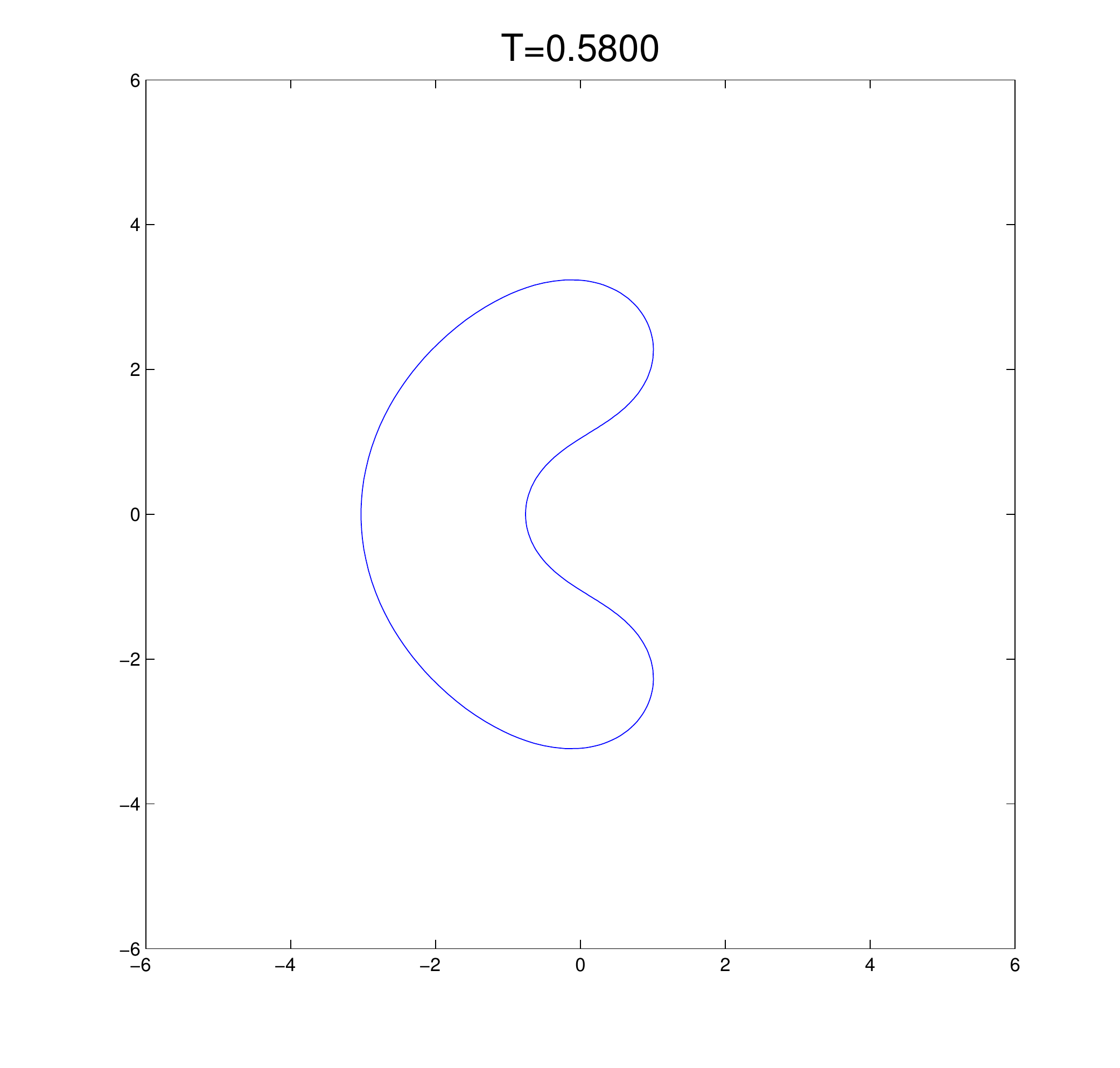}
}
\hspace{0mm}
\subfloat[]{
\includegraphics[width=0.5\textwidth]{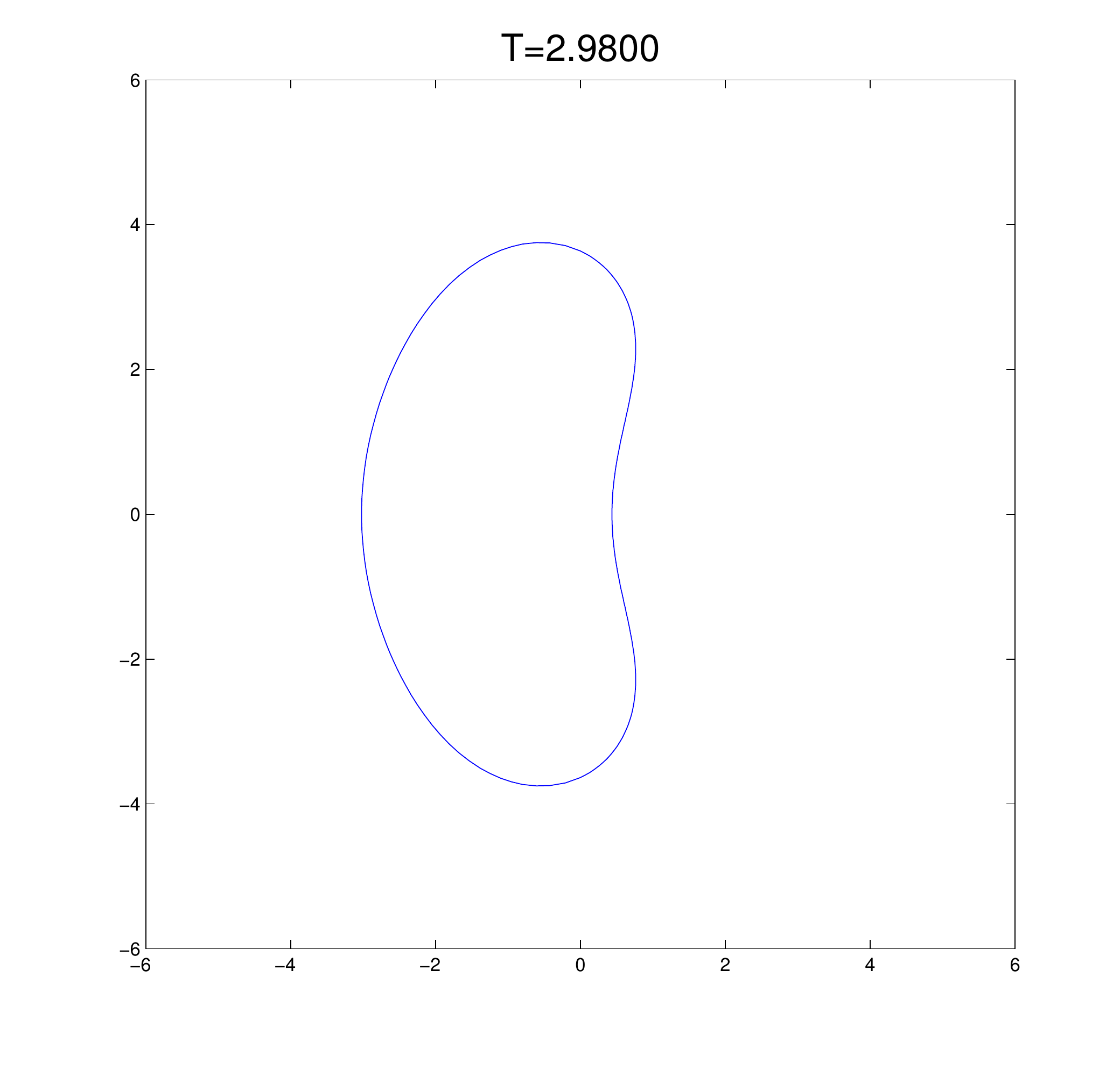}
}
\subfloat[]{
\includegraphics[width=0.5\textwidth]{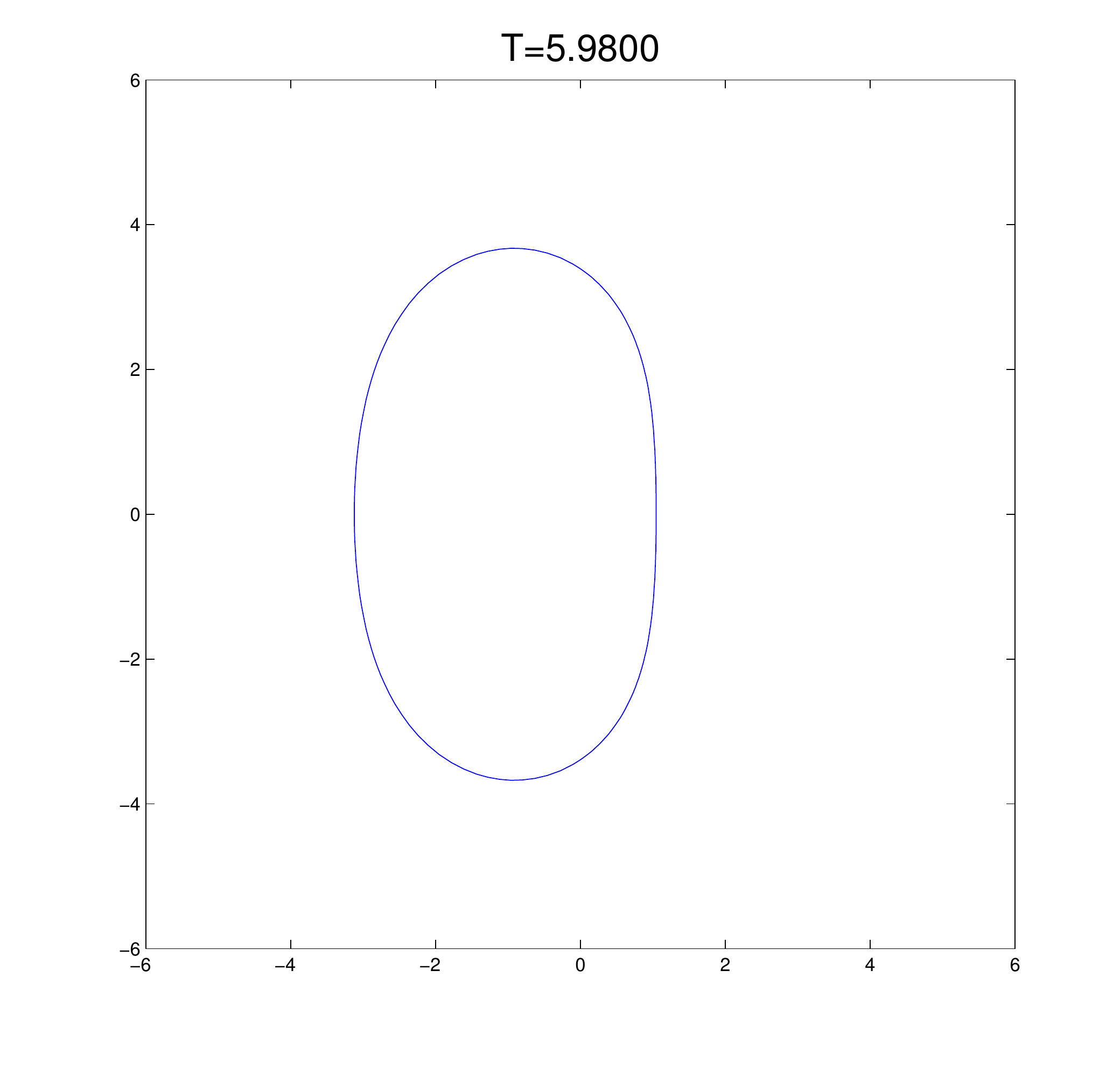}
}
\caption{C shape evolution (to be continued)}
\end{figure}

\begin{figure}[h!]\ContinuedFloat
\subfloat[T = 35.38]{
\includegraphics[width=0.5\textwidth]{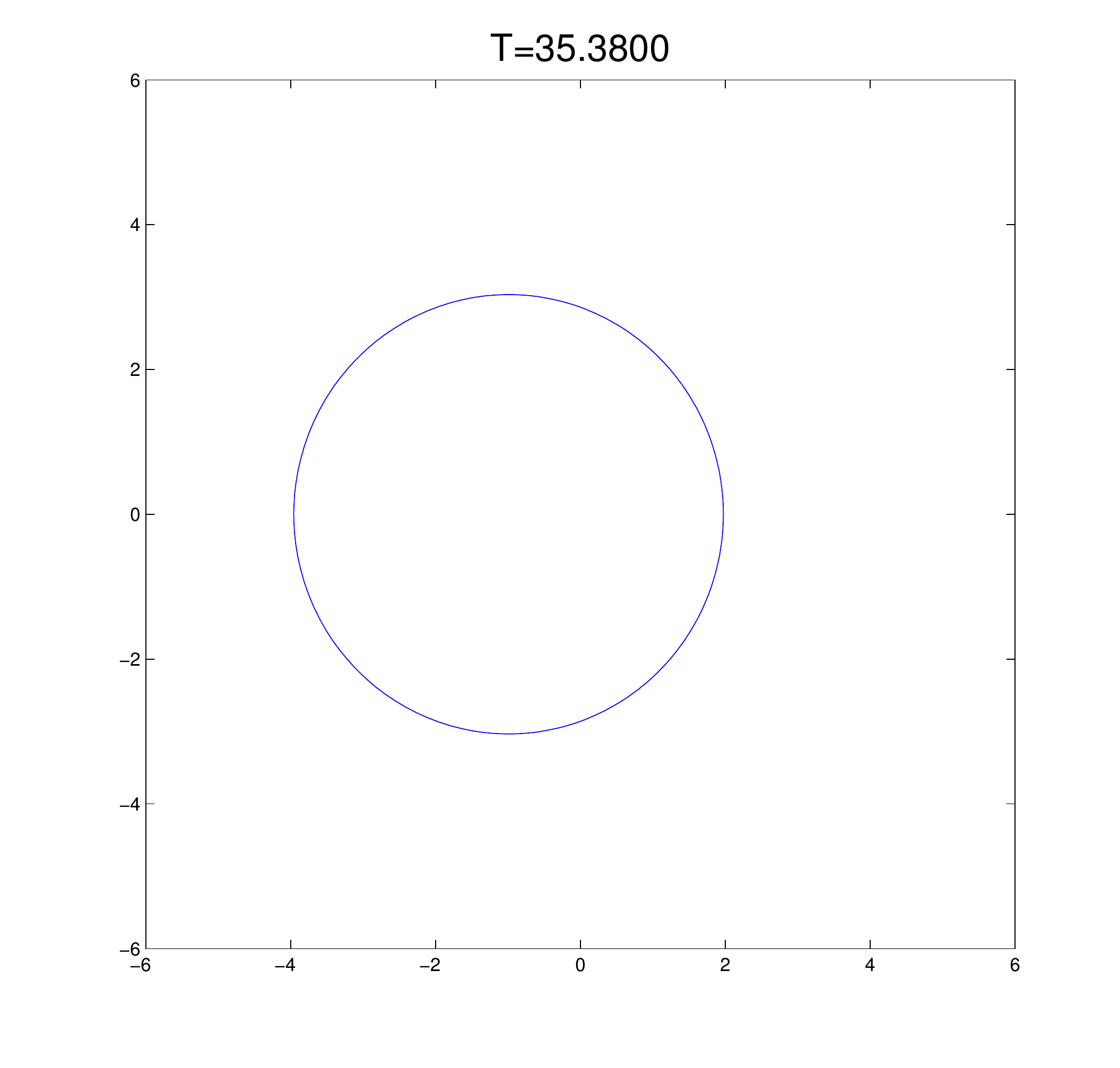}
}
\subfloat[From T=0 to T=35.38]{
\includegraphics[width=0.5\textwidth]{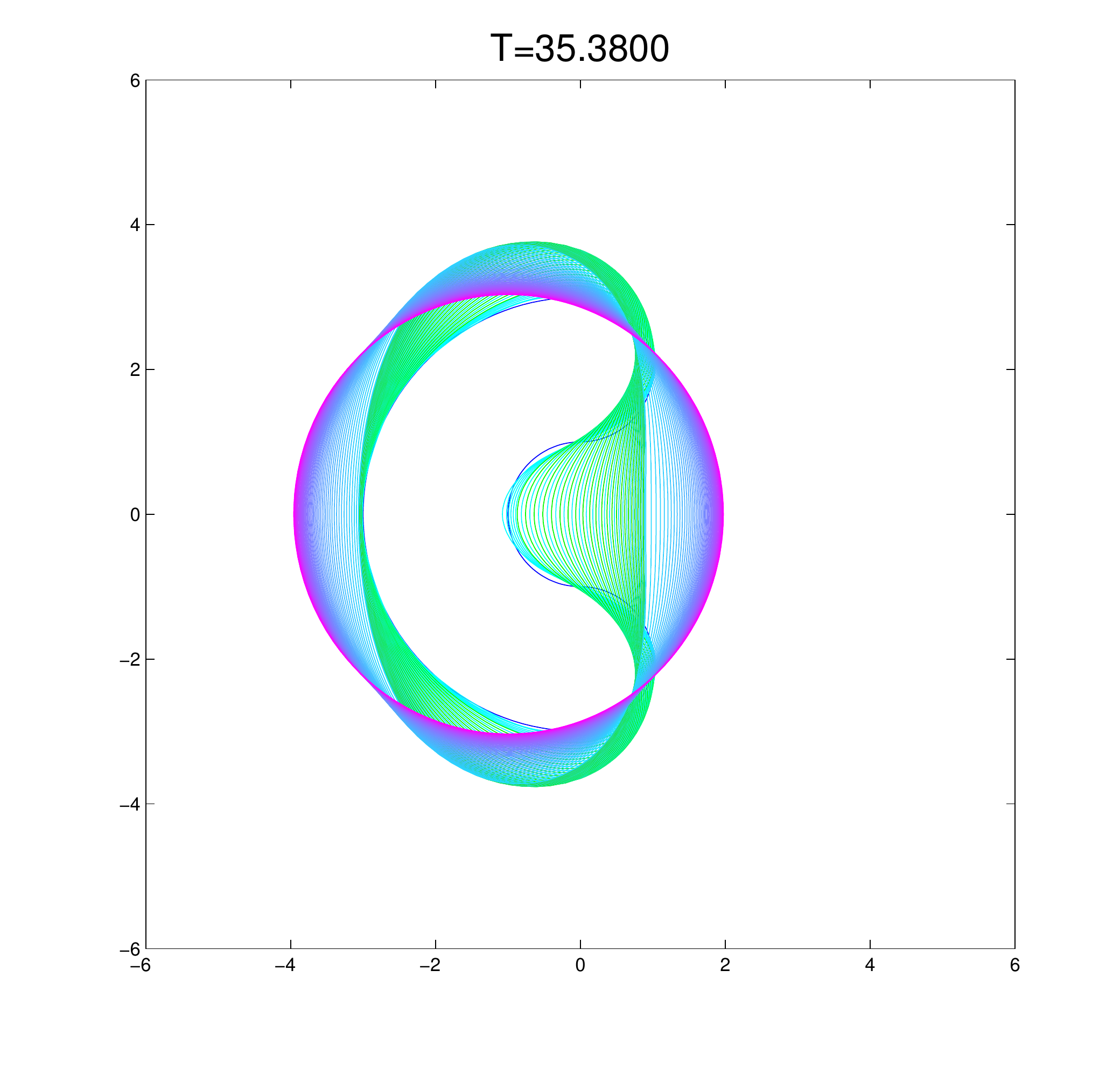}
}
\hspace{0mm}
\subfloat[Conservation of the Length $A(\Gamma)$]{
\includegraphics[width=0.5\textwidth]{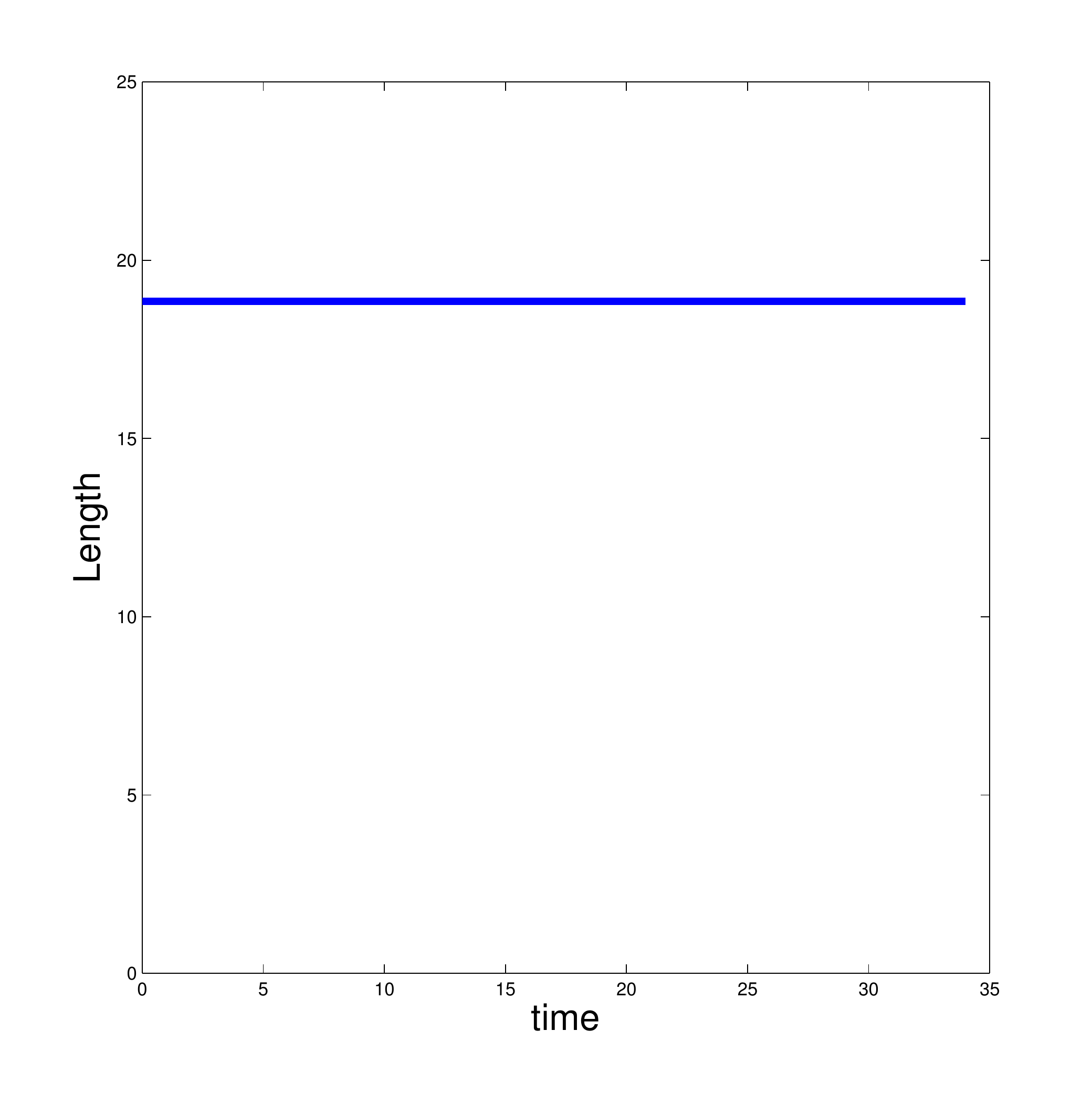}
}
\subfloat[Decrease of Willmore Energy $W(\Gamma)$]{
\includegraphics[width=0.5\textwidth]{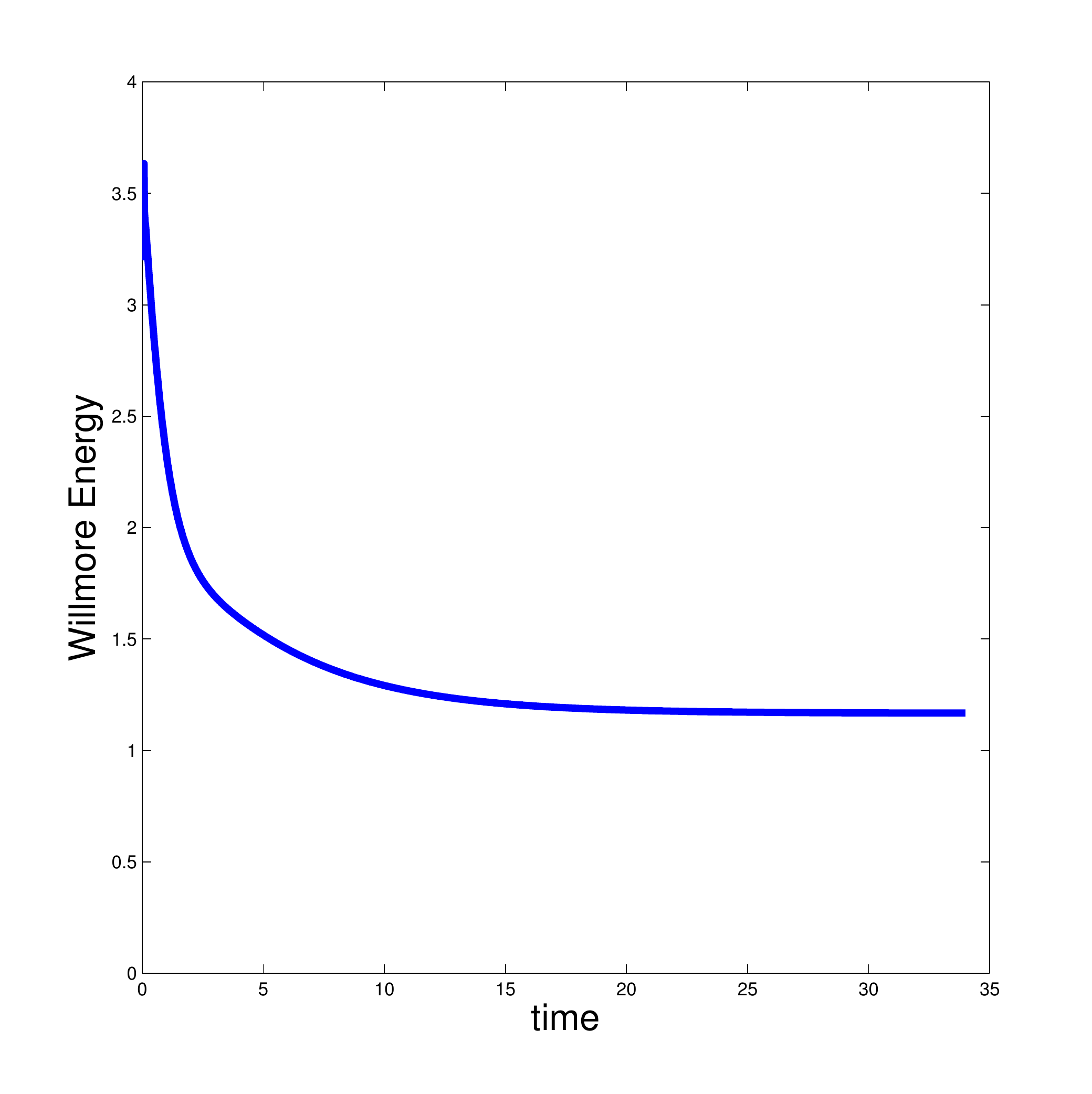}
}
\caption{(continue) C shape evolution}
\end{figure}

\clearpage
\subsubsection*{Example 4 - Model 2} %1 long 1 short
Now we begin to demonstrate some examples of Model 2, where the barriers/obstacles ${\bf B}$ is included. In all these examples, the dash lines or the grey block areas represent the region of the obstacles ${\bf B}$.

\begin{figure}[h!]
\centering
\subfloat[Initial shape - ellipse]{
\includegraphics[width=0.5\textwidth]{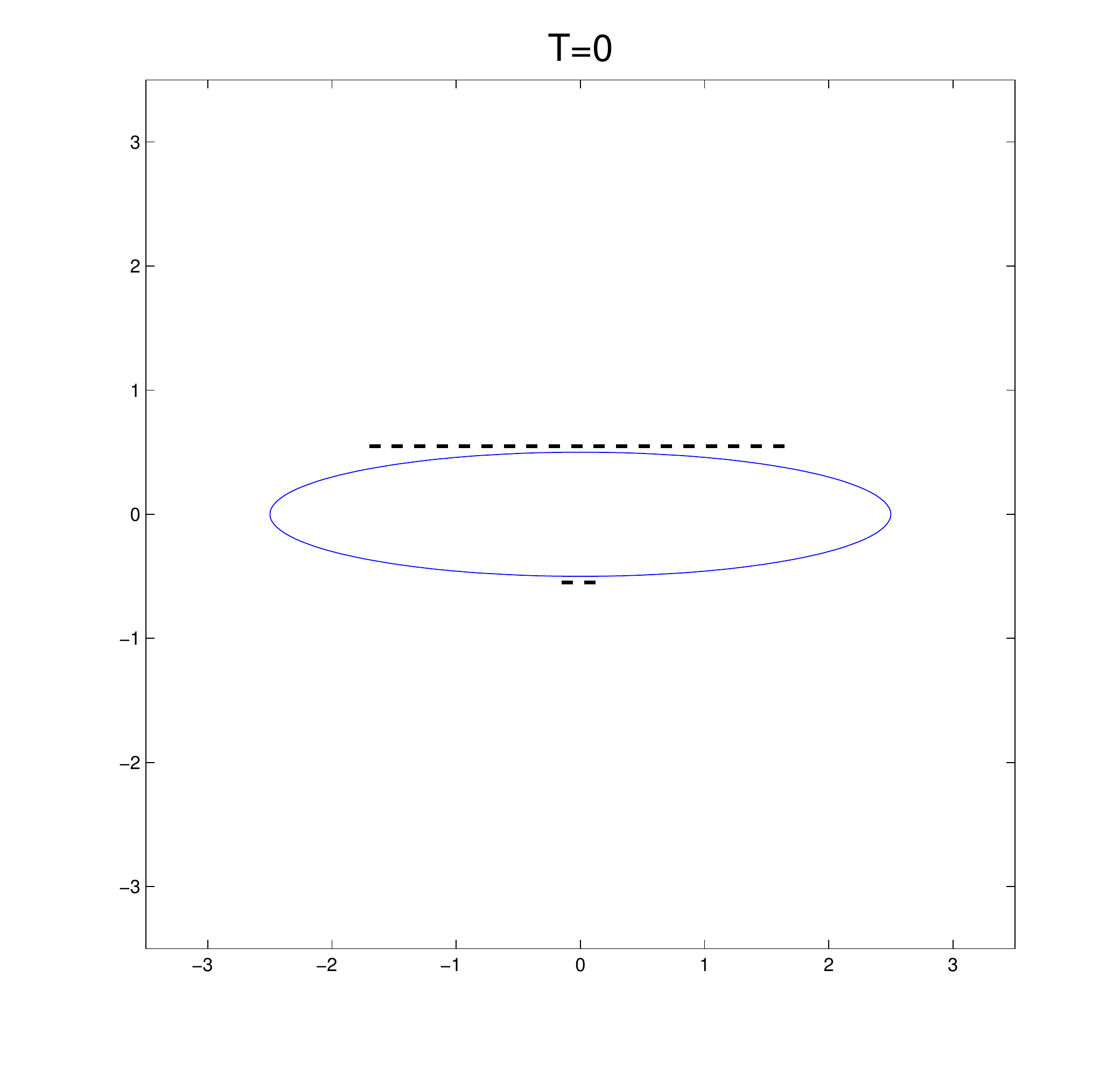}
}
\subfloat[]{
\includegraphics[width=0.5\textwidth]{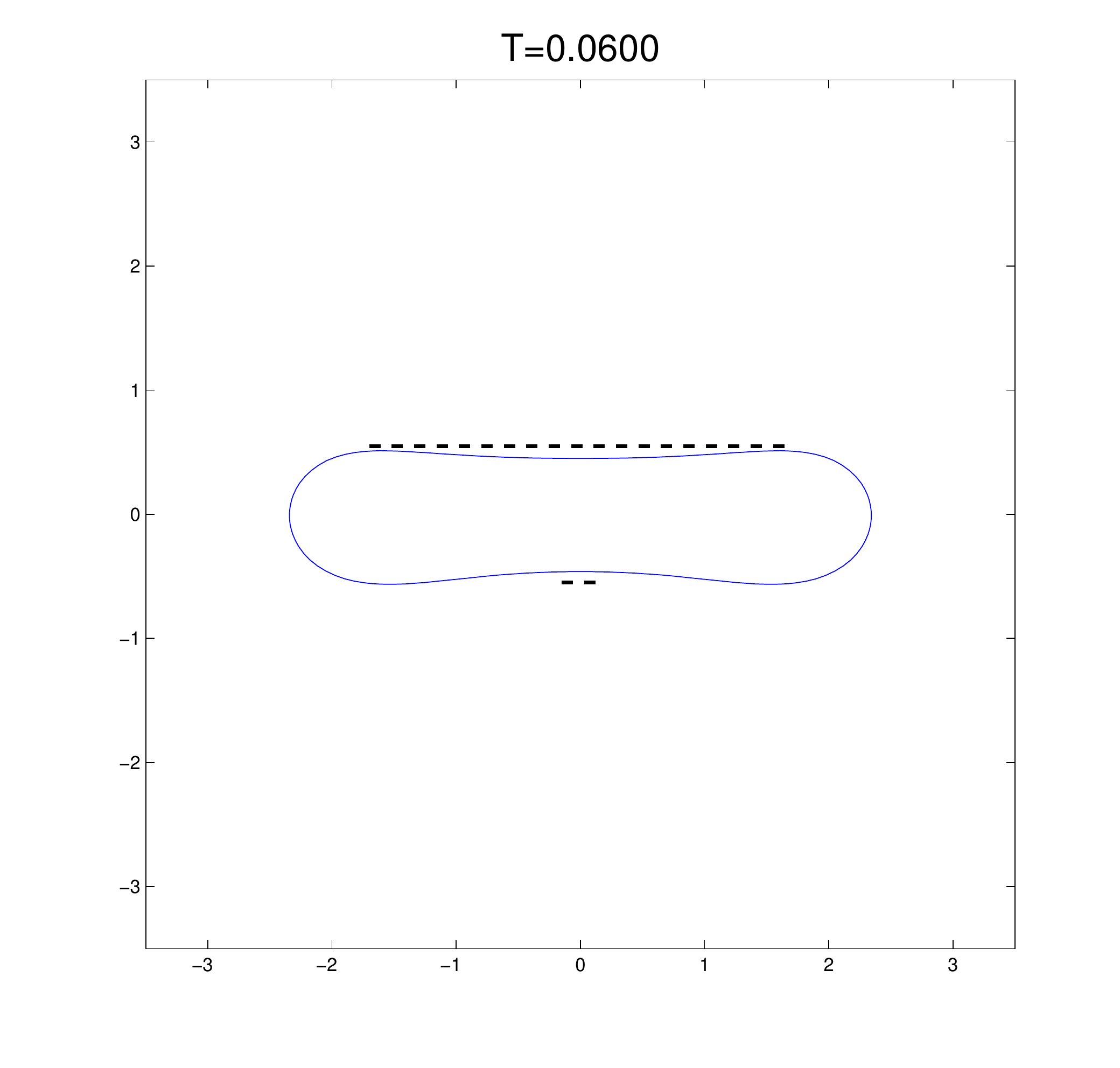}
}
\hspace{0mm}
\subfloat[]{
\includegraphics[width=0.5\textwidth]{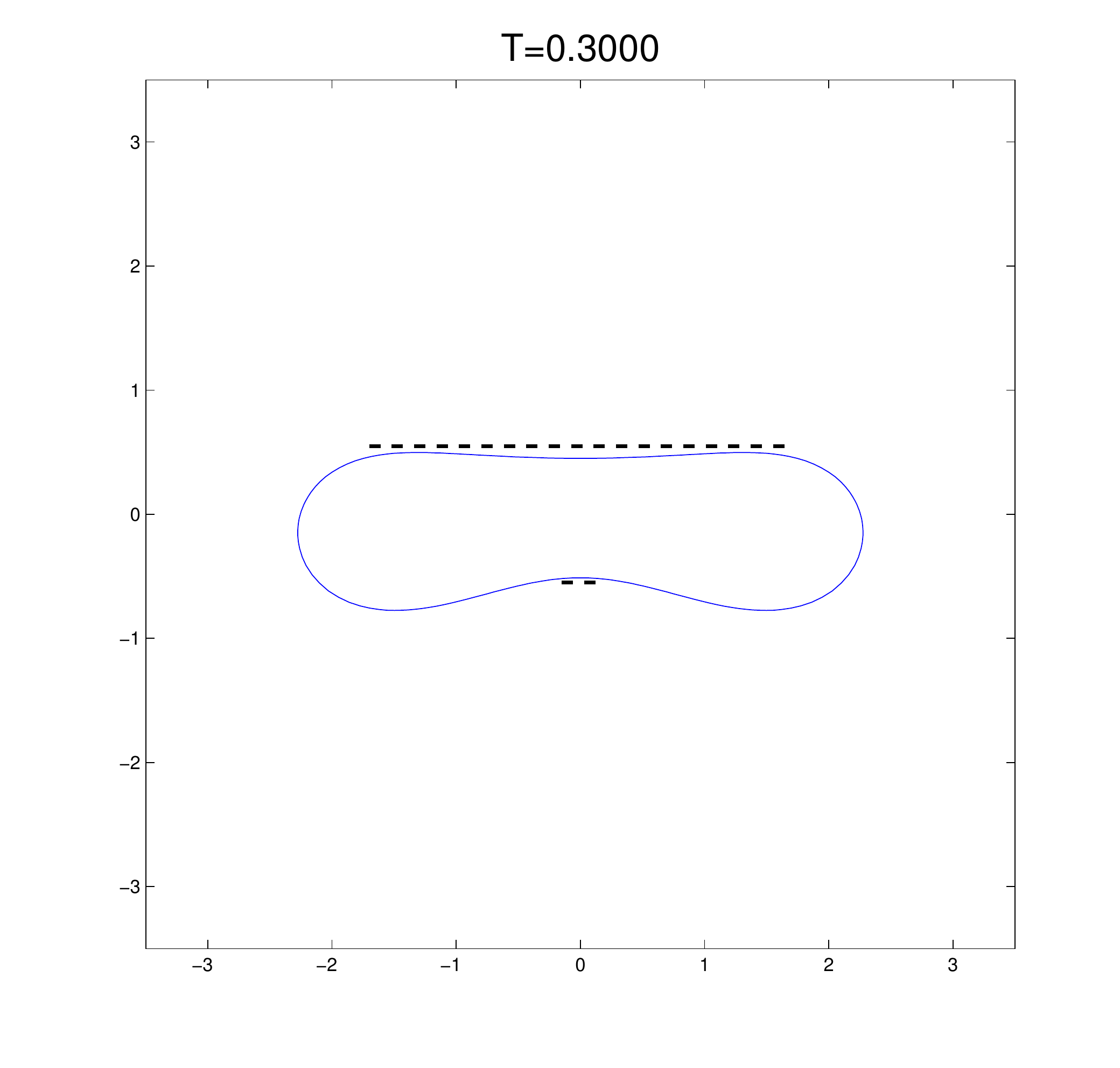}
}
\subfloat[]{
\includegraphics[width=0.5\textwidth]{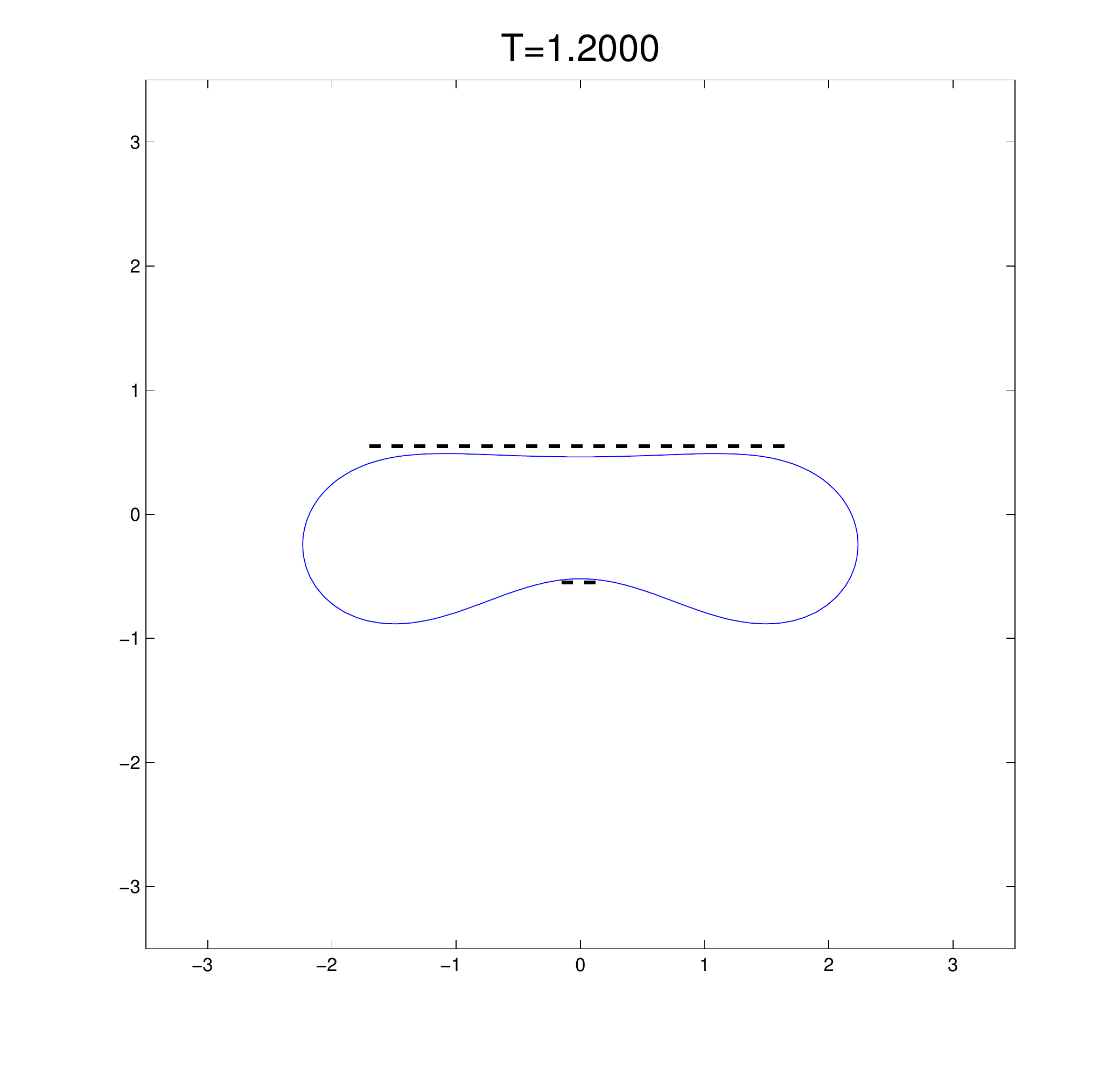}
}
\caption{Example 4 (to be continued)}
\end{figure}
\begin{figure}[h!]\ContinuedFloat
\centering
\subfloat{
\includegraphics[width=0.5\textwidth]{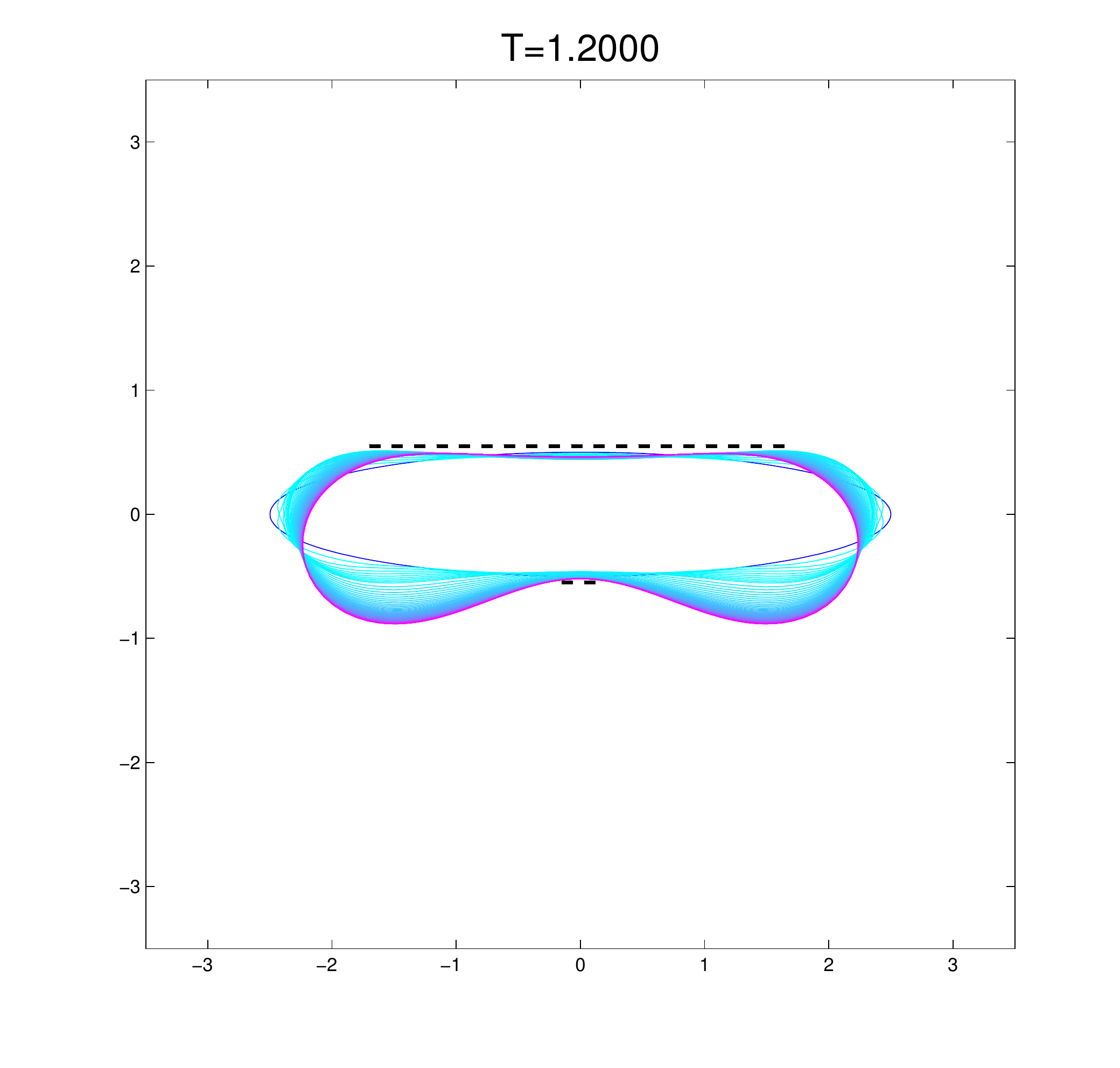}
}
\hspace{0mm}
\subfloat[Conservation of the Length $A(\Gamma)$]{
\includegraphics[width=0.5\textwidth]{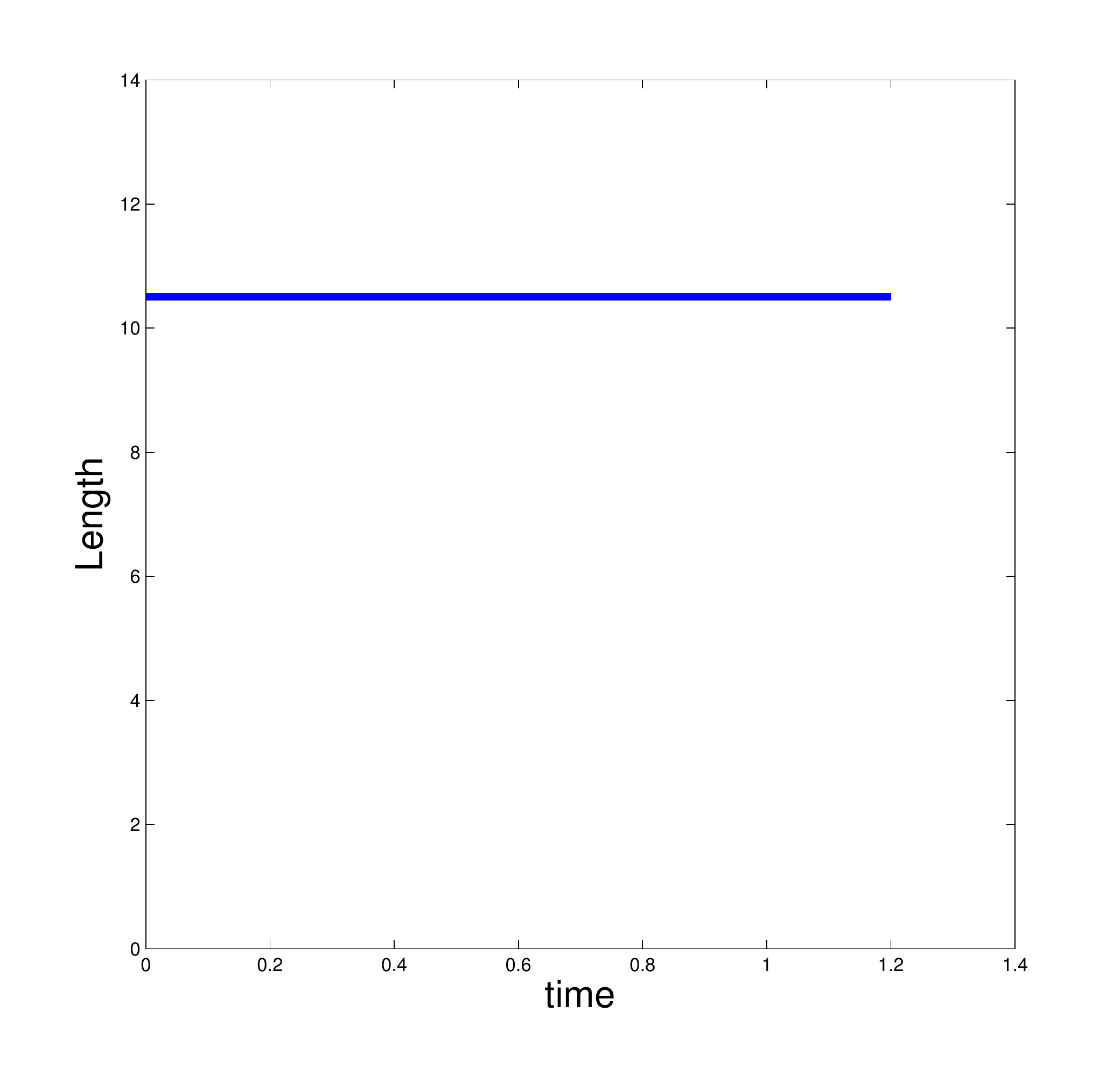}
}
\subfloat[Decrease of Energy $W(\Gamma)$]{
\includegraphics[width=0.5\textwidth]{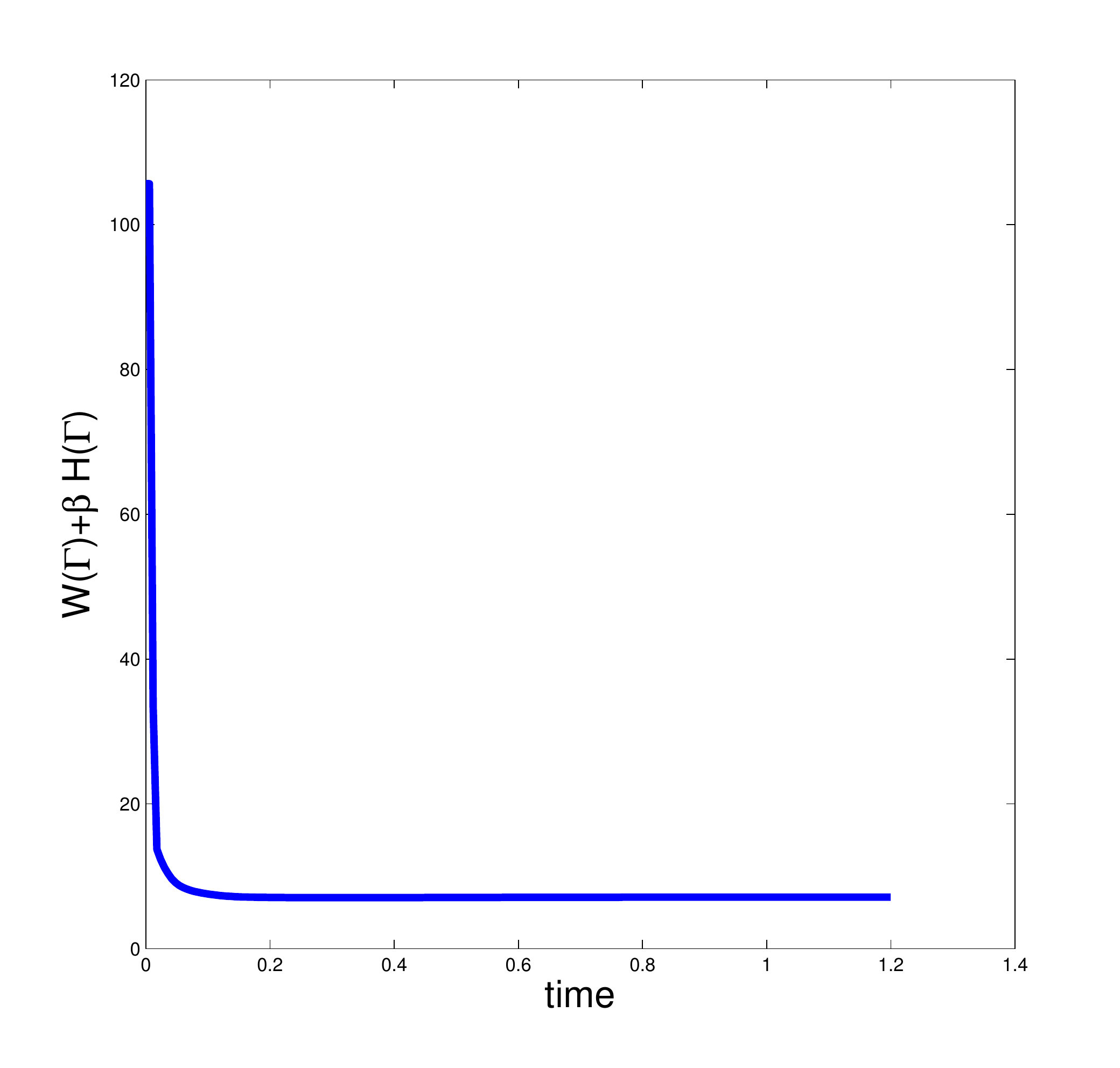}
}
\caption{(continue) Example 4}
\end{figure}

\clearpage
\subsubsection*{Example 5 - Model 2} % get out of the bound
This is an interesting example. When the barriers are put on one side of an ellipse, the ellipse will get rid of the bounds and then become a circle. This experiment is pretty realistic.
\begin{figure}[h!]
\centering
\subfloat[T=0]{
\includegraphics[width=0.5\textwidth]{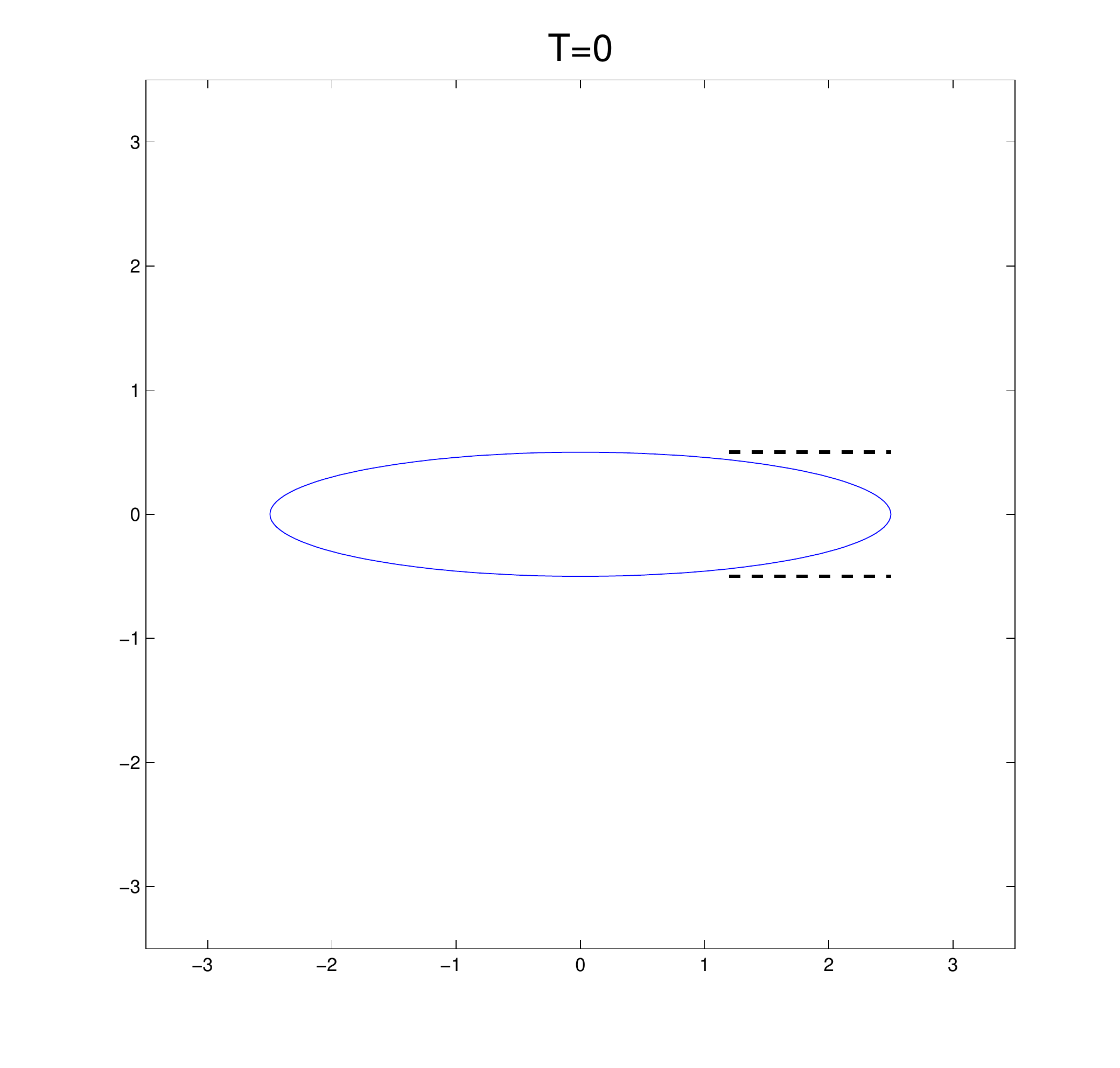}
}
\subfloat[]{
\includegraphics[width=0.5\textwidth]{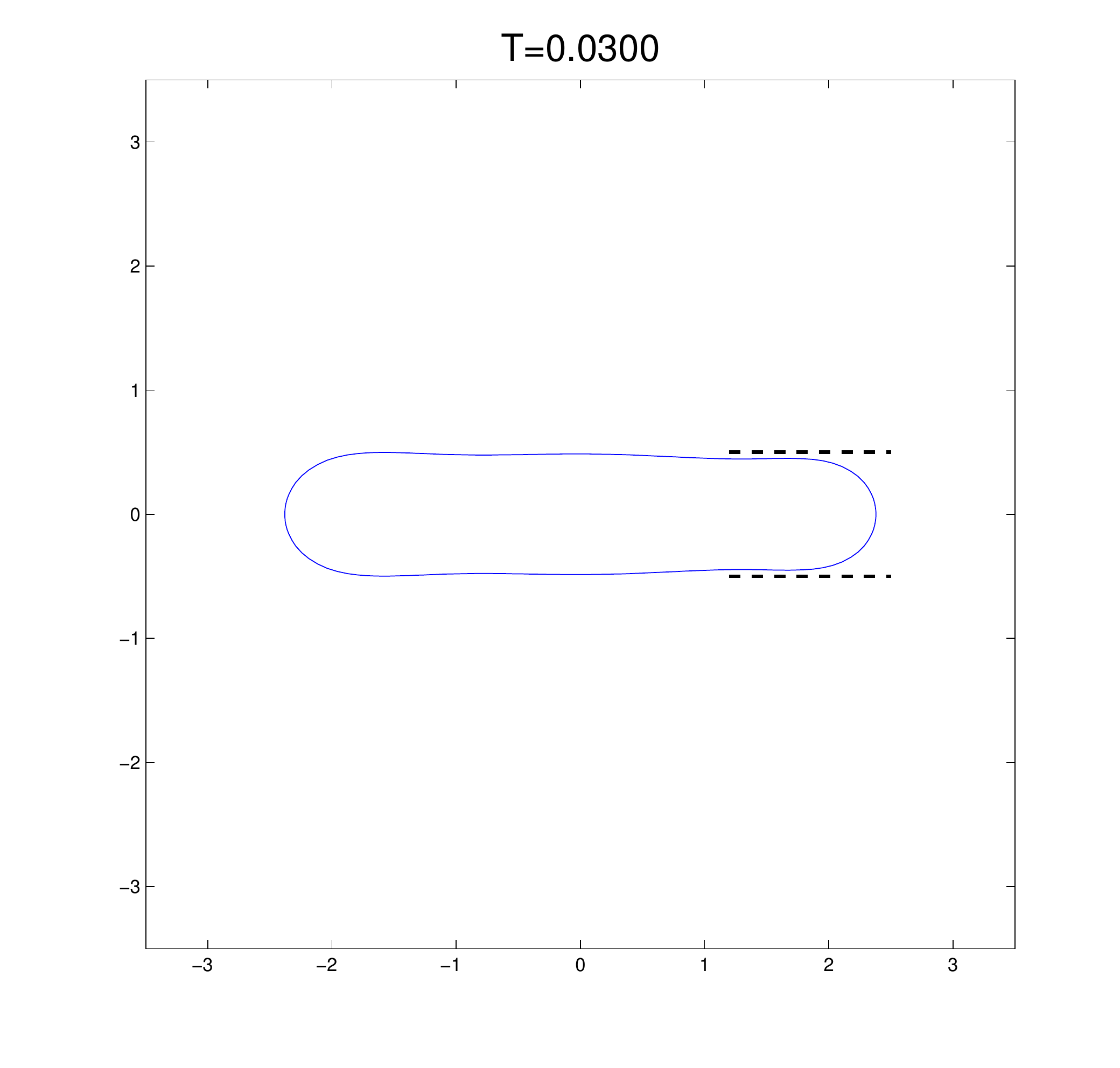}
}
\hspace{0mm}
\subfloat[]{
\includegraphics[width=0.5\textwidth]{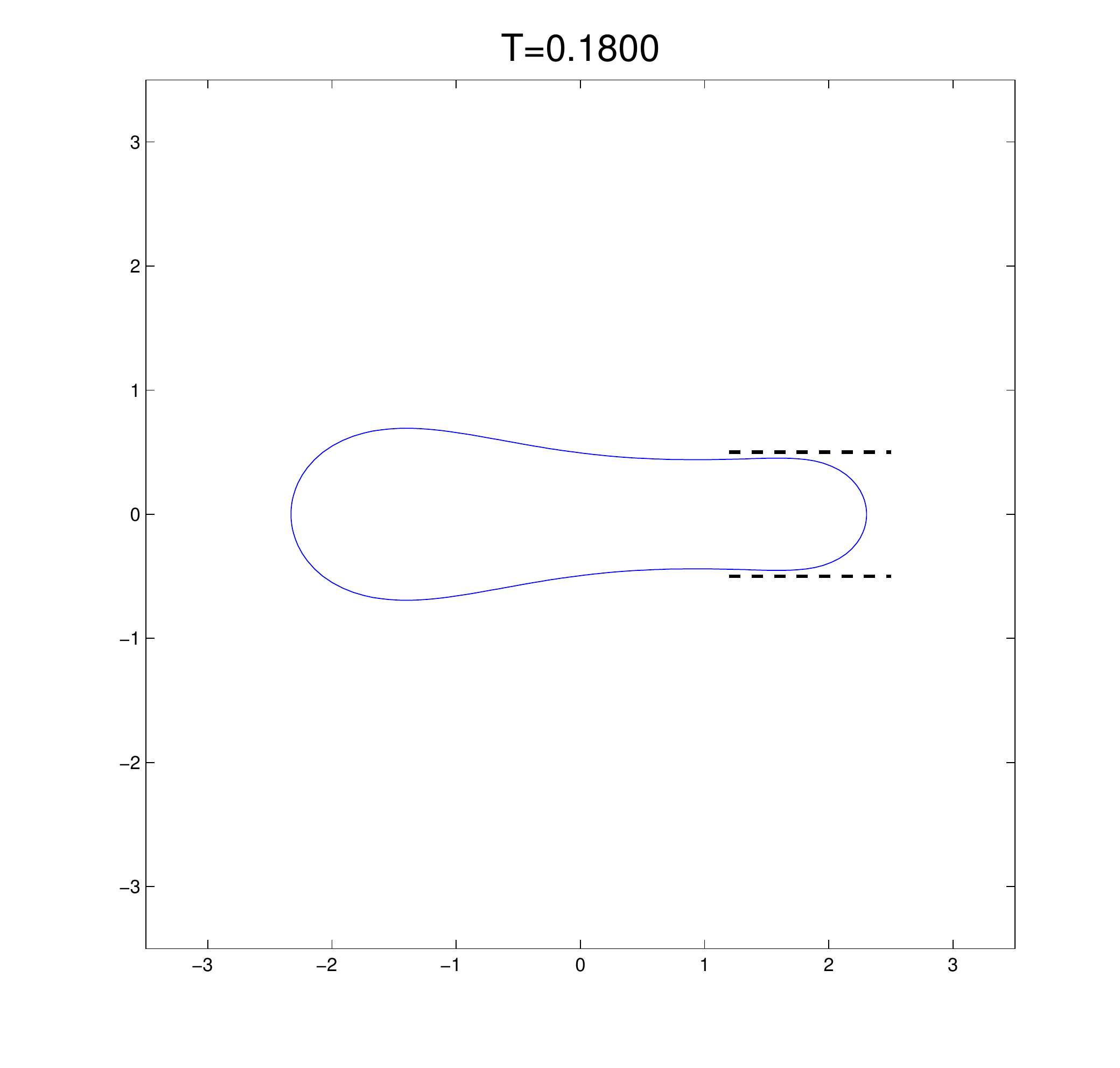}
}
\subfloat[]{
\includegraphics[width=0.5\textwidth]{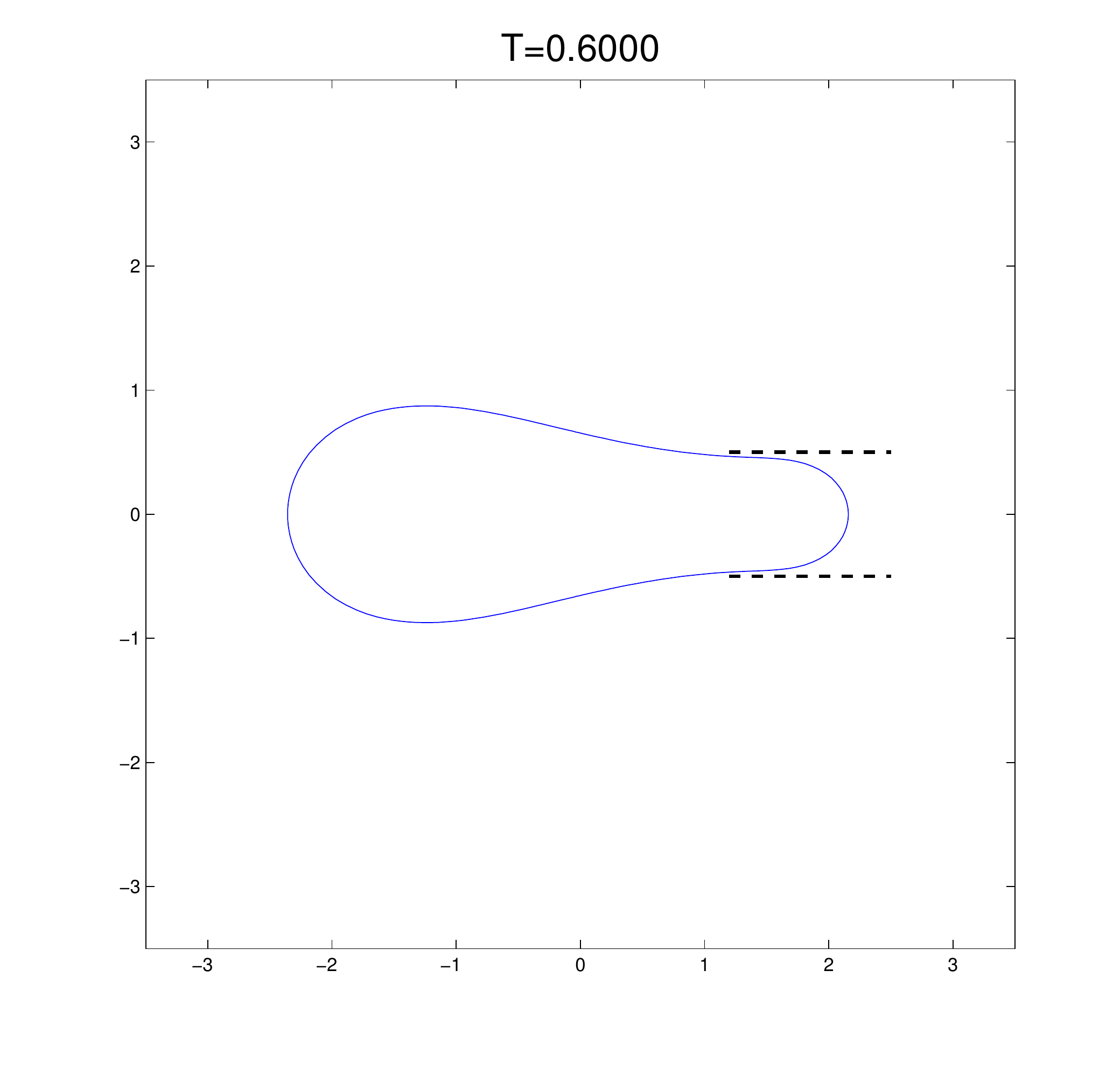}
}
\caption{Example 5 (to be continued)}	
\end{figure}

\begin{figure}[h!]\ContinuedFloat
\centering
\subfloat[]{
\includegraphics[width=0.43\textwidth]{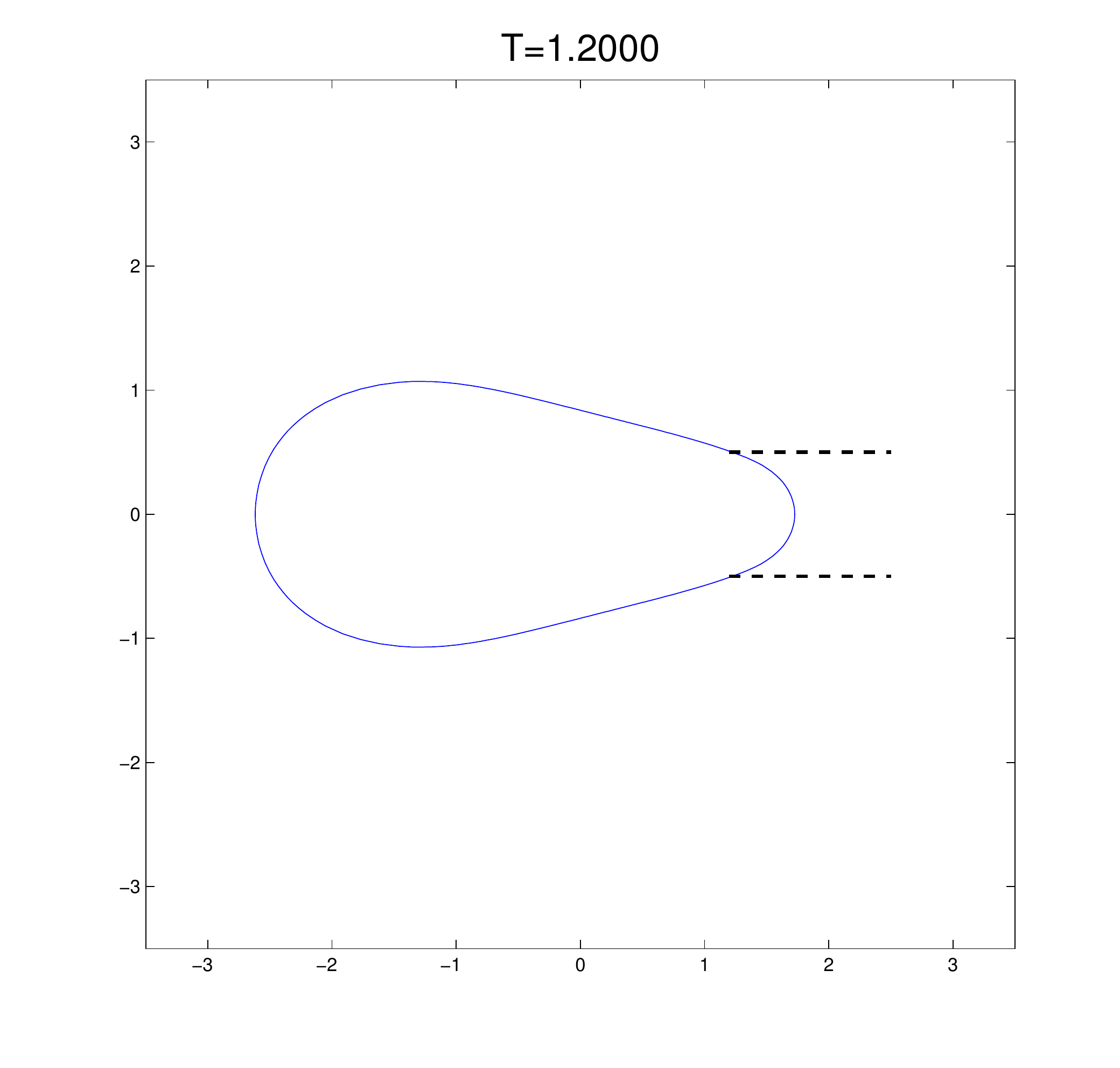}
}
\subfloat[]{
\includegraphics[width=0.43\textwidth]{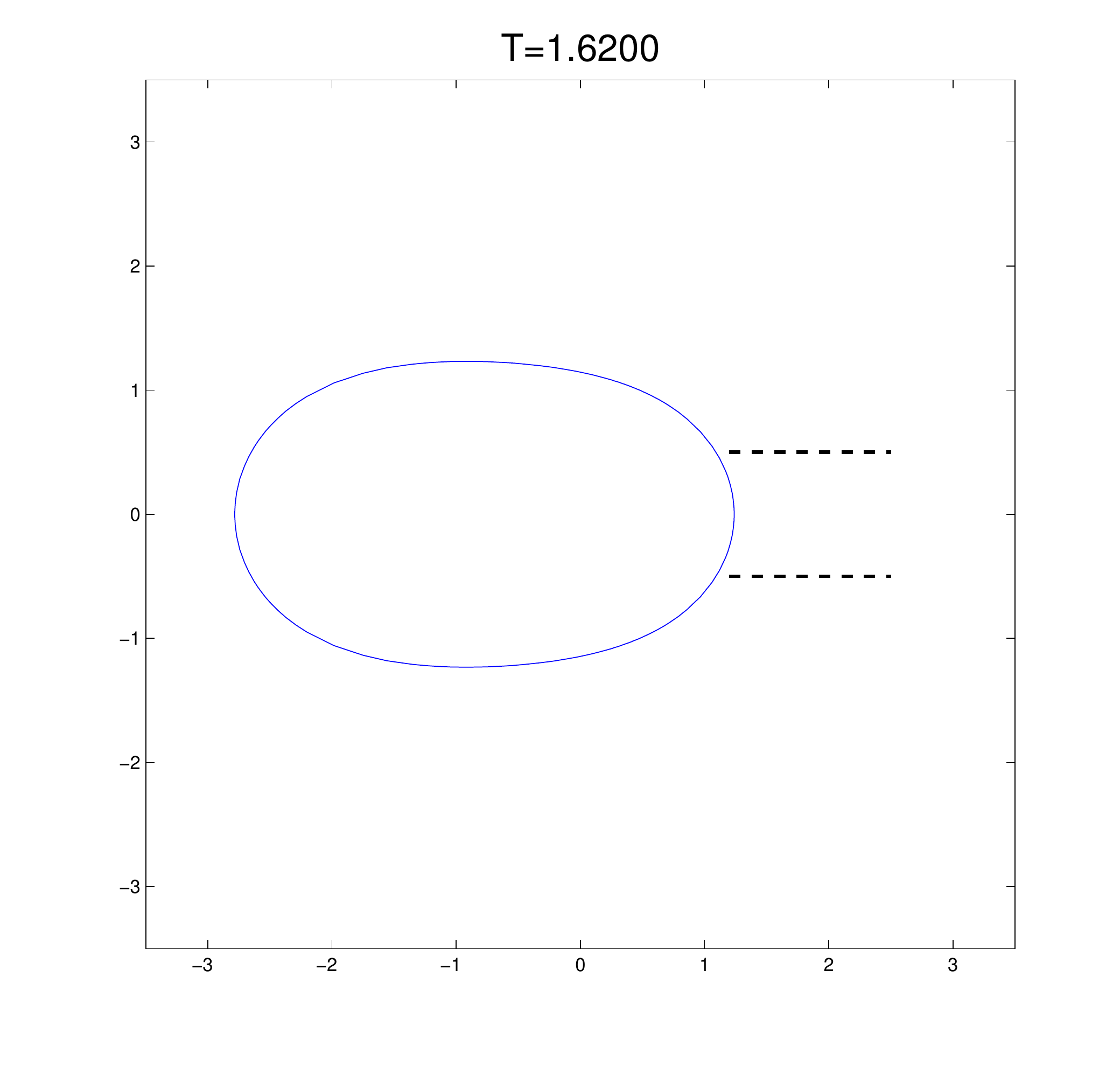}
}
\hspace{0mm}
\subfloat[]{
\includegraphics[width=0.43\textwidth]{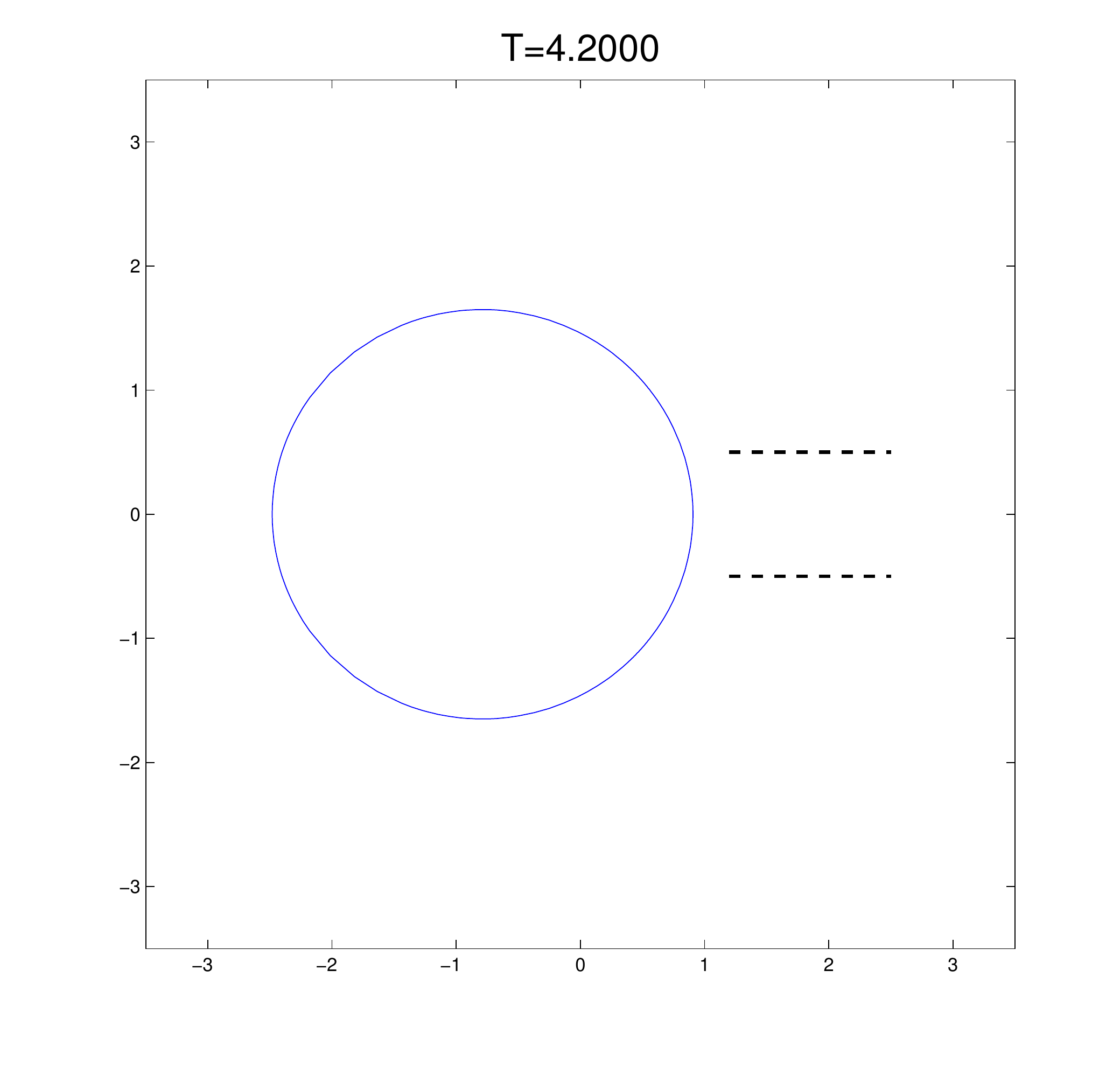}
}
\subfloat[]{
\includegraphics[width=0.43\textwidth]{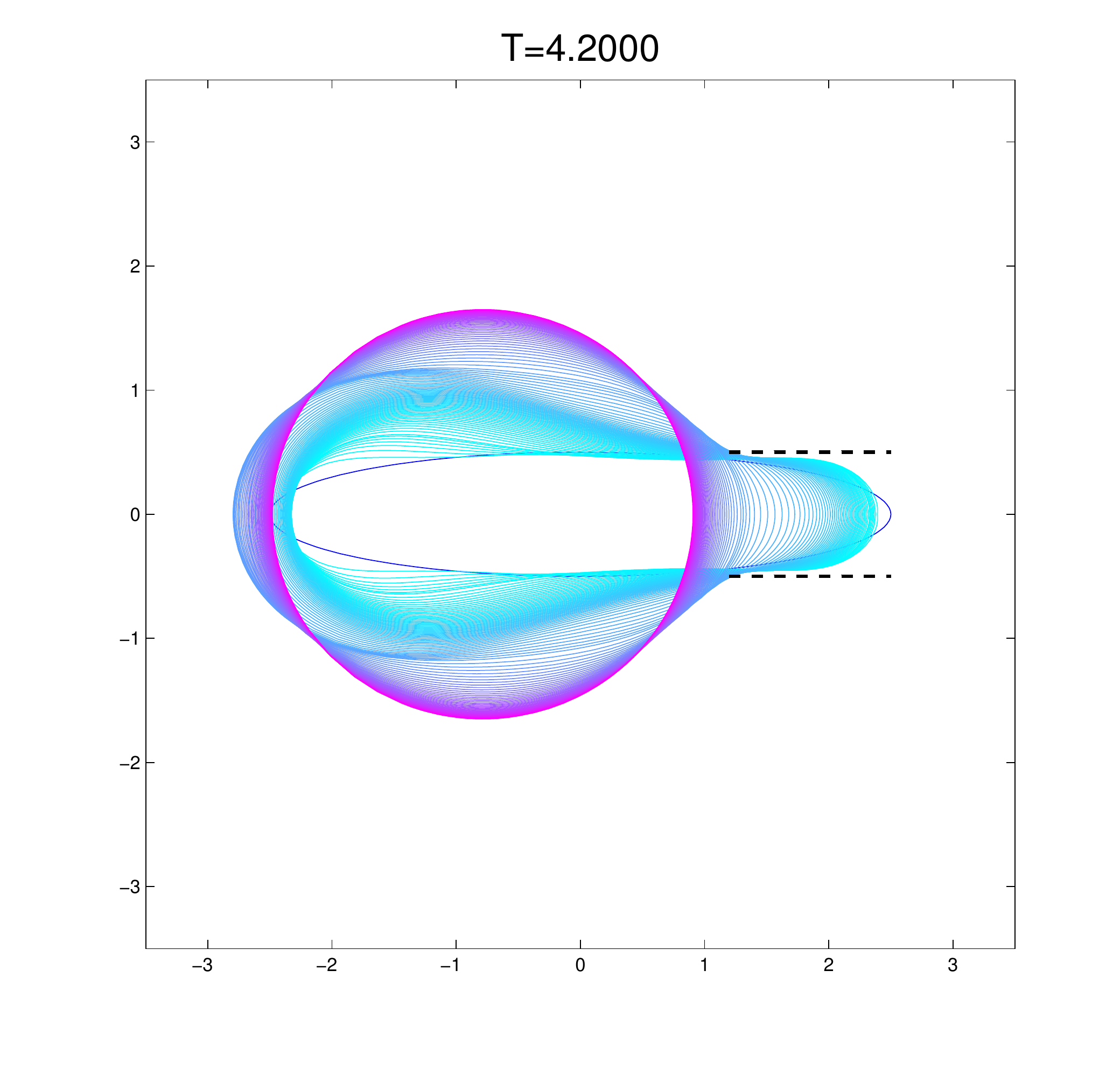}
}
\hspace{0mm}
\subfloat[Conservation of the Length]{
\includegraphics[width=0.43\textwidth]{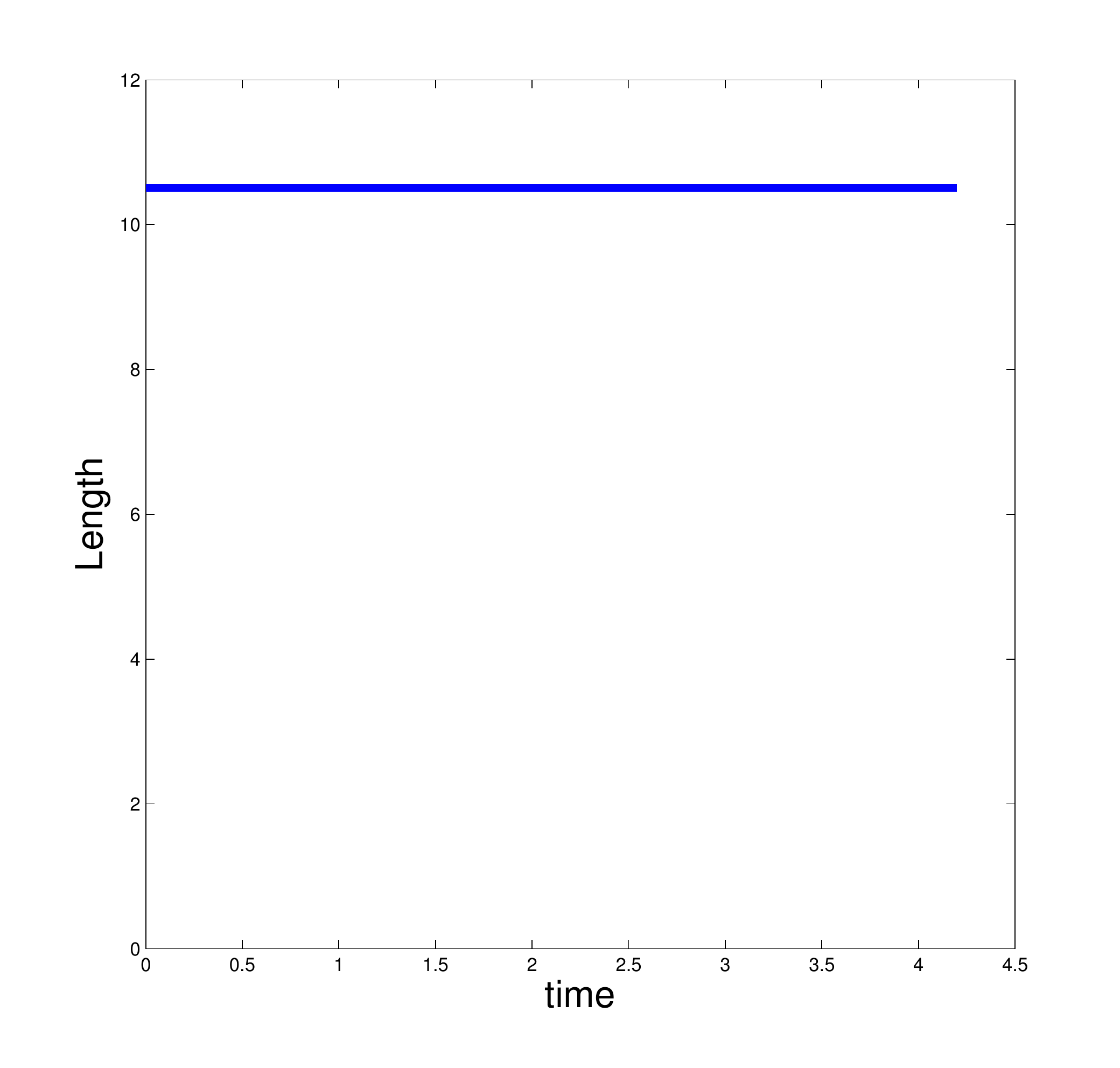}
}
\subfloat[Energy vs Time]{
\includegraphics[width=0.43\textwidth]{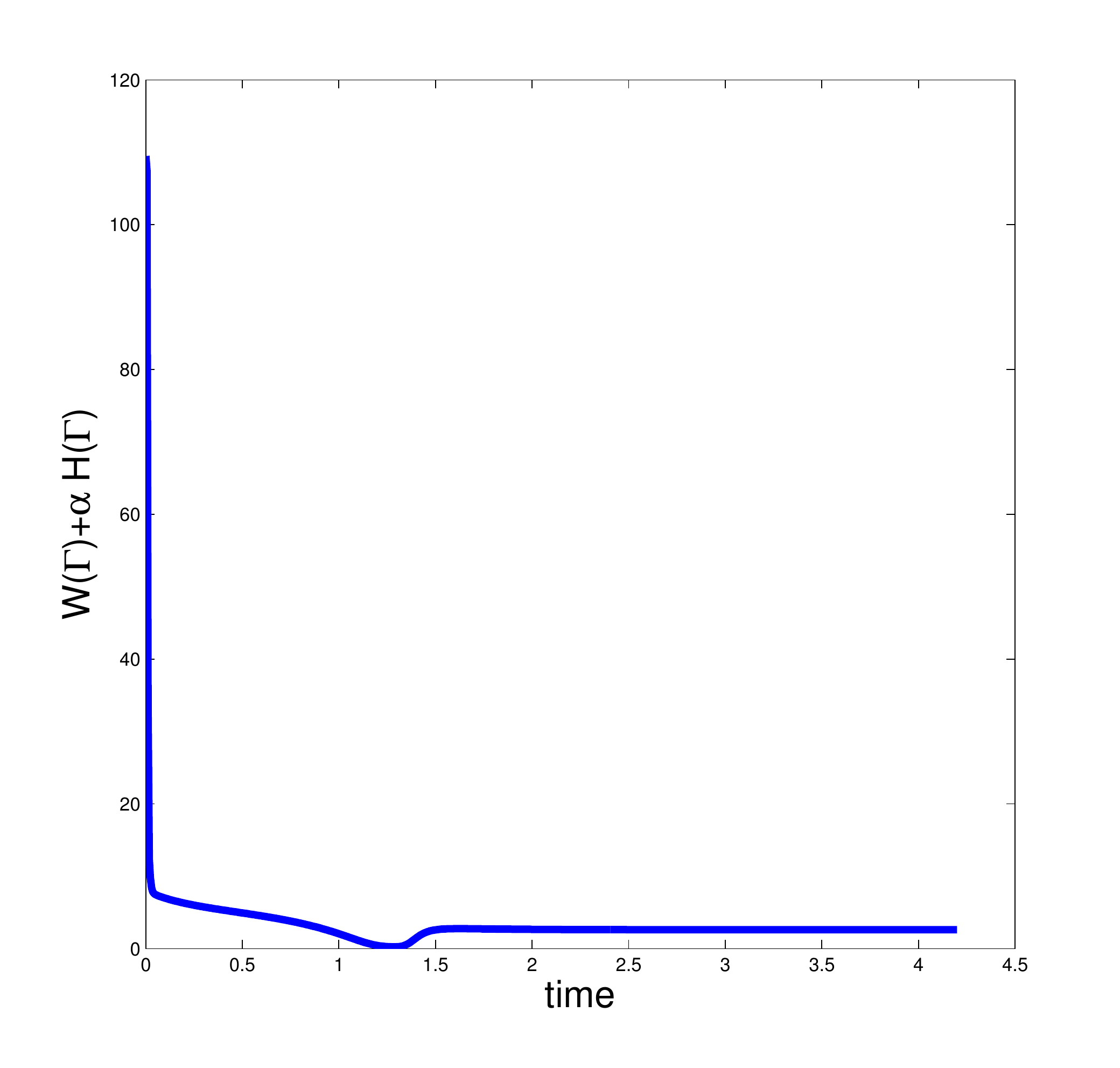}
}
\caption{(continue) Example 5}
\end{figure}

\clearpage
\subsubsection*{Example 6 - Model 2} % spinning
\begin{figure}[h!]
\centering
\subfloat[T=0]{
\includegraphics[width=0.5\textwidth]{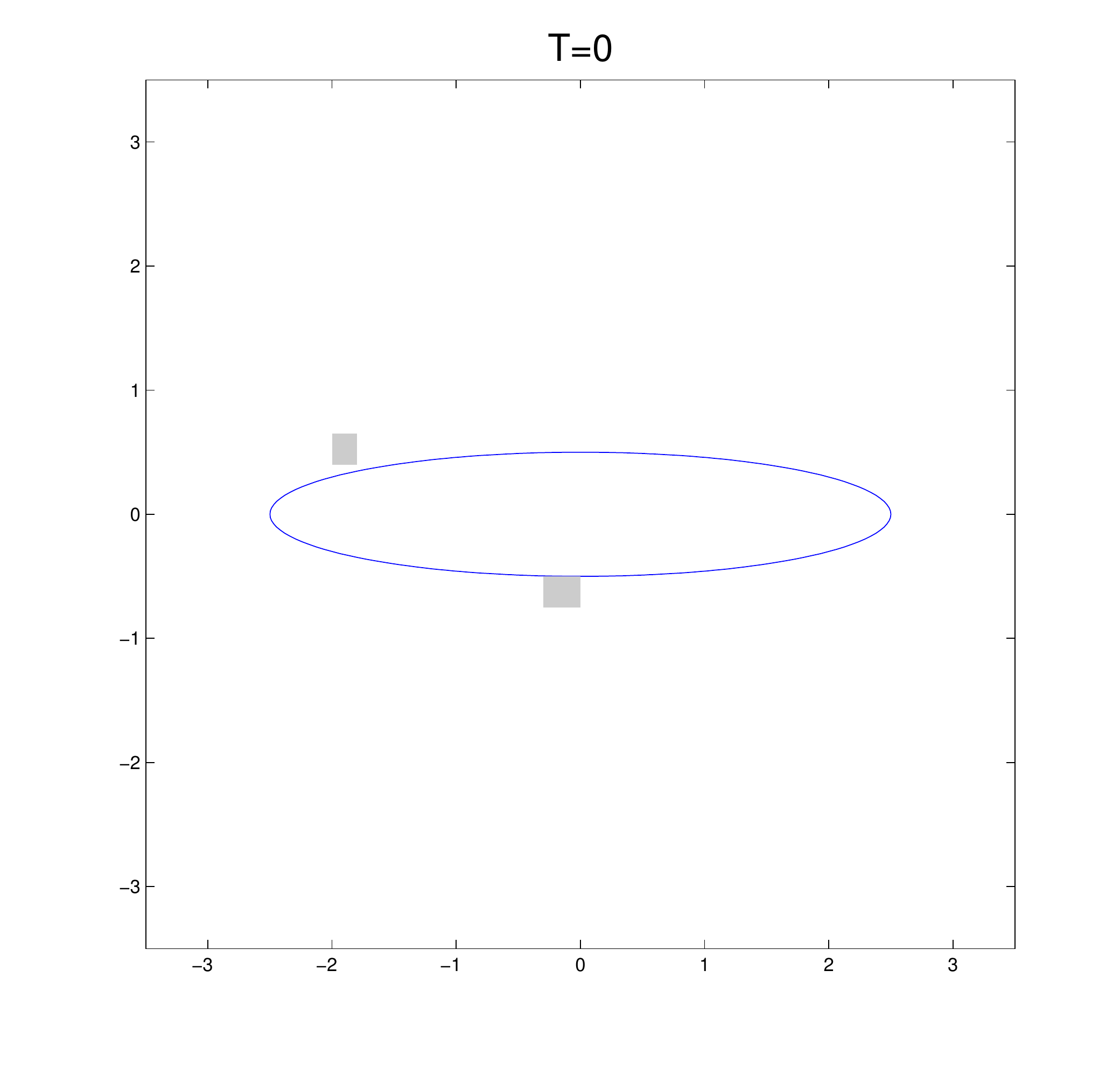}
}
\subfloat[]{
\includegraphics[width=0.5\textwidth]{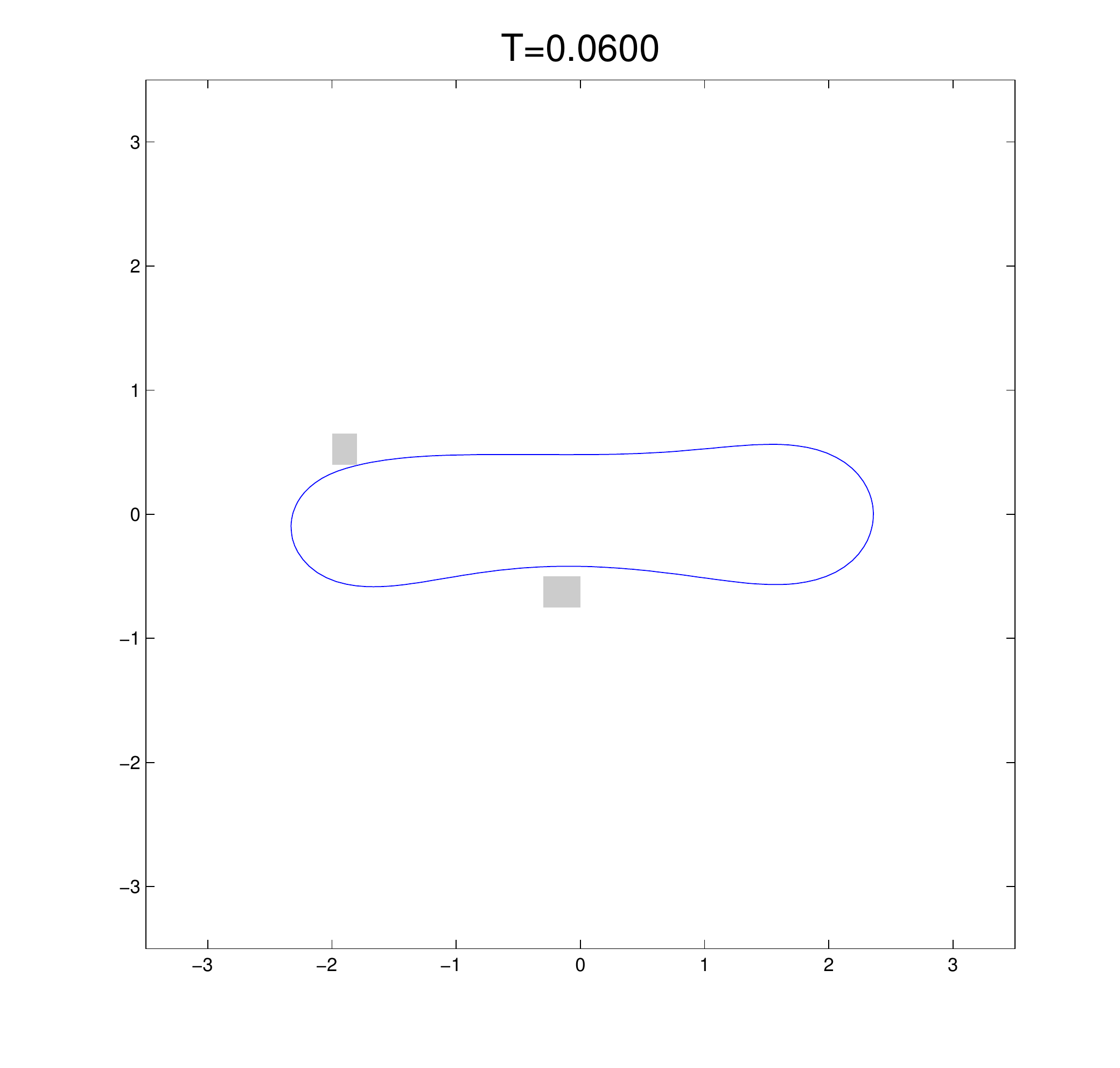}
}
\hspace{0mm}
\subfloat[]{
\includegraphics[width=0.5\textwidth]{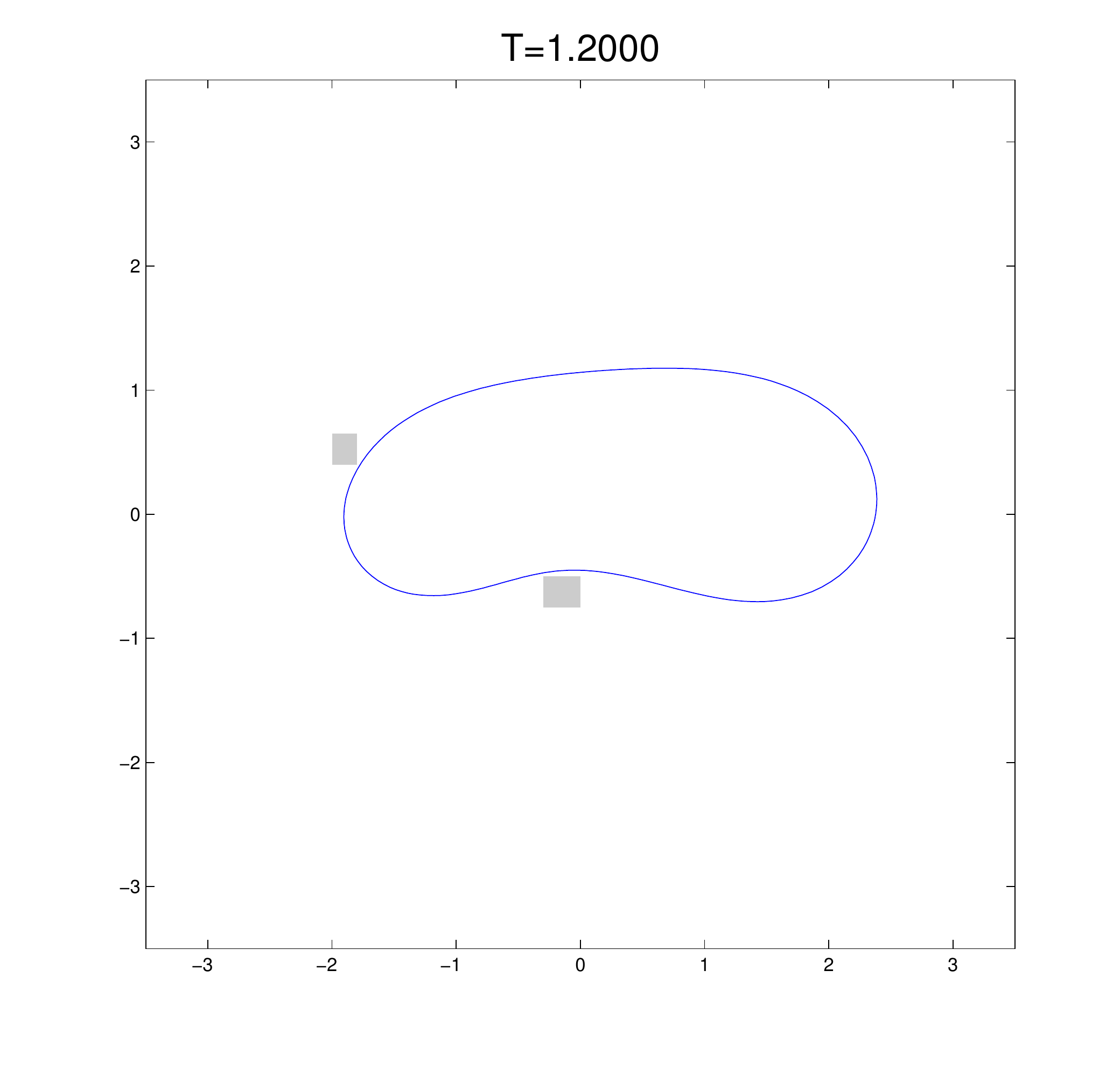}
}
\subfloat[]{
\includegraphics[width=0.5\textwidth]{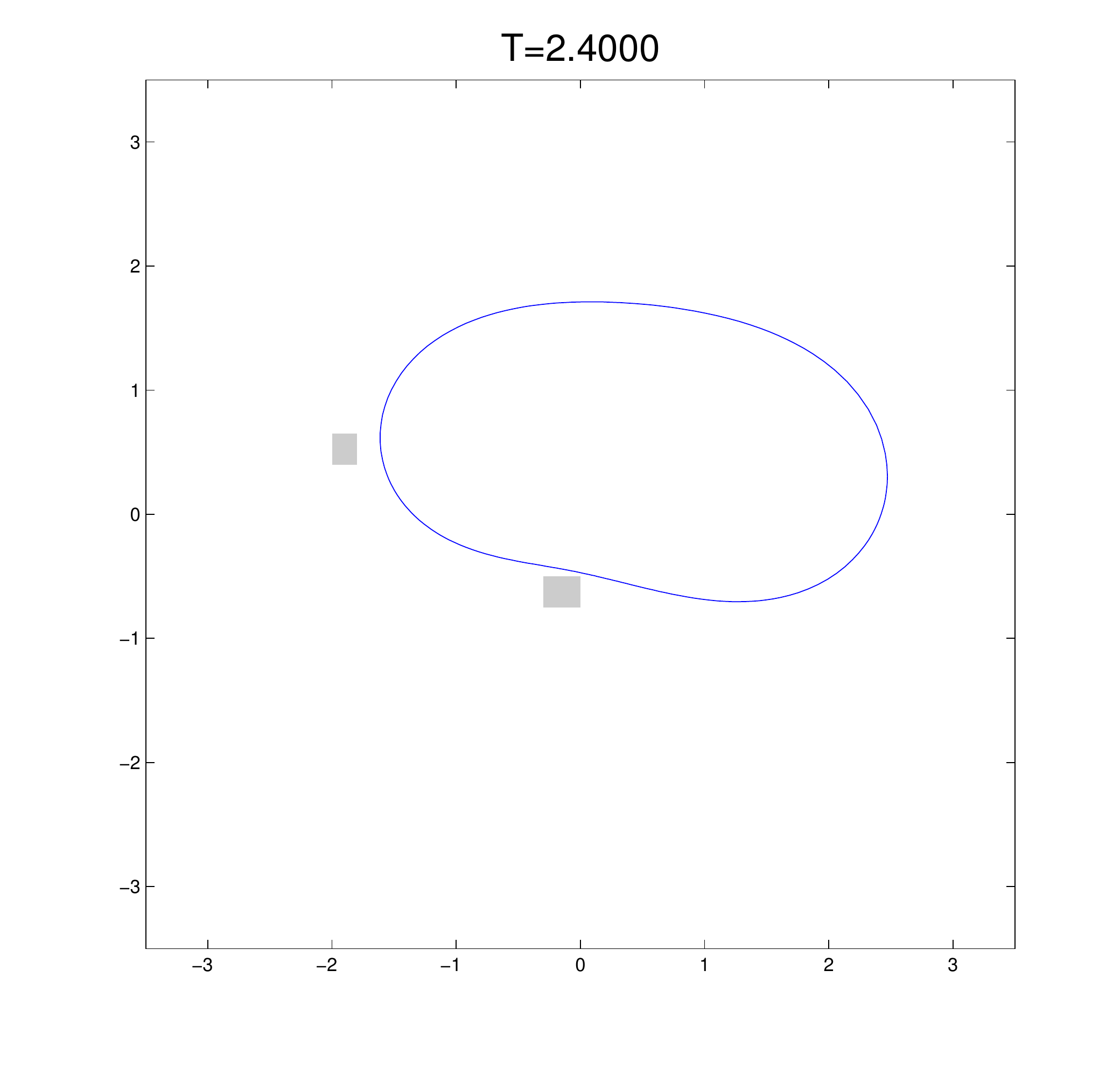}
}
\caption{Example 6 (to be continued)}	
\end{figure}
\begin{figure}[h!]\ContinuedFloat
\centering
\subfloat[]{
\includegraphics[width=0.45\textwidth]{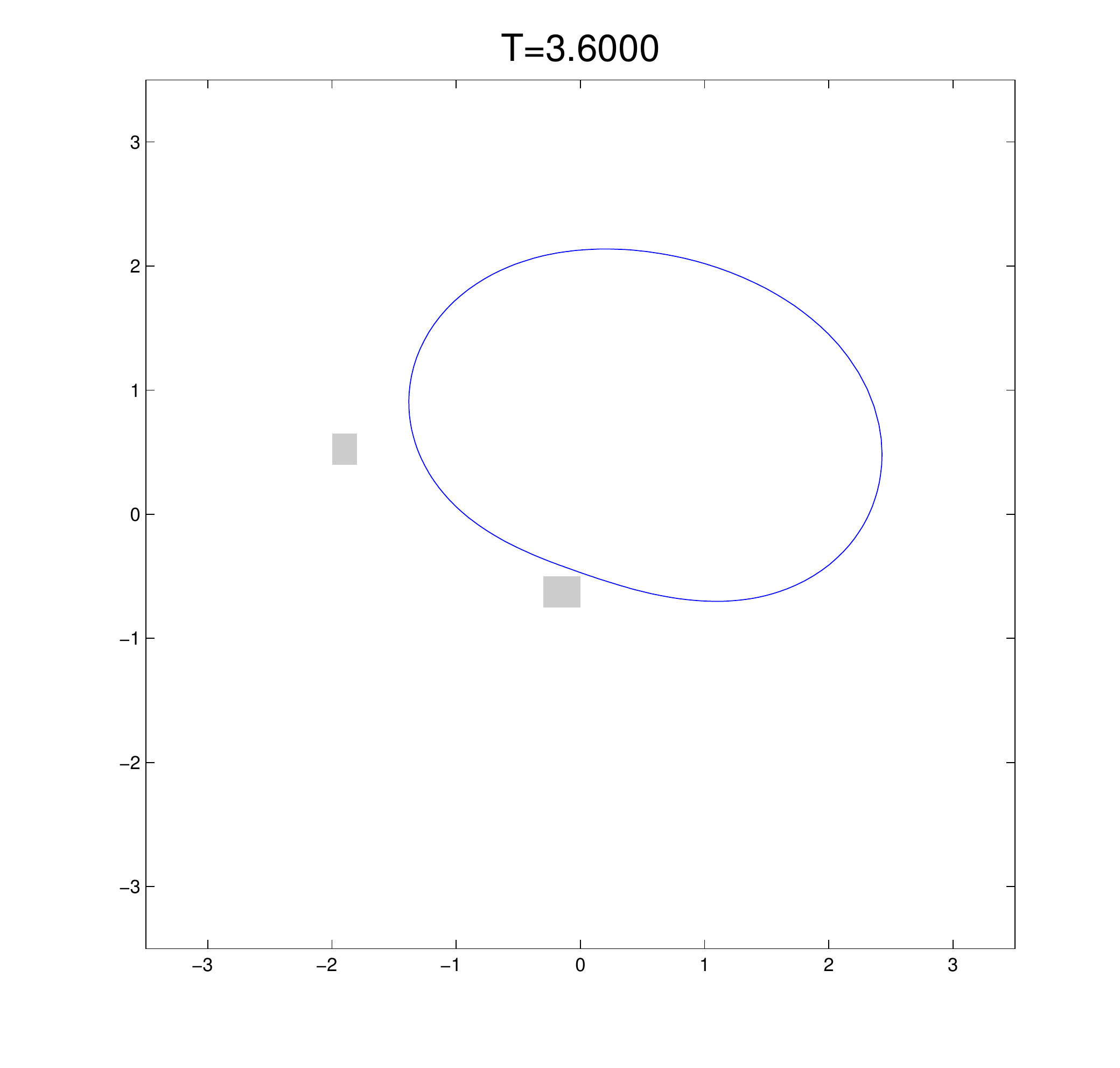}
}
\subfloat[]{
\includegraphics[width=0.45\textwidth]{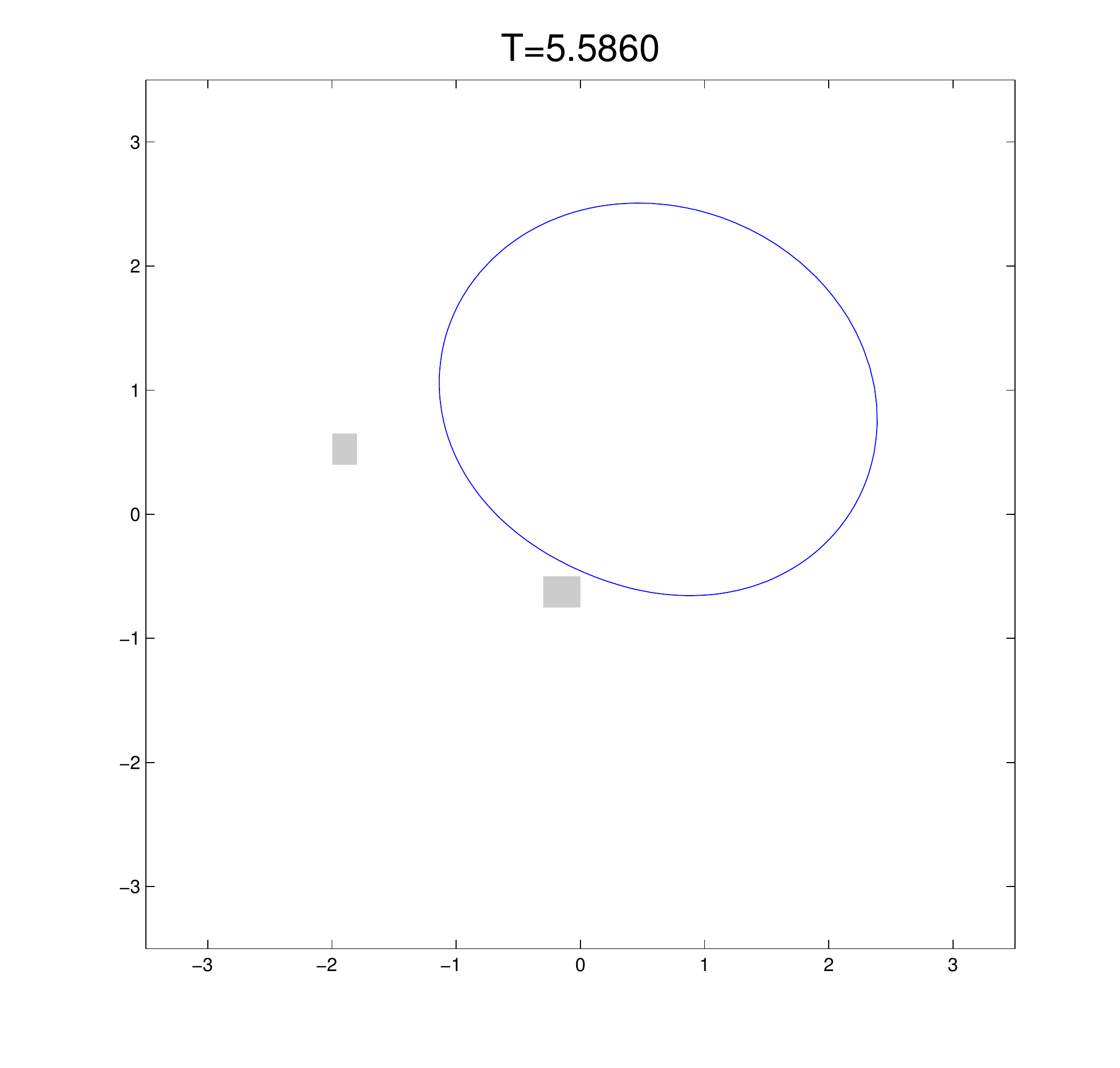}
}
\hspace{0mm}
\subfloat[]{
\includegraphics[width=0.45\textwidth]{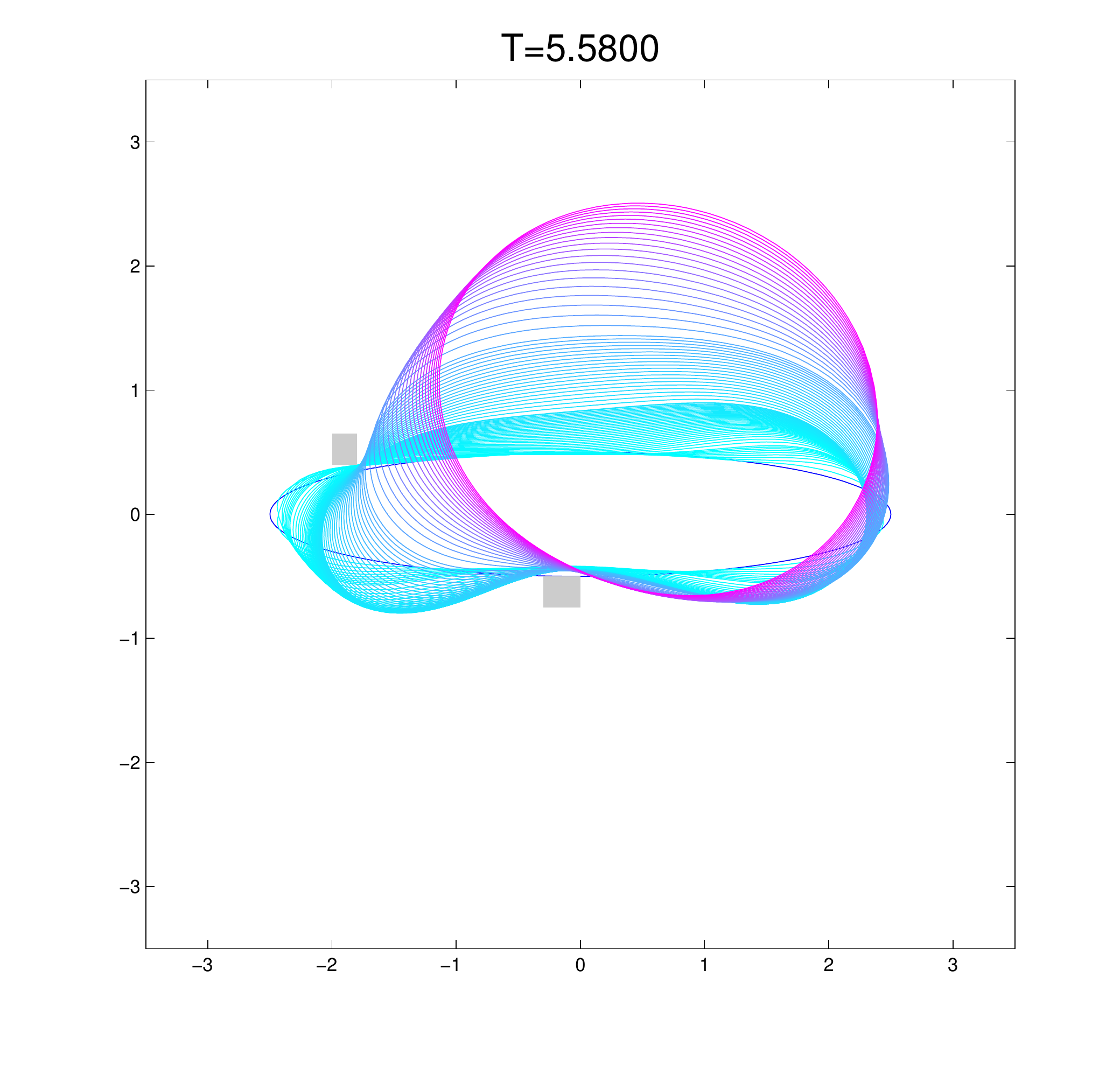}
}
\hspace{0mm}
\subfloat[Conservation of the Length]{
\includegraphics[width=0.45\textwidth]{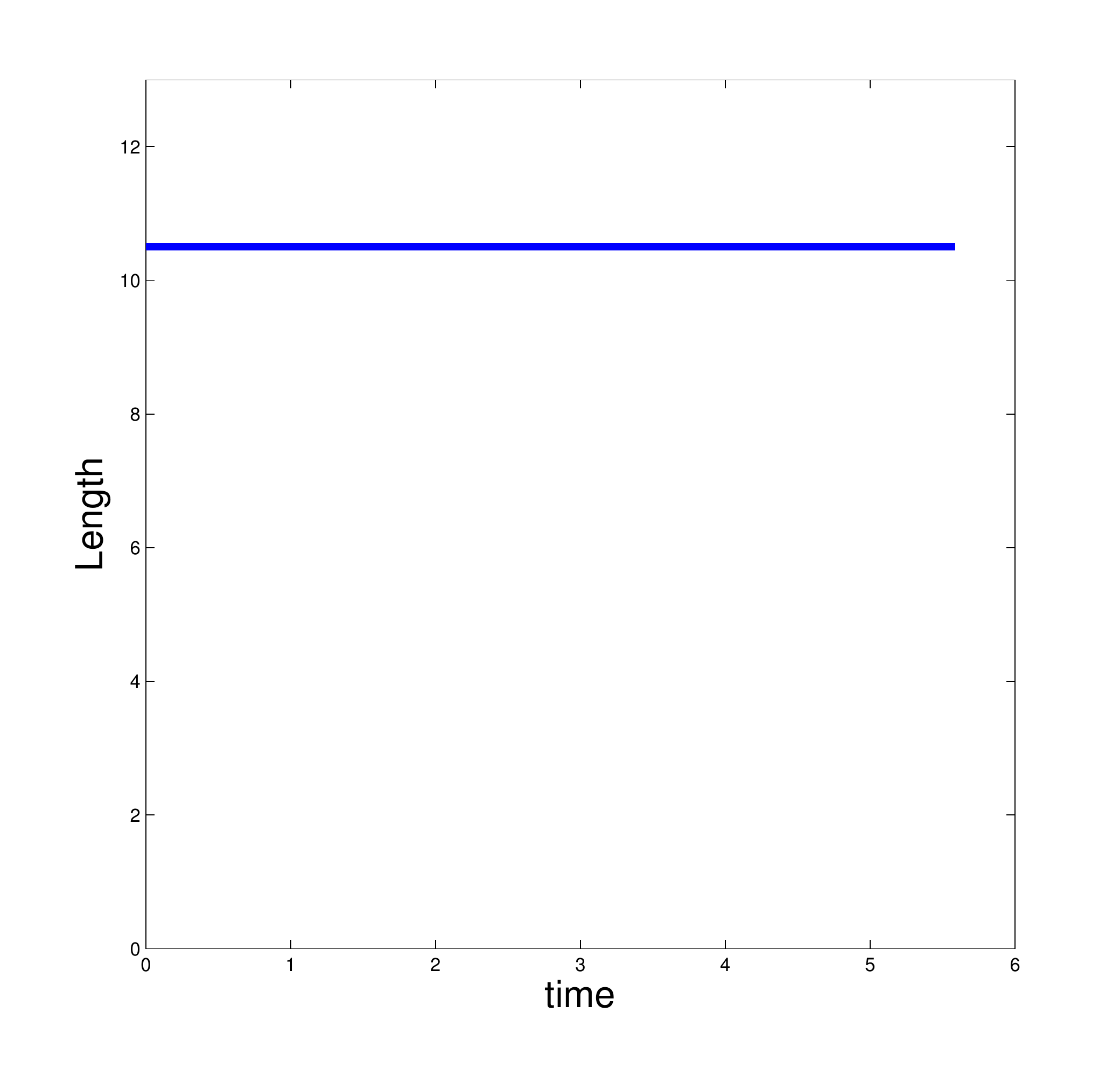}
}
\subfloat[Energy vs Time]{
\includegraphics[width=0.45\textwidth]{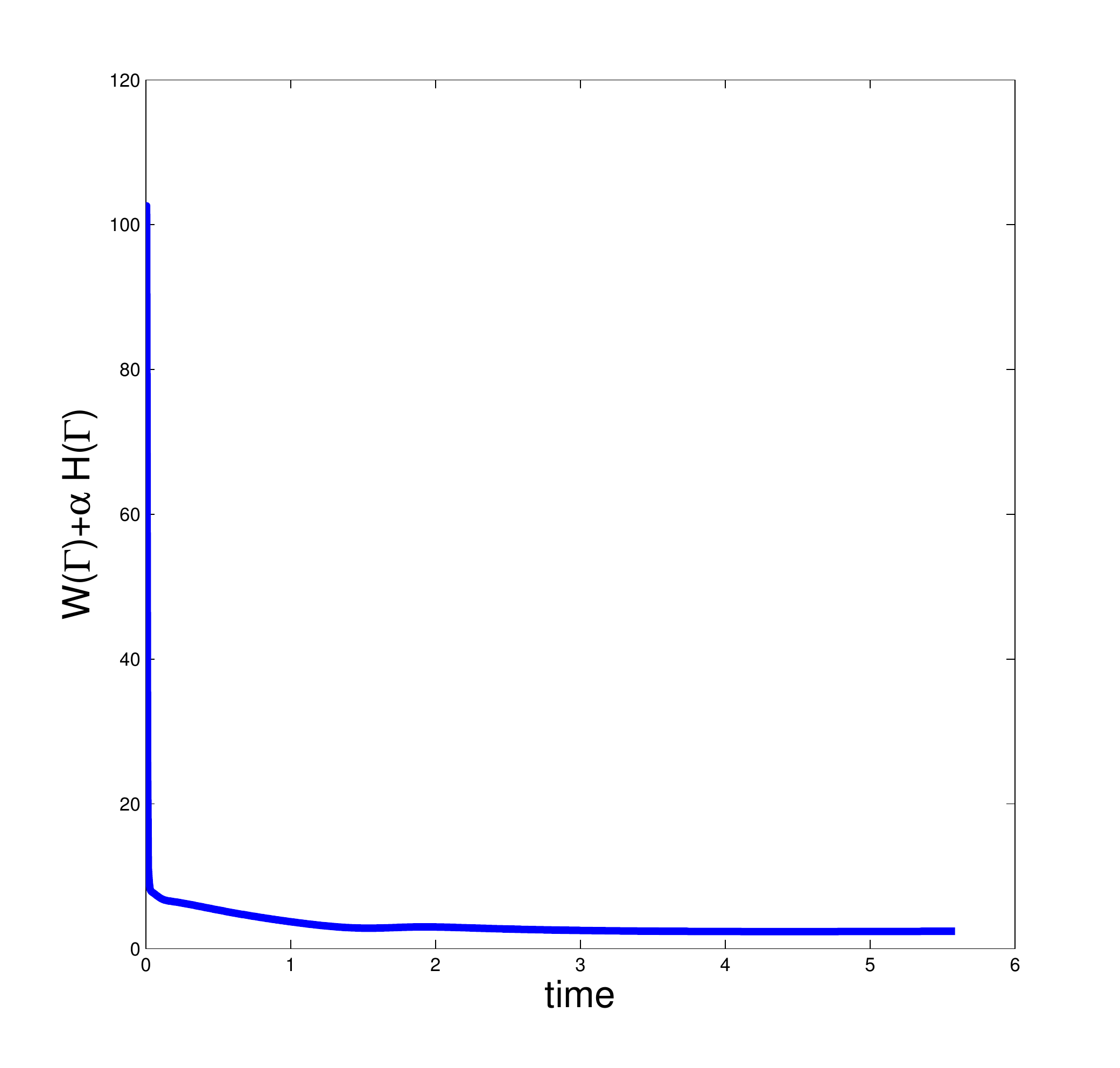}
}
\caption{(continue) Example 6}
\end{figure}

%\clearpage
%\subsubsection*{Example 7 - Model 3}
%Here we demonstrate an example which does not only contain the obstacles ${\bf B}$ but also contain the more than one objects $\Gamma_i$. Example 7 is a model consisting of a C shape, an ellipse and two bars.
%\begin{figure}[h!]
%
%\caption{Example 7}
%\label{fig:CE}	
%\end{figure}

The examples continue in next section. In Section 5.2, more examples of Model 2 and Model 3 are explained. In particular, those examples have their practical meanings applying to the shape evolution of Golgi stacks.

\clearpage
\section{Golgi Stacks Examples}

\subsection*{Detailed Biological Explanations}

\begin{figure}[h!]
\centering
\includegraphics[width=65mm]{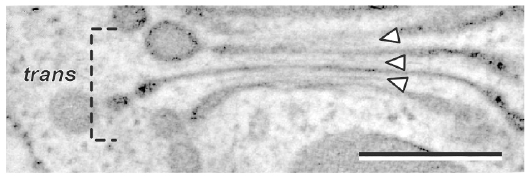}	
\caption{Tomographic slice image of trans Golgi cisternae in a
border cell and its inter-cisternal elements (arrowheads). Scale bar: 300 nm.}
\label{fig:ETthinning}
\end{figure}

There are different stages (cis-,med-,trans-) for the Golgi cisternaes. The morphological changes we aim to mimic are those for trans-cisternae in plant cells. Biologists have done some works on the ET analyses of Golgi, which observed that trans-cisternae are less thicker than interior cisternae \cite{thinning1,thinning2}. Prof. Kang (Life Science Department, CUHK) also provided an image of trans-cisternae in the border cell (Figure \ref{fig:ETthinning}), from which one can draw similar conclusion with the previous biological work.

To explain the observed morphological properties, we model the evolution of the Golgi membrane $\Gamma$ in $\mathbb{R}^2$ for simplicity. The idea of our modeling is roughly stated in Chapter 2. More specifically, suppose that the shapes of the cisternal membrane are determined to minimize the elastic energy of the membrane with two barriers lying above and below the cisternae. These barriers are used to model some inter-cisternal elements, so for a single cisternae they just lie above and below. This is how we place the obstacles ${\bf B}$. Besides, based on the condition that cisternal assembly is completed in the med- stage of Golgi, the surface area of the membrane only grows in the med- stage, but for trans-cisternae the number of molecules which consist the membrane is fixed. Hence, the constraint $A(\Gamma) = A(\Gamma_0)$ applies. Lastly, the maturation process of Golgi in the plant border cells involves lots of synthesis of the biological substance. Naturally, the luminal volume of trans-cisternae should increase, though biologists did not measure it specifically. We do not consider a condition on the volume in our math model, because there is too much uncertainty about the measurement of the increase. However, when the elastic energy decrease, the volume usually increase naturally.
\subsubsection*{Numerical Examples for Single Cisternae}
With these ideas, we use Model 2 to mimic the evolution of a single trans-cisternae with initial shape $\Gamma_0$ shown in Figure \ref{fig:single1}(a). The first set of numerical results are demonstrated in Figure \ref{fig:single1}.

\begin{figure}[h!]
\centering
\subfloat[T=0]{
\includegraphics[width=65mm]{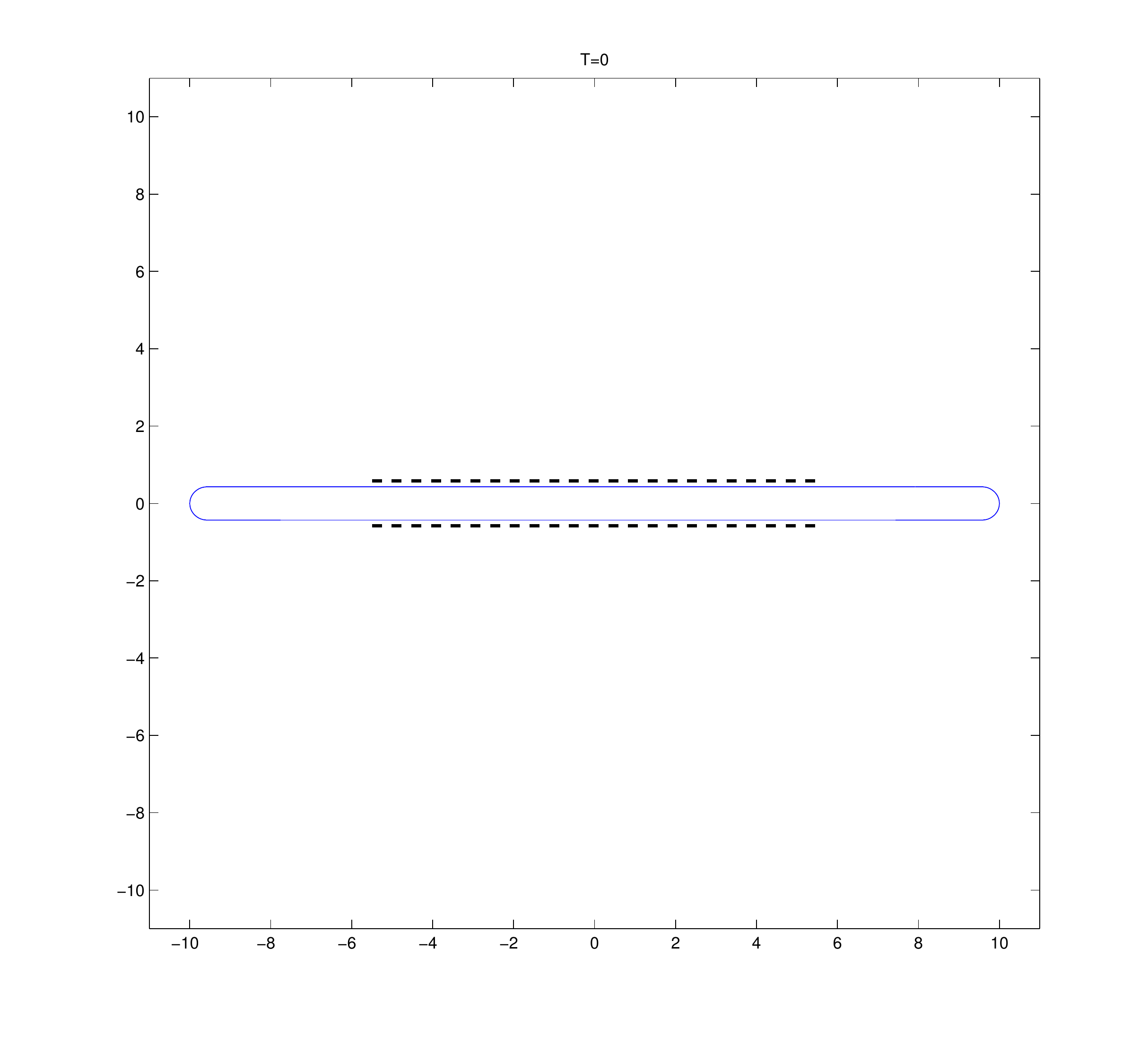}
}
\subfloat[T=0.15]{
\includegraphics[width=65mm]{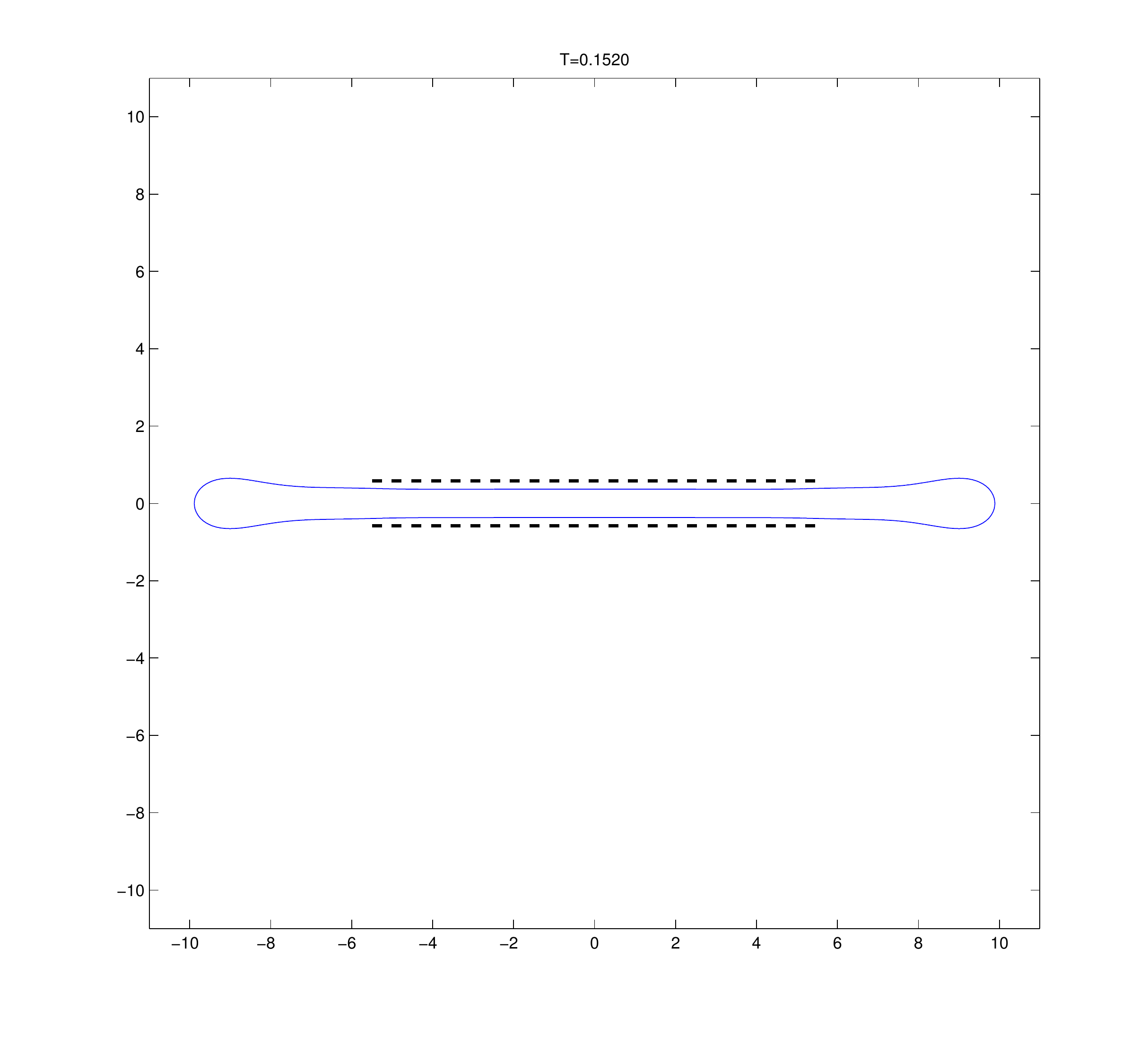}
}
\hspace{0mm}
\subfloat[T=0.79]{
\includegraphics[width=65mm]{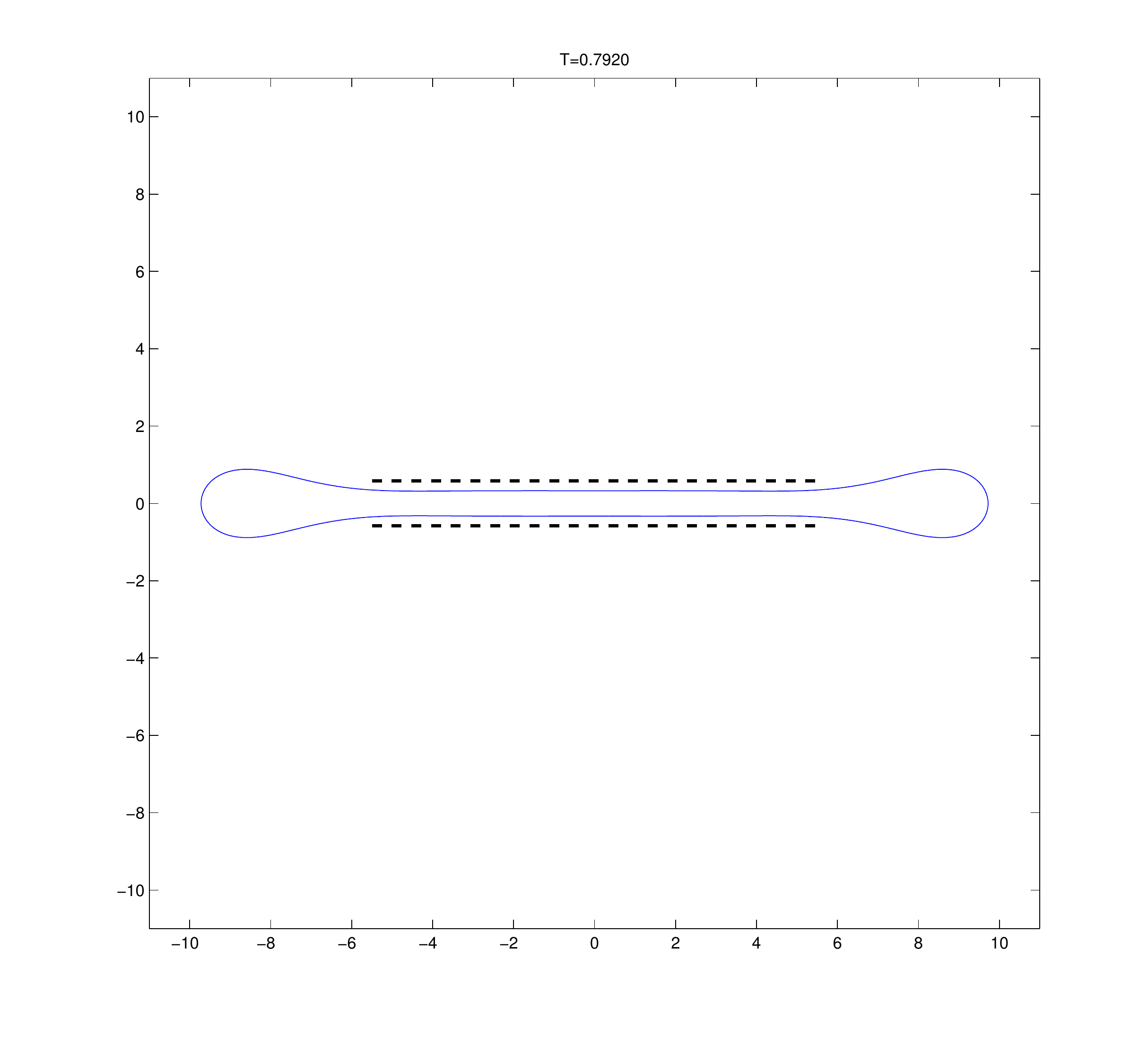}
}
\subfloat[T=1.59]{
\includegraphics[width=65mm]{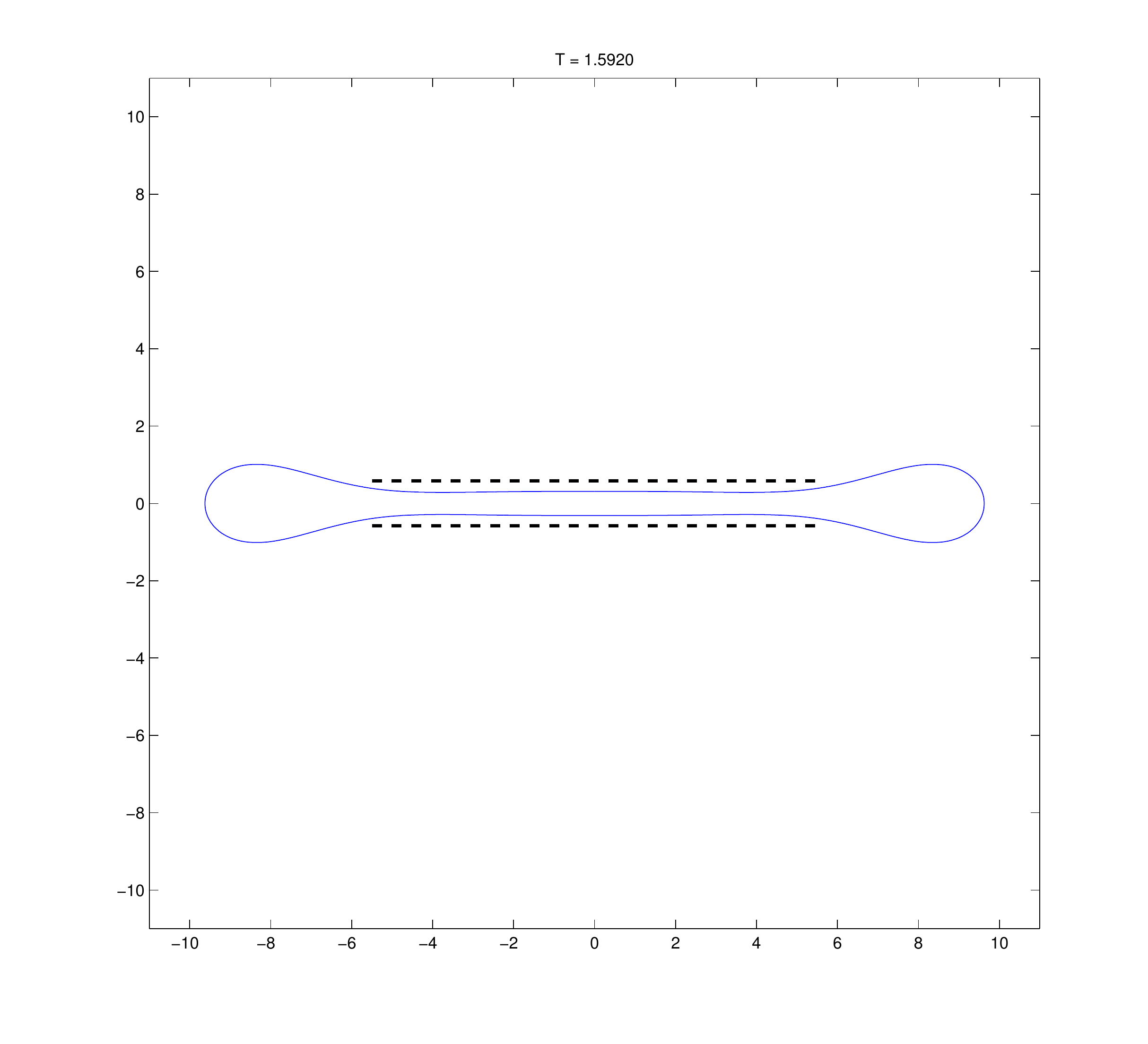}
}
\end{figure}
\begin{figure}[h!]\ContinuedFloat
\centering
\subfloat[T=2.39]{
\includegraphics[width=65mm]{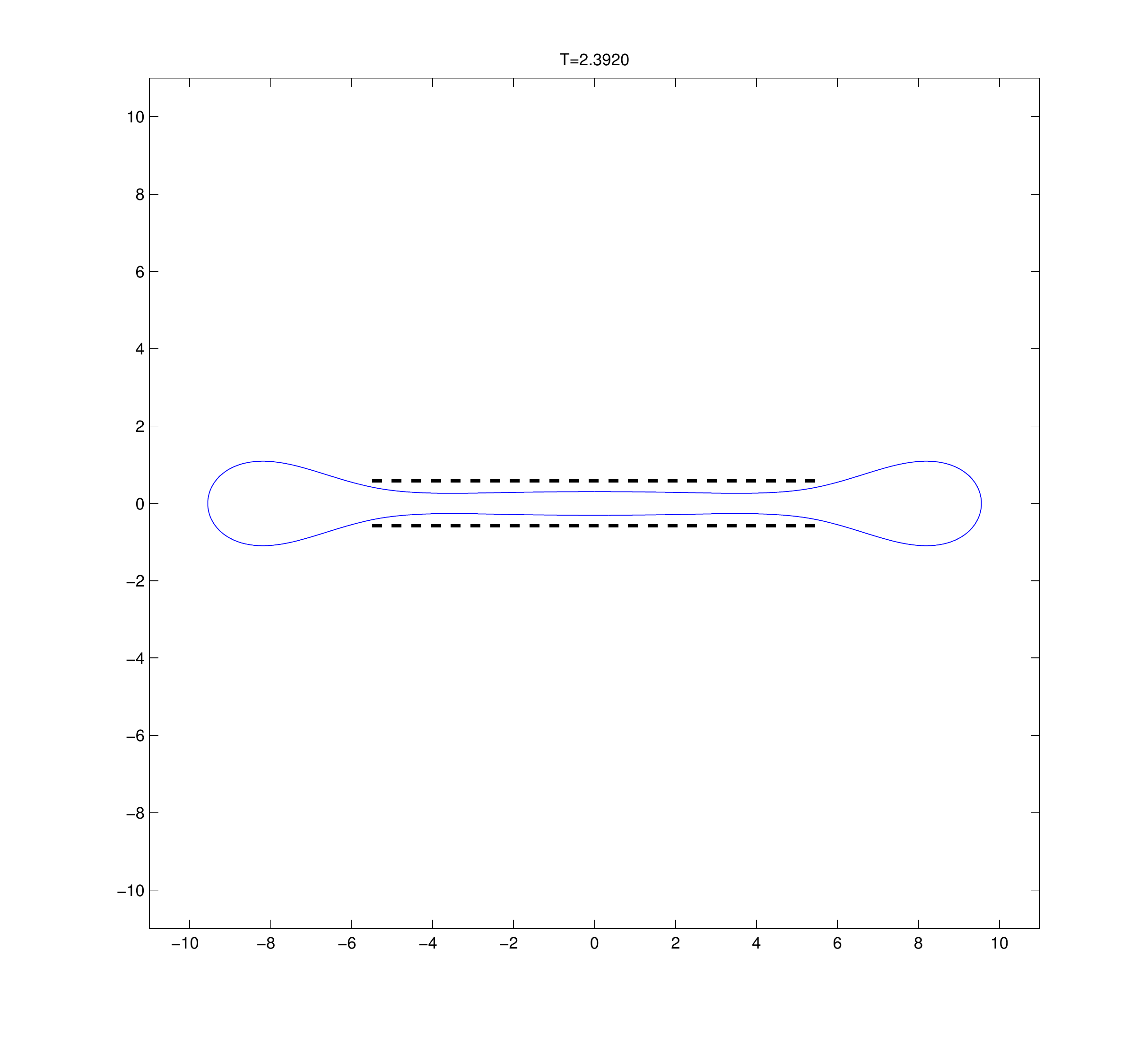}
}
\subfloat[T=3.19]{
\includegraphics[width=65mm]{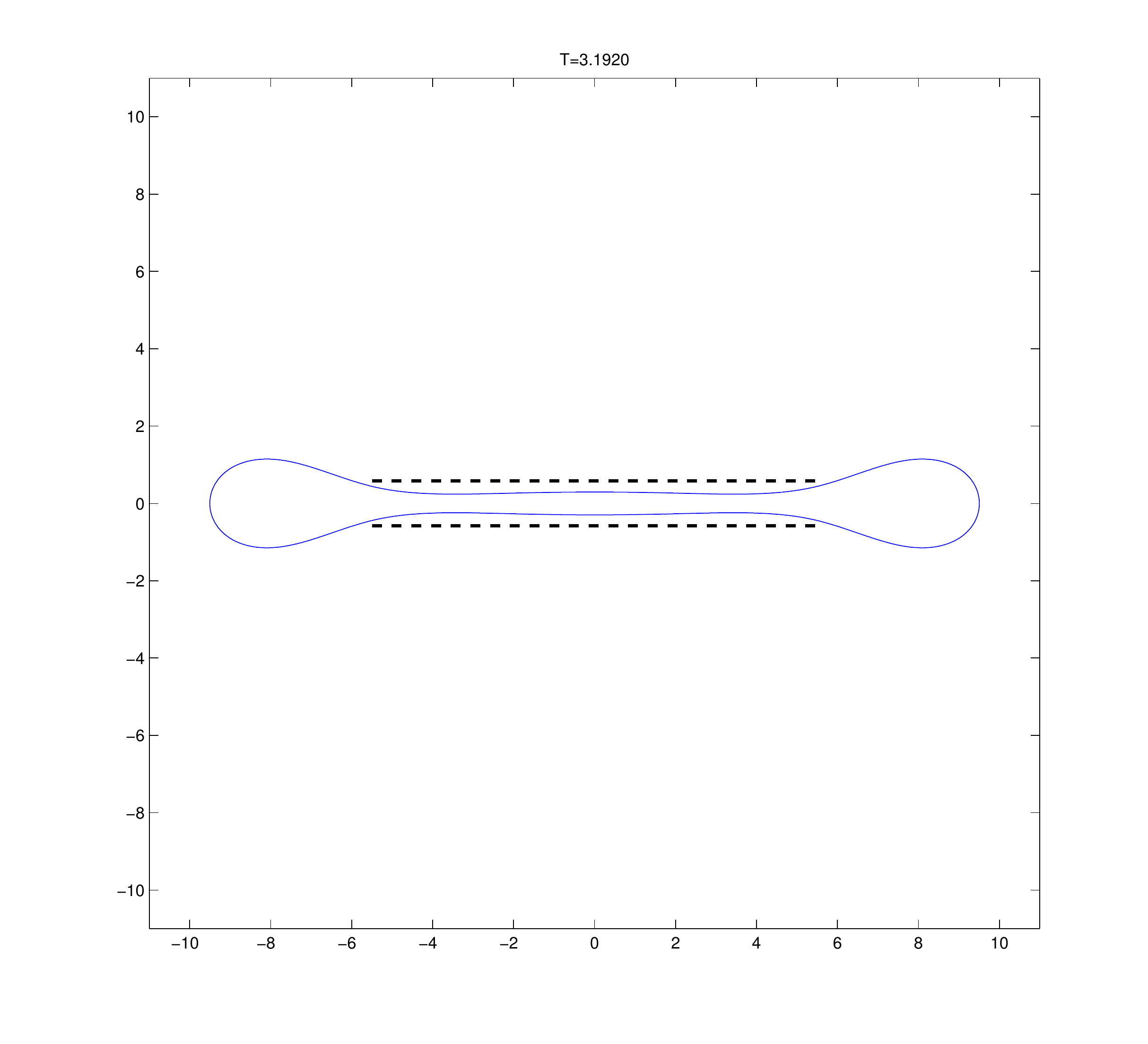}
}
\hspace{0mm}
\subfloat[T=5.19]{
\includegraphics[width=65mm]{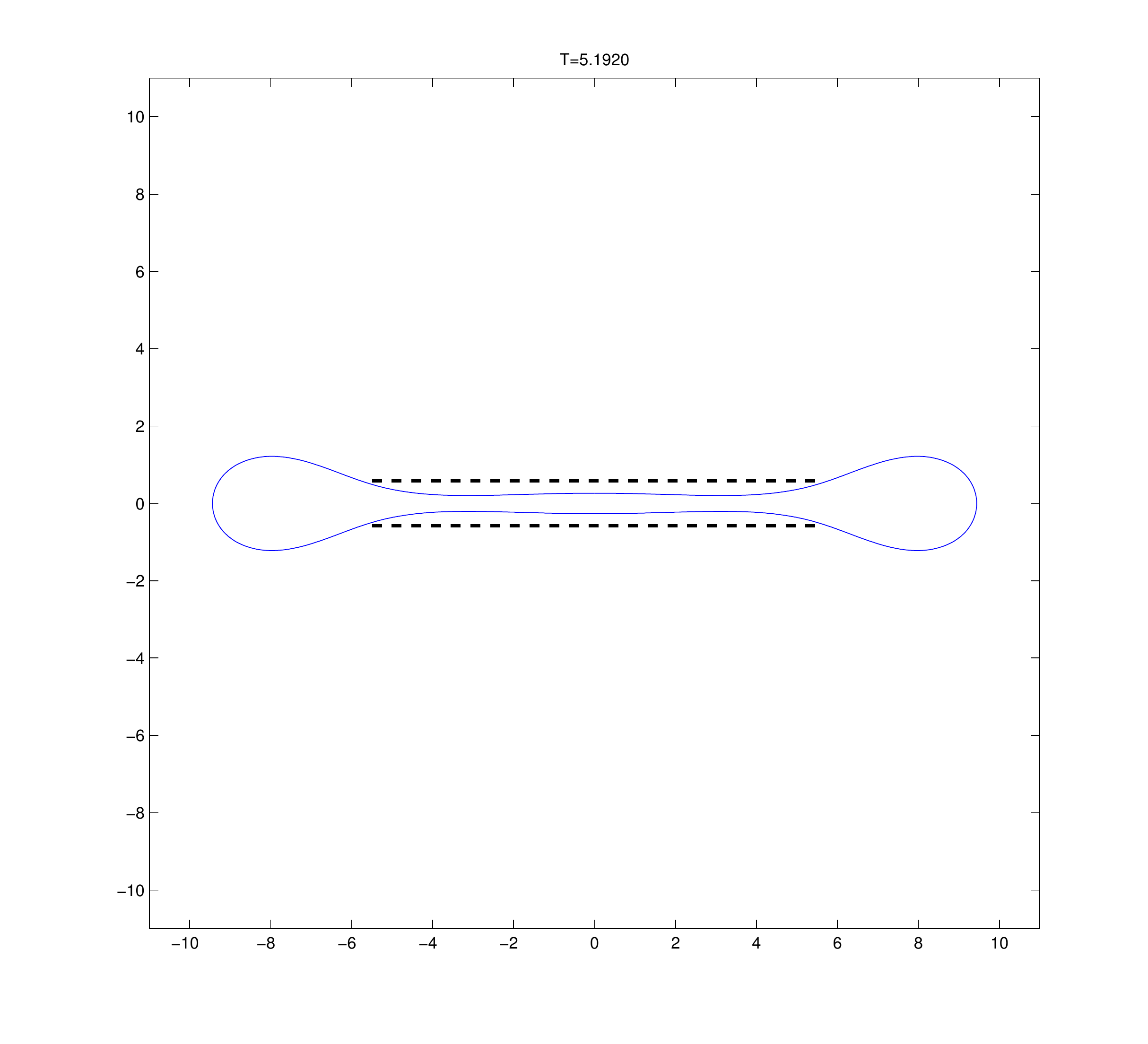}
}
\subfloat[T=6.40]{
\includegraphics[width=65mm]{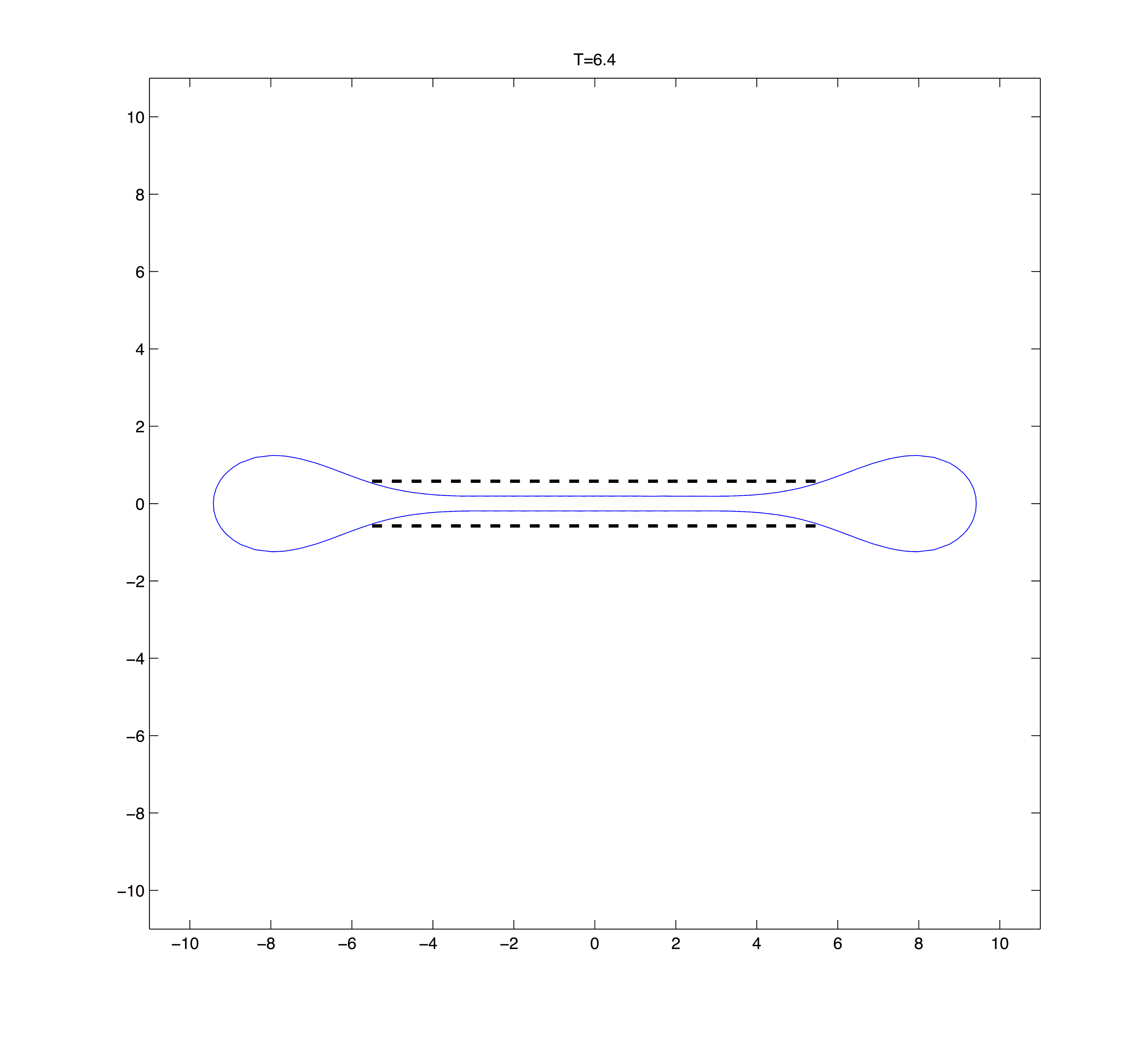}
}

\caption{Modeling of trans-cisternae. The initial shape is set as (a), and the length of the shape $(\Gamma_t)$ is conserved. Besides, the dash lines represent the barriers/obstacles ${\bf B}$. The process from (a) to (h) shows the morphological change of $\Gamma$ when the Willmore energy, equivalent to elastic energy, decreases.}
\label{fig:single1}
\end{figure}

\newpage
\subsubsection*{Observation and Conclusion}
By observing the numerical experiment shown in Figure \ref{fig:single1}, we arrive at the following conclusions.

First, When the obstacles ${\bf B}$ were placed to restrain the vertical expansion of the cisternal membrane $\Gamma$, the growth was limited to the peripheral regions of the cisternae. As a comparison, if we remove the obstacles, the cisternae grows to a peanut shape (as shown in Figure \ref{fig:circle}(b)-(c)) and then it grows isotropically. The cisternae without barriers finally become a sphere (Figure \ref{fig:circle}(h)), which is the optimized shape when the length of $\Gamma$ in $\mathbb{R}^2$ is fixed. In conclusion, Figure \ref{fig:single1} shows an mathematical example that matches the phenomenon that the shapes of trans-cisternae evolve out to the cisternal margin by confining the extension of the central domain and decreasing the elastic energy. Note that naturally it is more stable if the energy is low.

\begin{figure}[h!]
\centering
\subfloat[Initial Shape]{
\includegraphics[width=60mm]{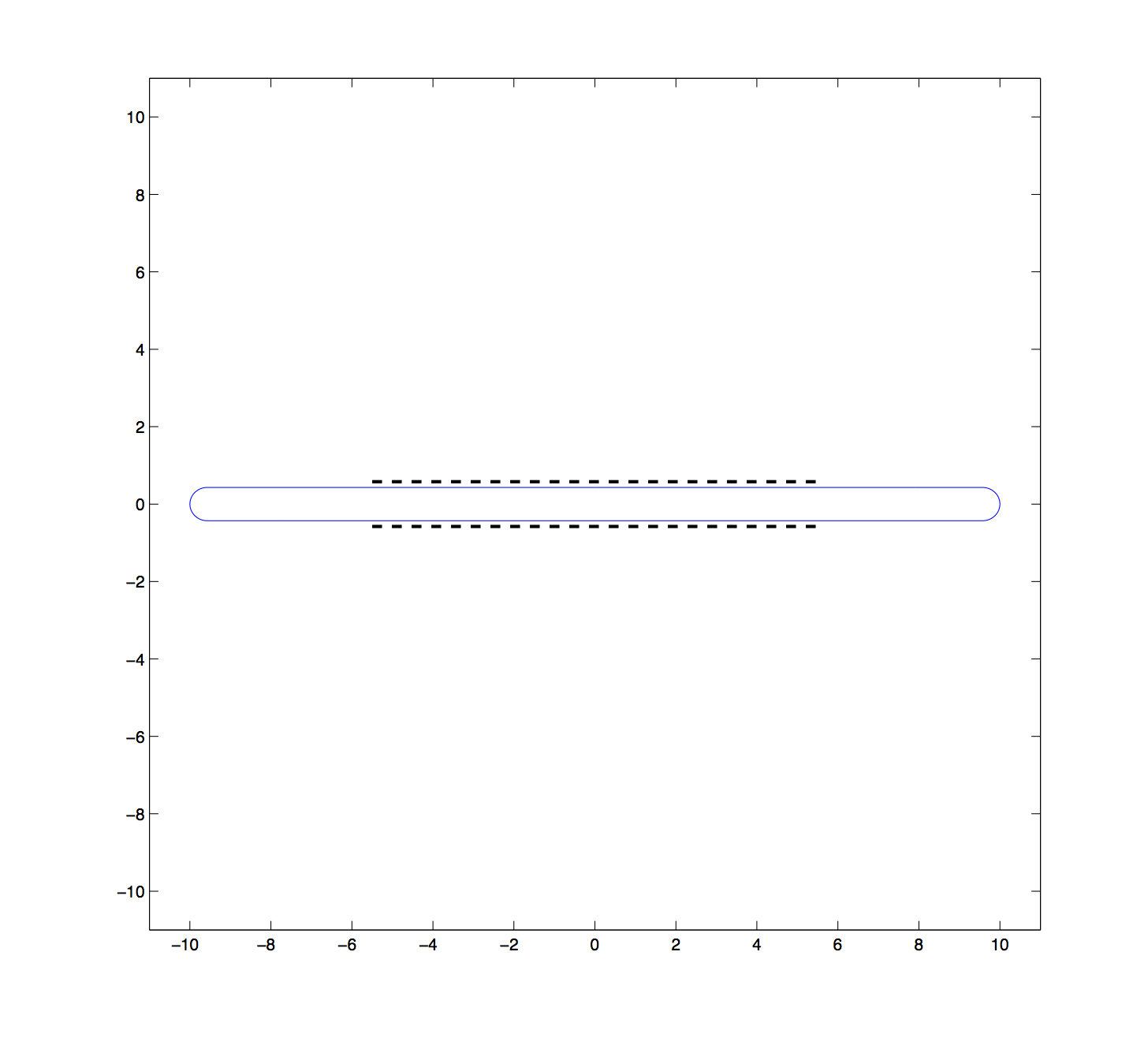}
}
\subfloat[Final Shape]{
\includegraphics[width=60mm]{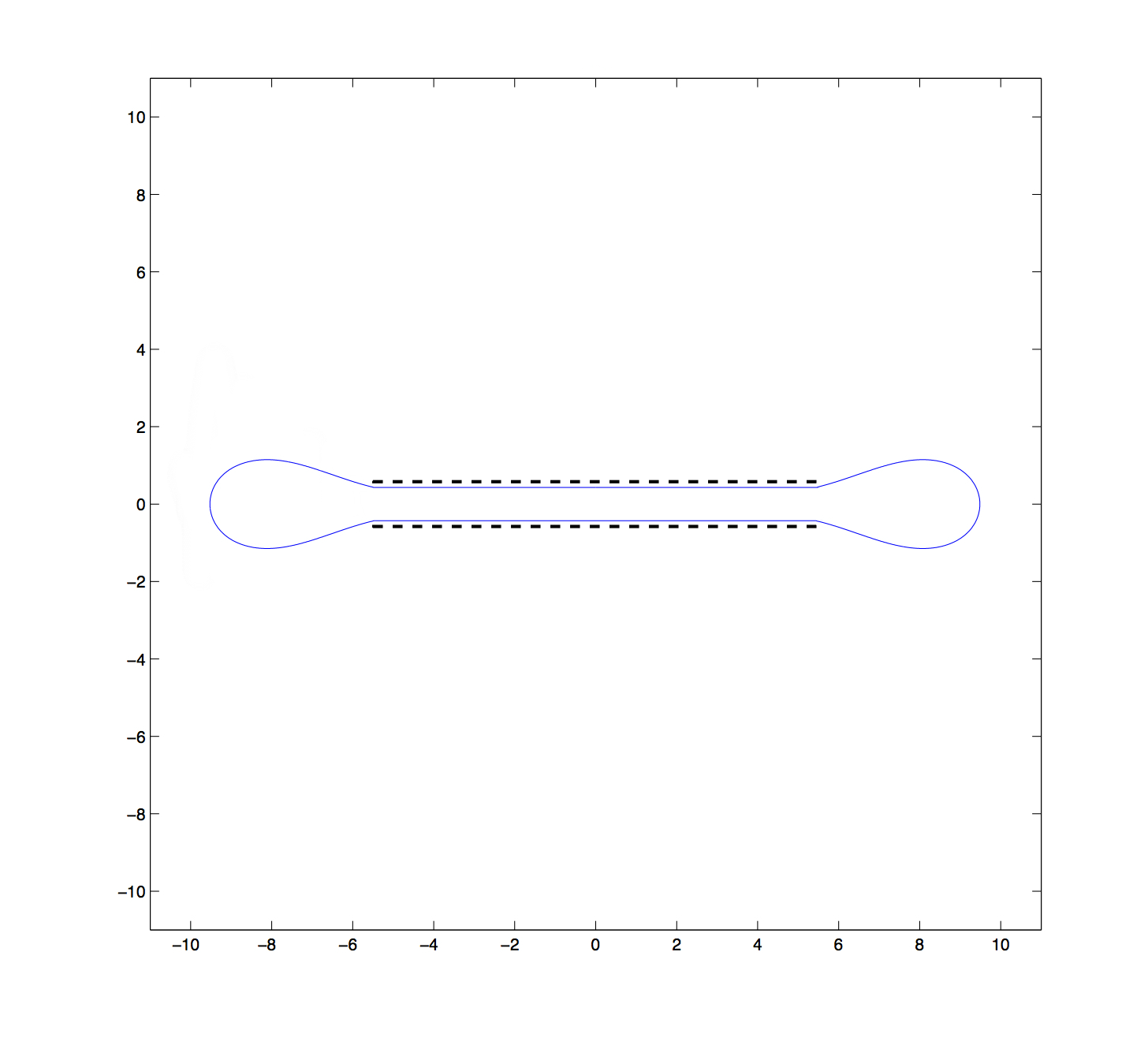}
}
\caption{Counter Example of the Central Thinning}
\label{fig:countra}
\end{figure}

Second, it is intriguing to find that the central domain of the cisternae gets thinner when the elastic energy decreases. It is not an essential phenomenon when the marginal domain is swelling. Figure \ref{fig:countra} demonstrates a counter example. In this example, the central domain is also limited by the barriers. The marginal part also swells. The volume also increases. The length is also fixed. However, the central domain dose not become thinner. Hence, we here relate the thinning of the central part of the trans-cisternae, which is observed by previous biological works \cite{thinning1,thinning2}, to the decrease of the elastic energy of the cisternal membrane.

\newpage
\subsubsection{A Single Cisternae with Moving Barriers $\bf B$}

Usually, in the math models, the position of the obstacles ${\bf B}$ is fixed. However, inspired by the hypothesis mentioned by Prof. Kang, that the intercisternal elements (such as Golgi matrix) maintain the same distance with the membrane since some of its components are embedded in the membrane. Hence, it could be more realistic to keep the distance between the membrane and the barriers. 
\begin{figure}[h!]
\centering
\subfloat[T=0]{
\includegraphics[width=65mm]{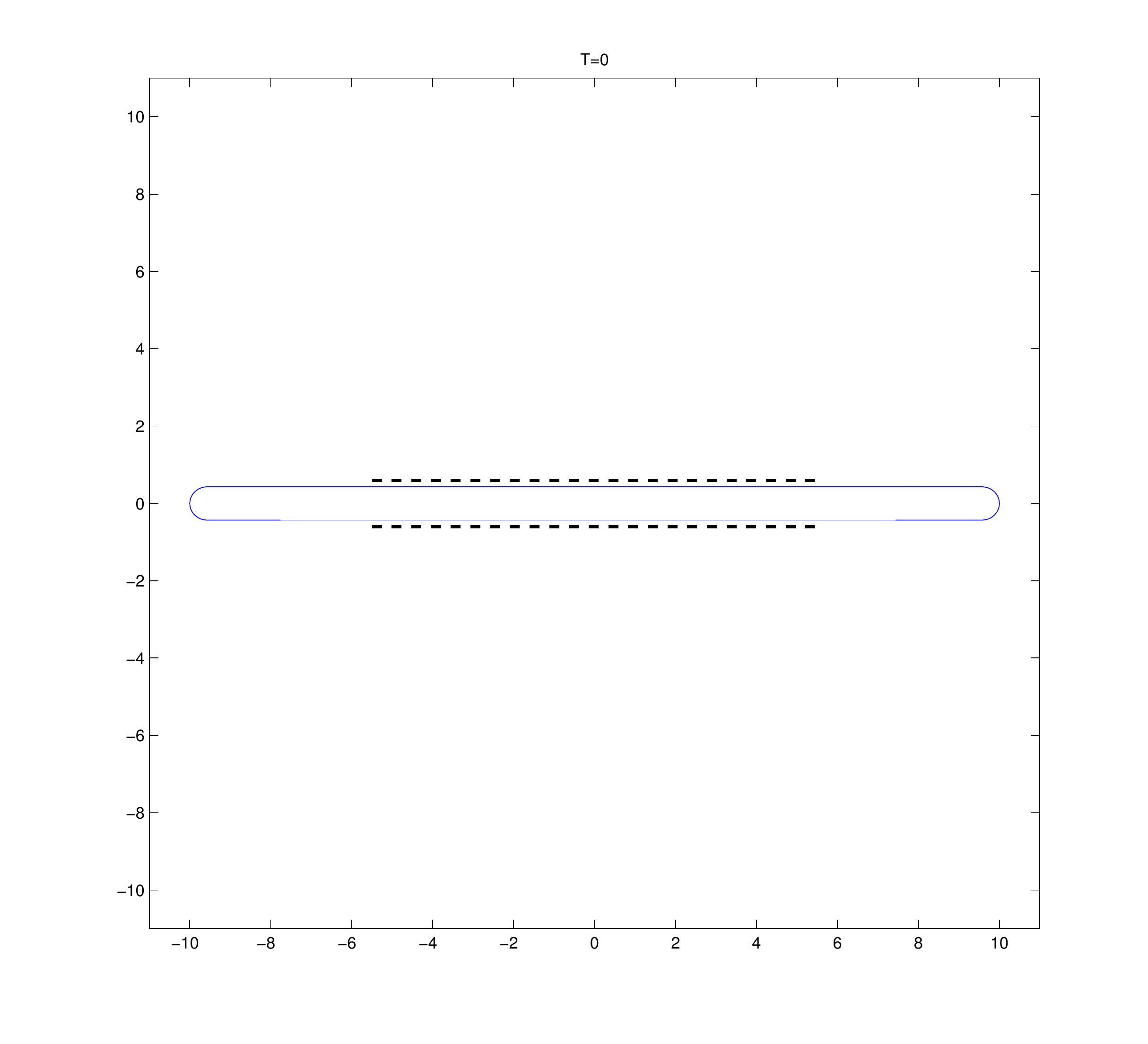}
}
\subfloat[T=0.08]{
\includegraphics[width=65mm]{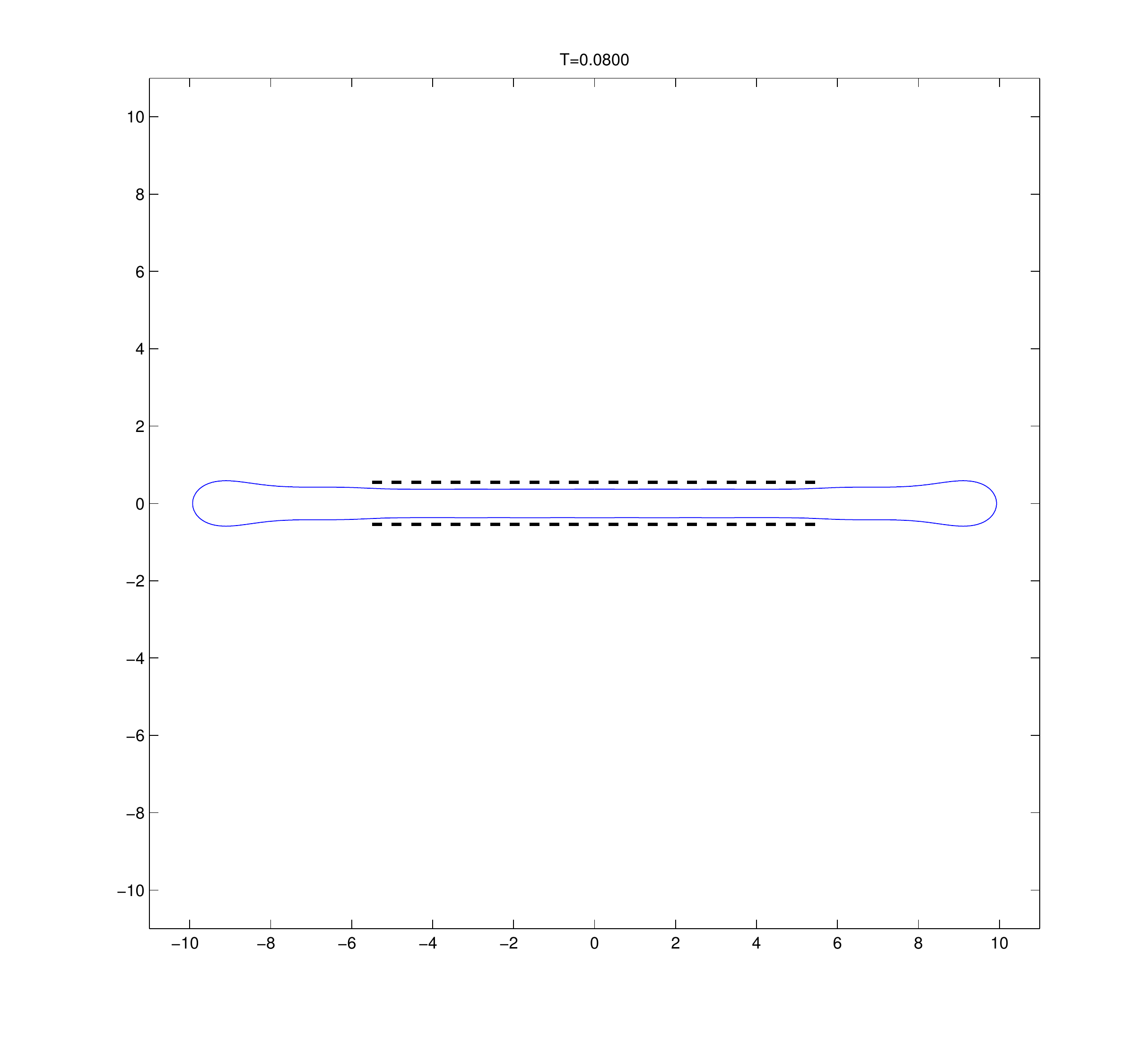}
}
\hspace{0mm}
\subfloat[T=0.16]{
\includegraphics[width=65mm]{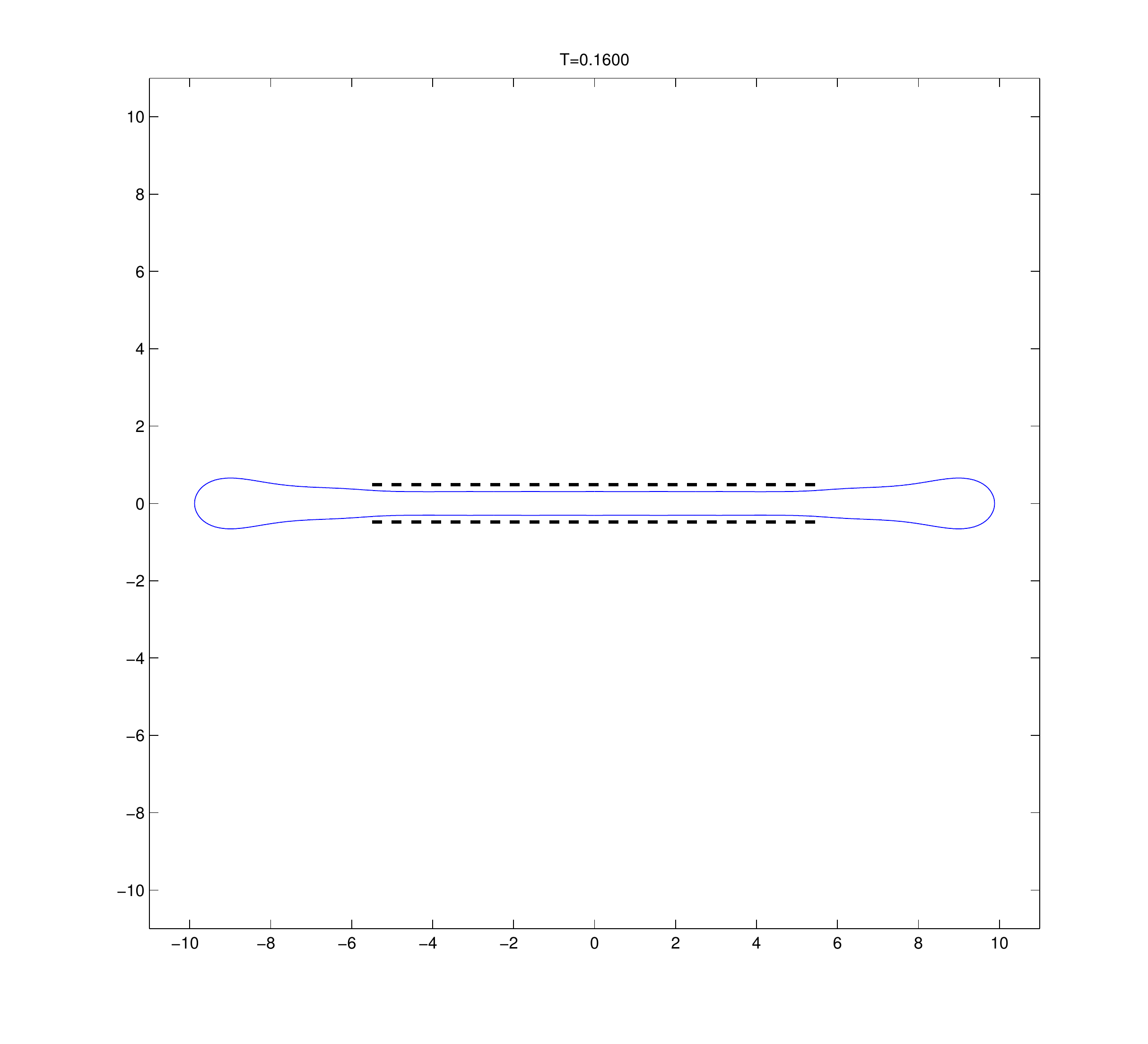}
}
\subfloat[T=0.24]{
\includegraphics[width=65mm]{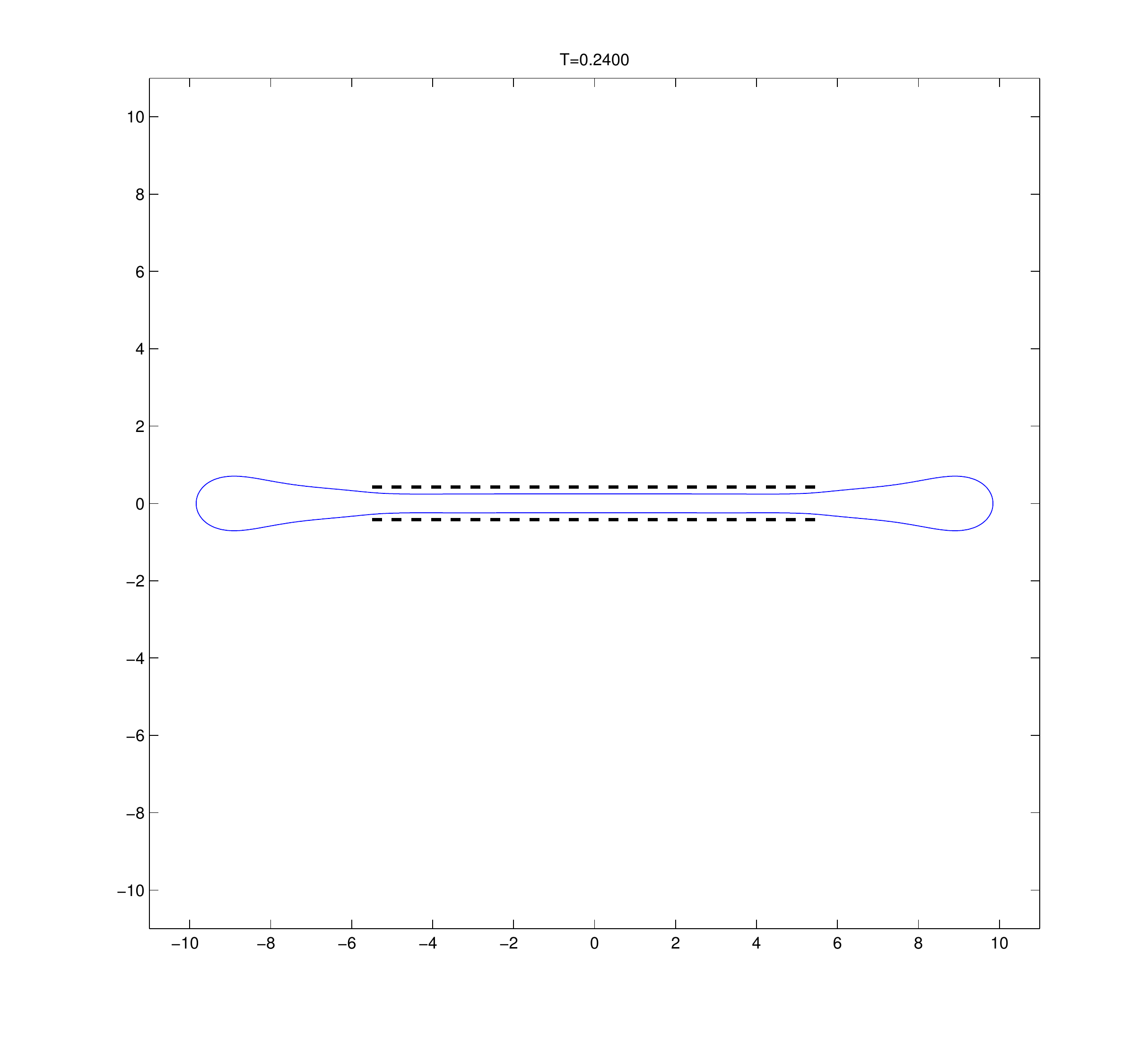}
}
\caption{The barrier ${\bf B}$ is moving dependent on $\Gamma(t)$. When the central part of the cisternae is getting thinner, the dash lines (${\bf B}$) are getting closer. (to be continued)}
\label{fig:MovingB}
\end{figure}

\begin{figure}[h!]
\centering
\subfloat[T=0.32]{
\includegraphics[width=65mm]{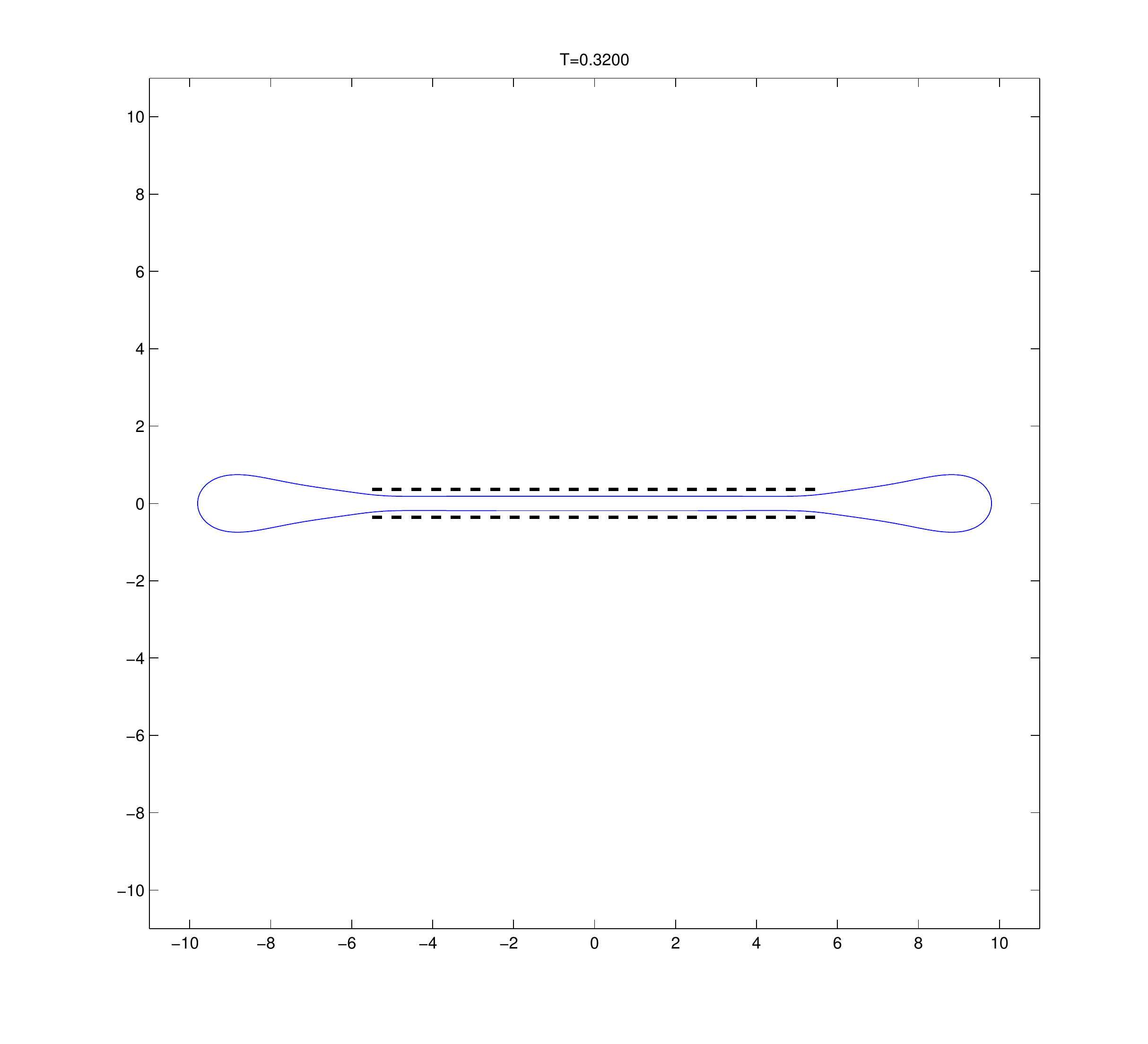}
}
\subfloat[T=0.36]{
\includegraphics[width=65mm]{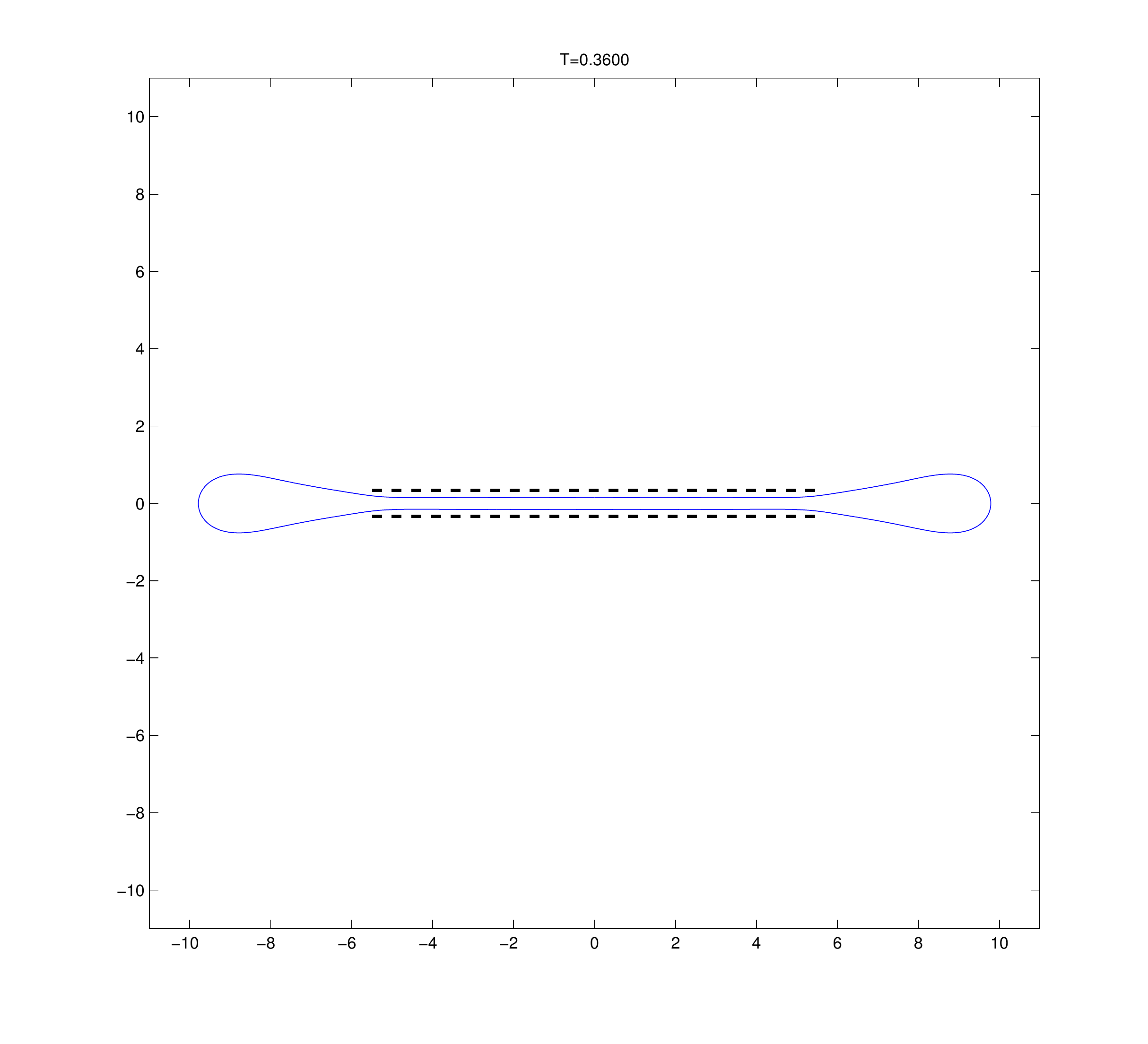}
}
\hspace{0mm}
\subfloat[T=0.40]{
\includegraphics[width=65mm]{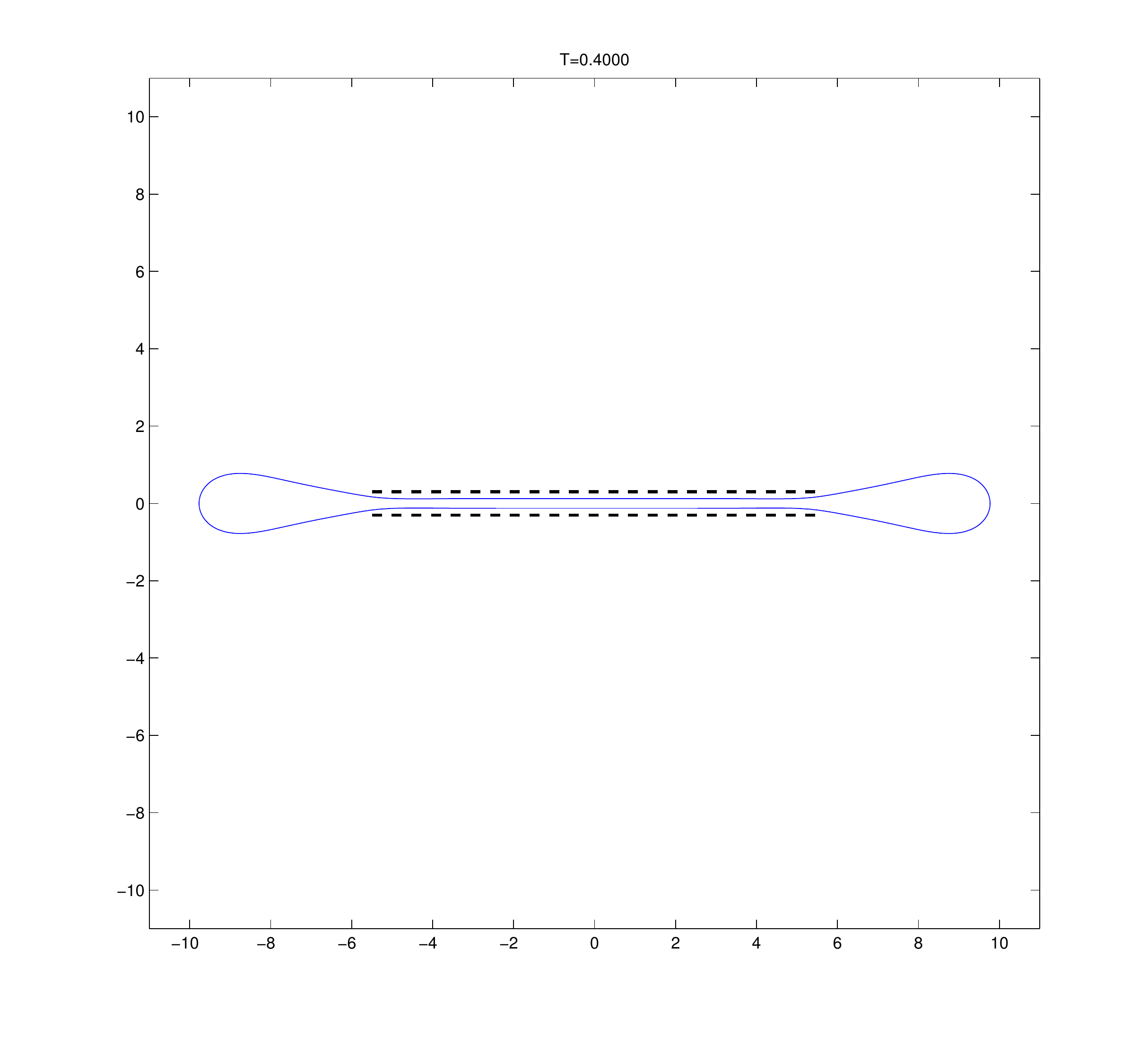}
}
\subfloat[T=0.48]{
\includegraphics[width=65mm]{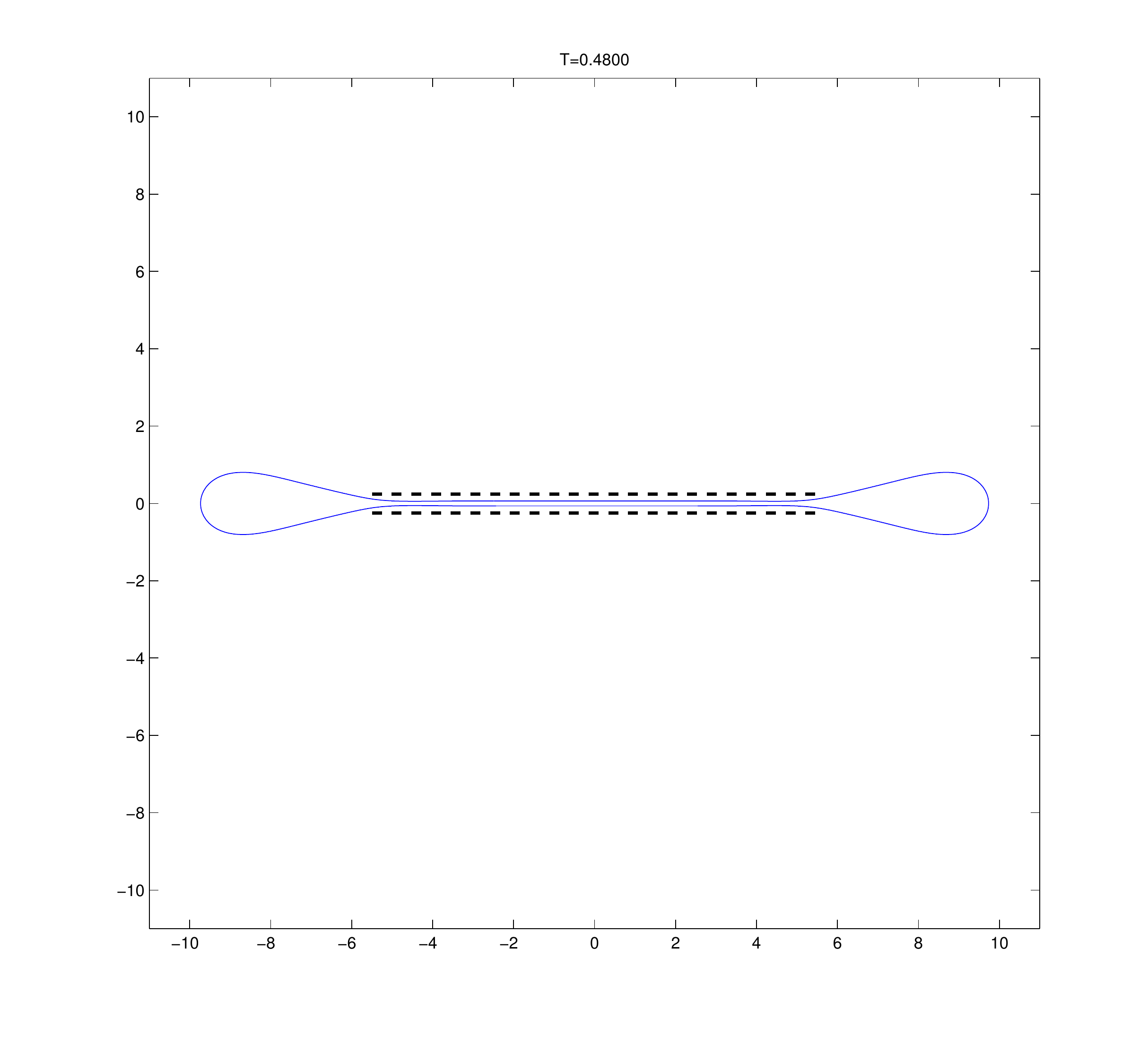}
}
\caption{(continue) The barrier ${\bf B}$ is moving dependent on $\Gamma(t)$. When the central part of the cisternae is getting thinner, the dash lines (${\bf B}$) are getting closer.}
\end{figure}
\newpage
\subsubsection*{A Single Cisternae without Barriers}
\begin{figure}[h!]
\centering
\subfloat[Initial shape]{
\includegraphics[width=58mm]{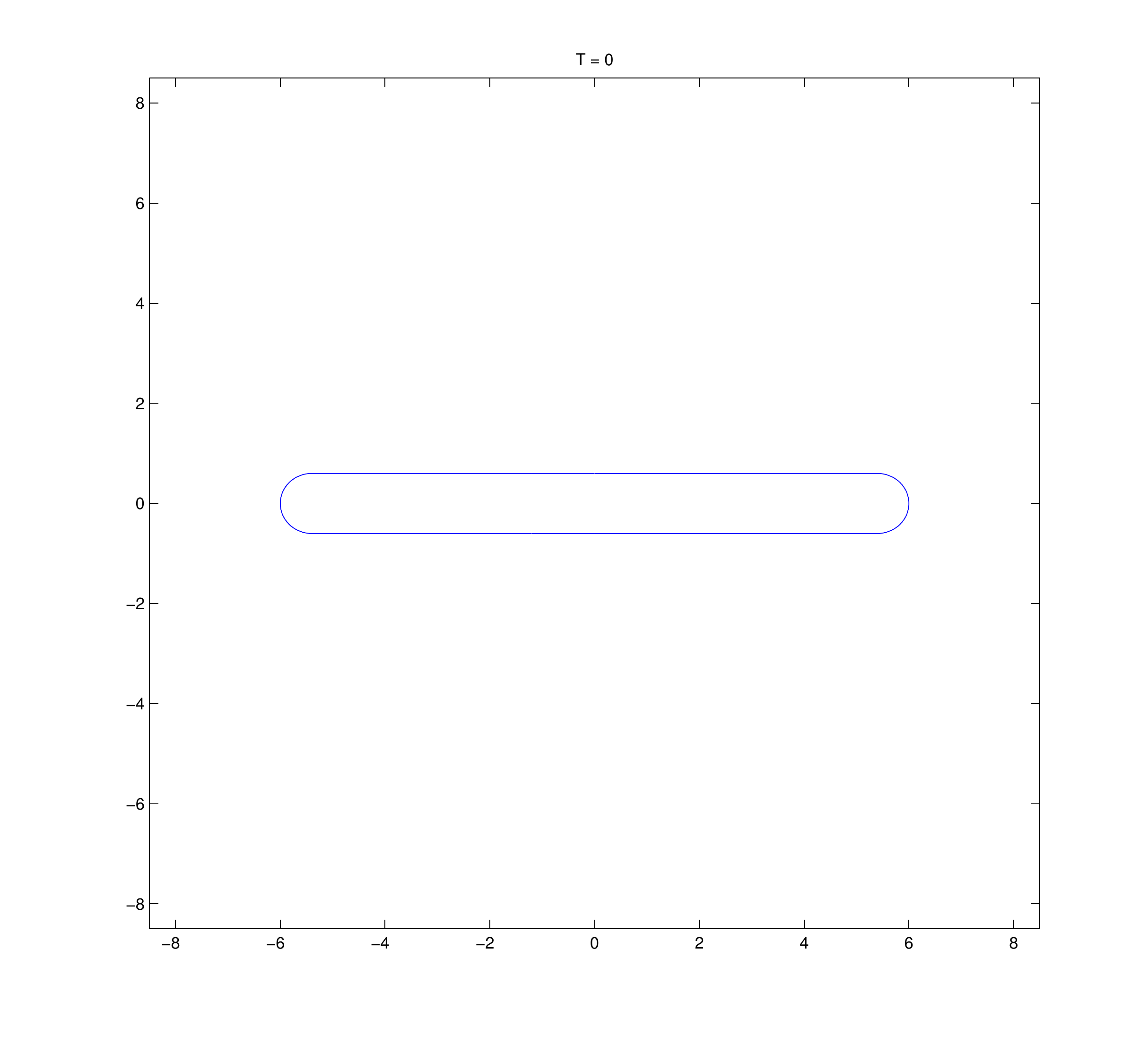}
}
\subfloat[$ $]{
\includegraphics[width=58mm]{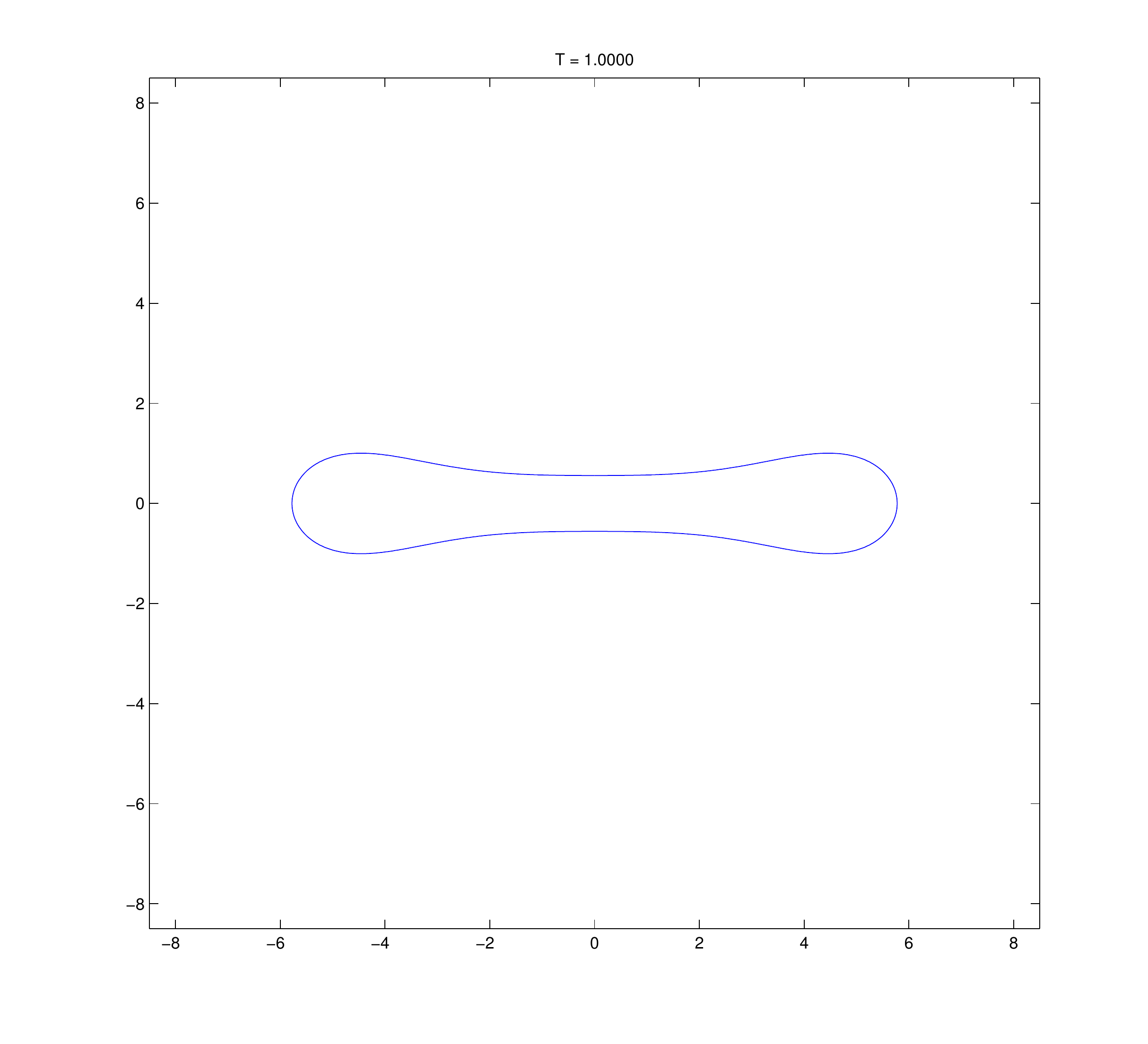}
}
\hspace{0mm}
\subfloat[$ $]{
\includegraphics[width=58mm]{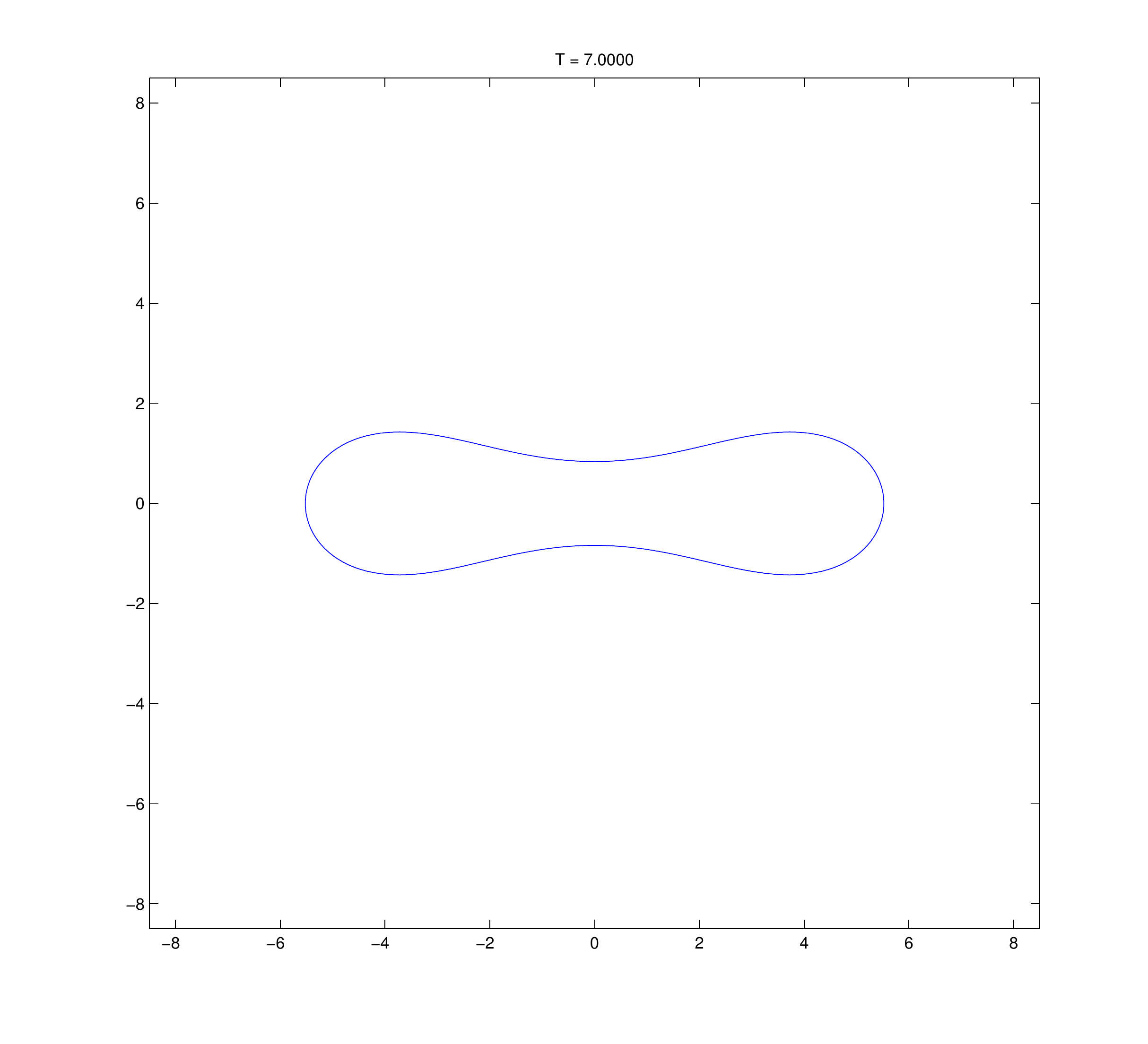}
}
\subfloat[$ $]{
\includegraphics[width=58mm]{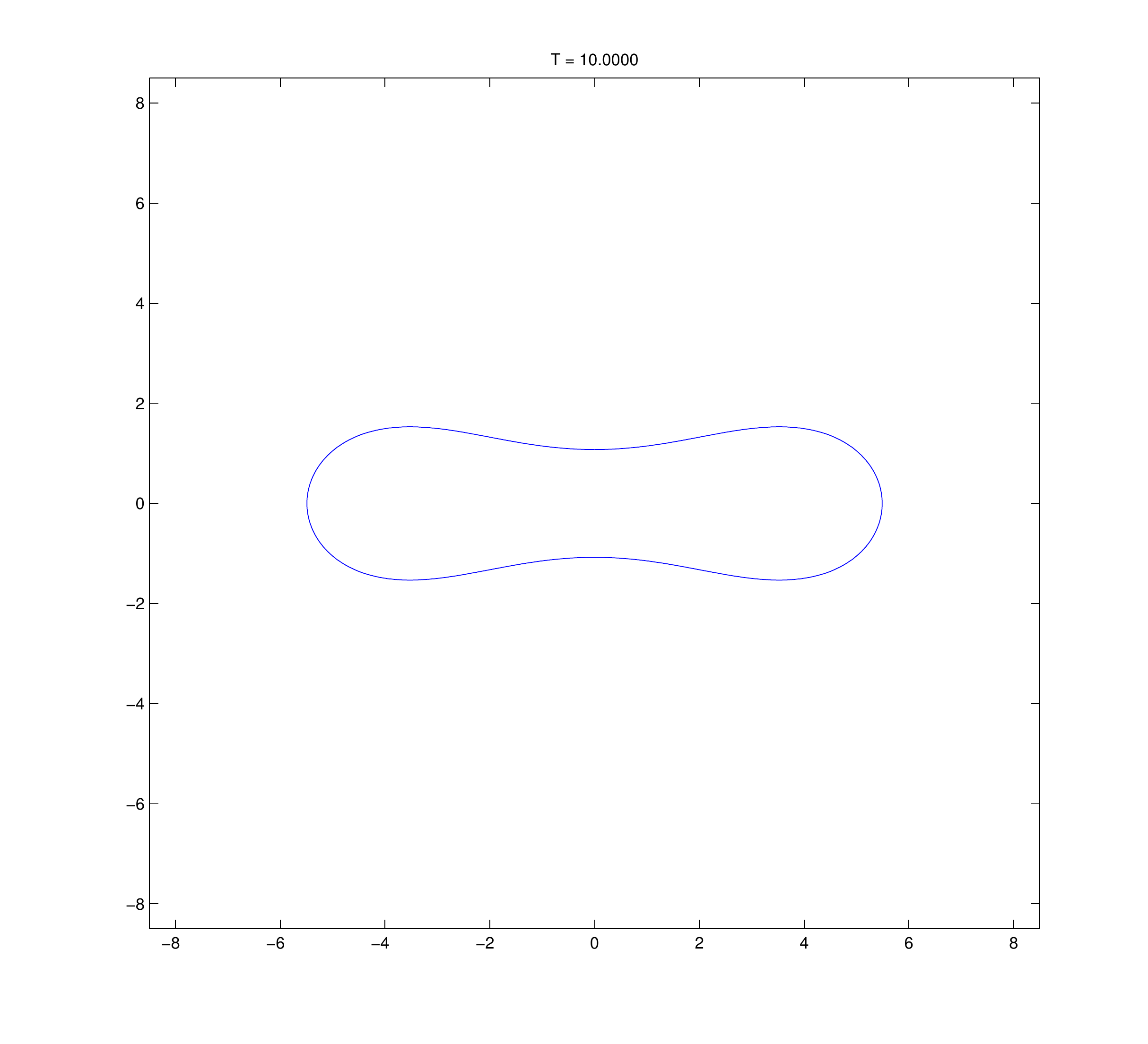}
}
\hspace{0mm}
\subfloat[$ $]{
\includegraphics[width=58mm]{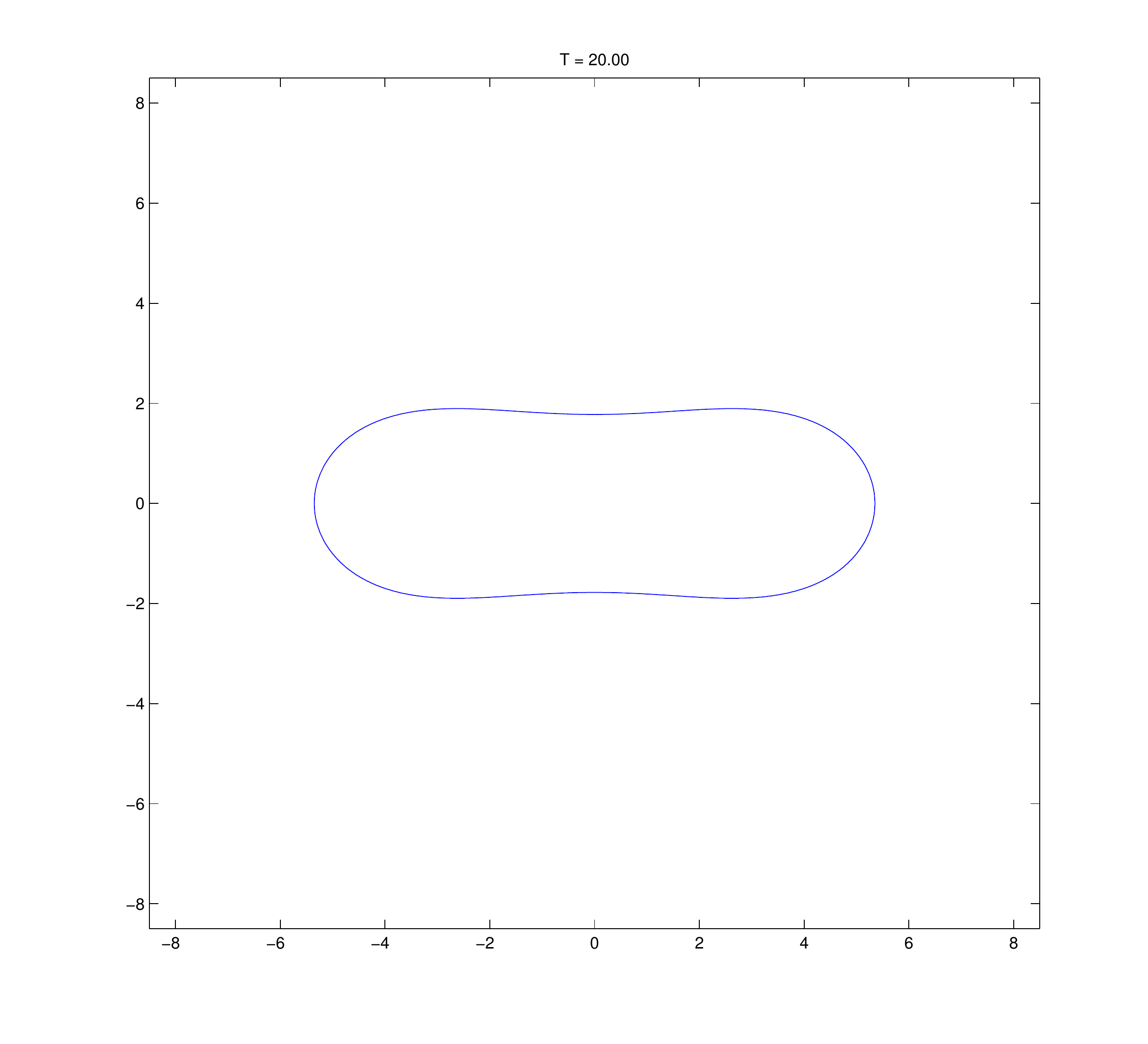}
}
\subfloat[$ $]{
\includegraphics[width=58mm]{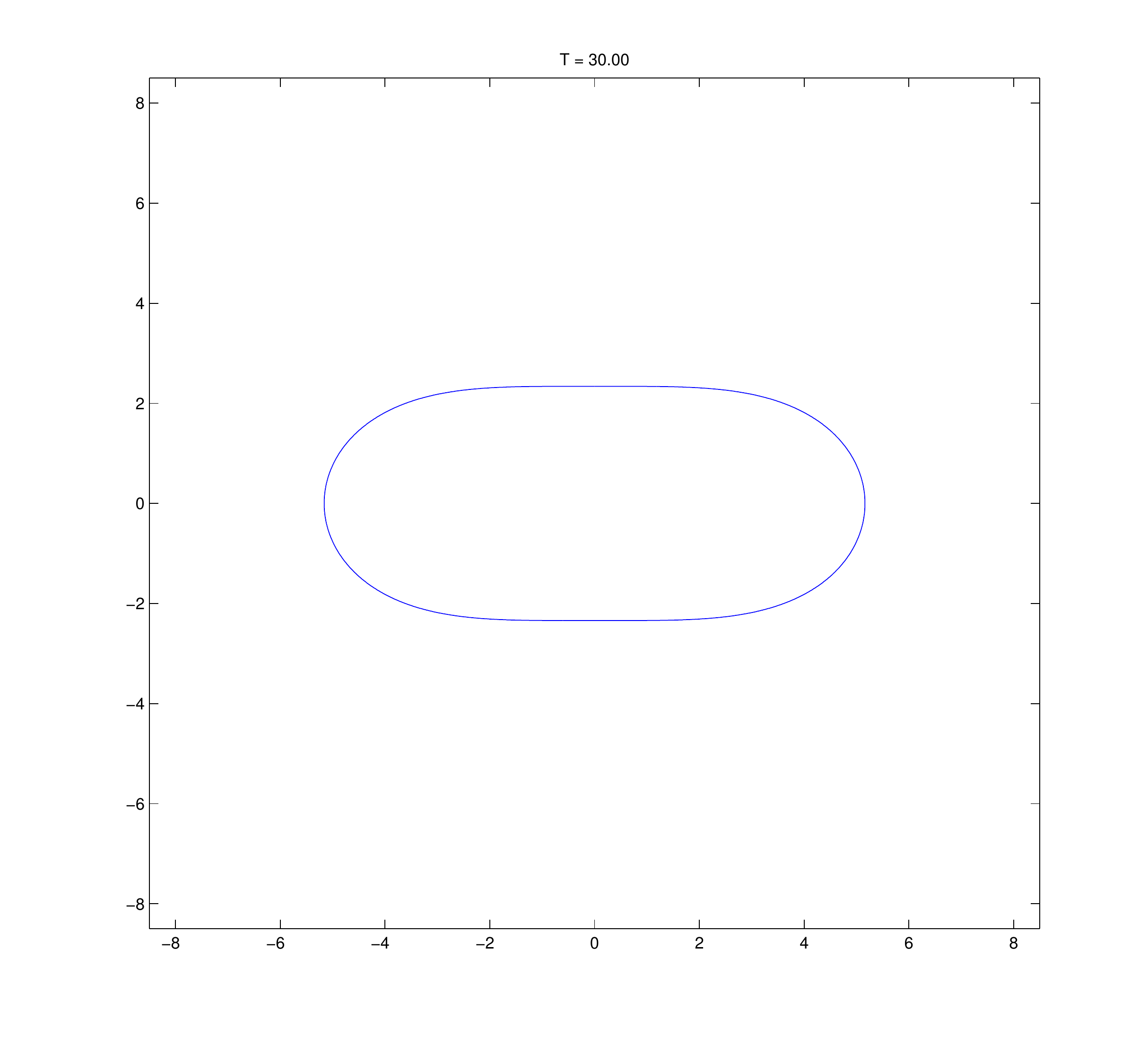}
}
\caption{A Single Cisternae without Barriers}
\label{fig:circle}
\end{figure}

\begin{figure}[h!]\ContinuedFloat
\centering
\subfloat[$ $]{
\includegraphics[width=65mm]{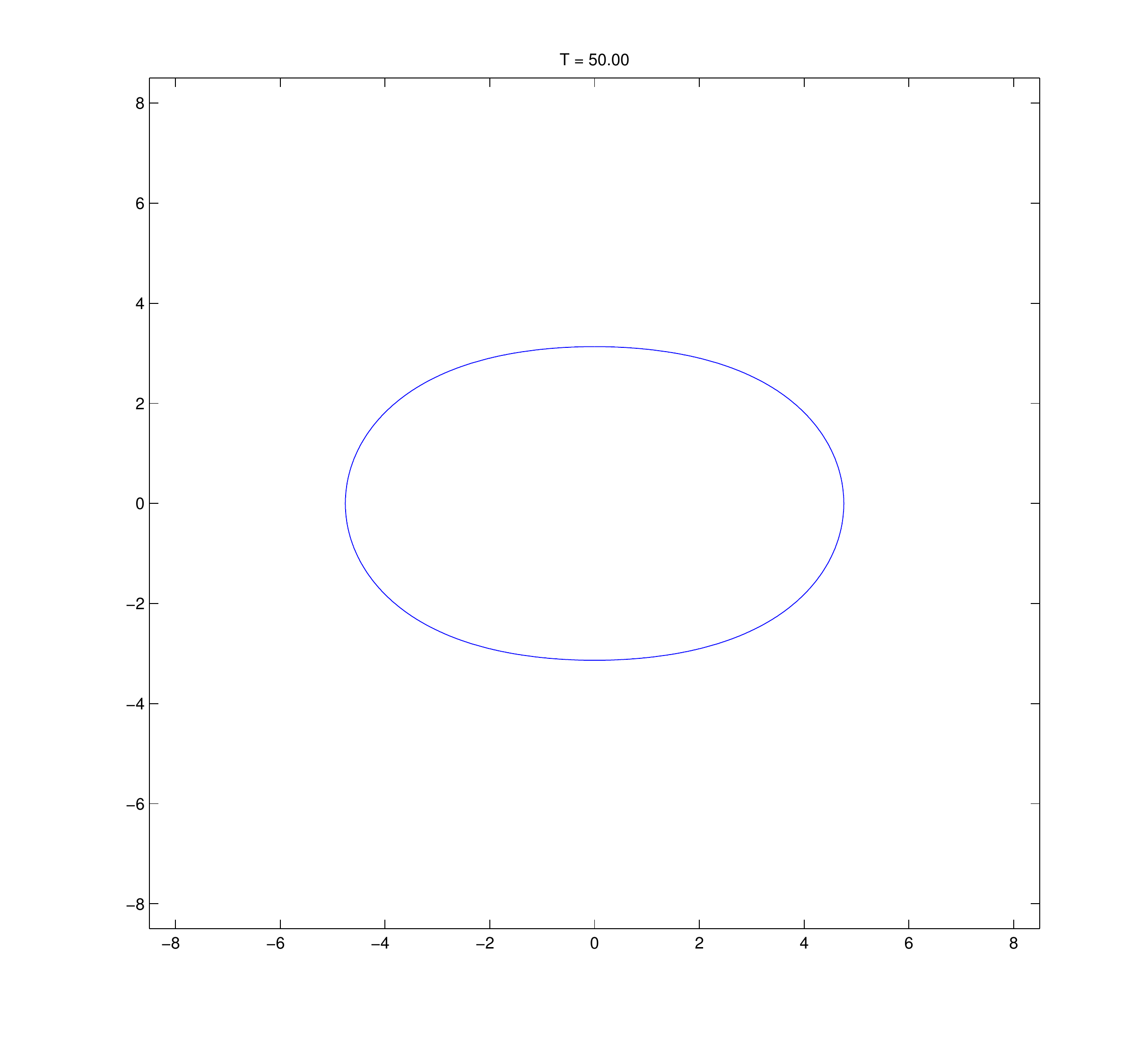}
}
\subfloat[Optimized shape]{
\includegraphics[width=65mm]{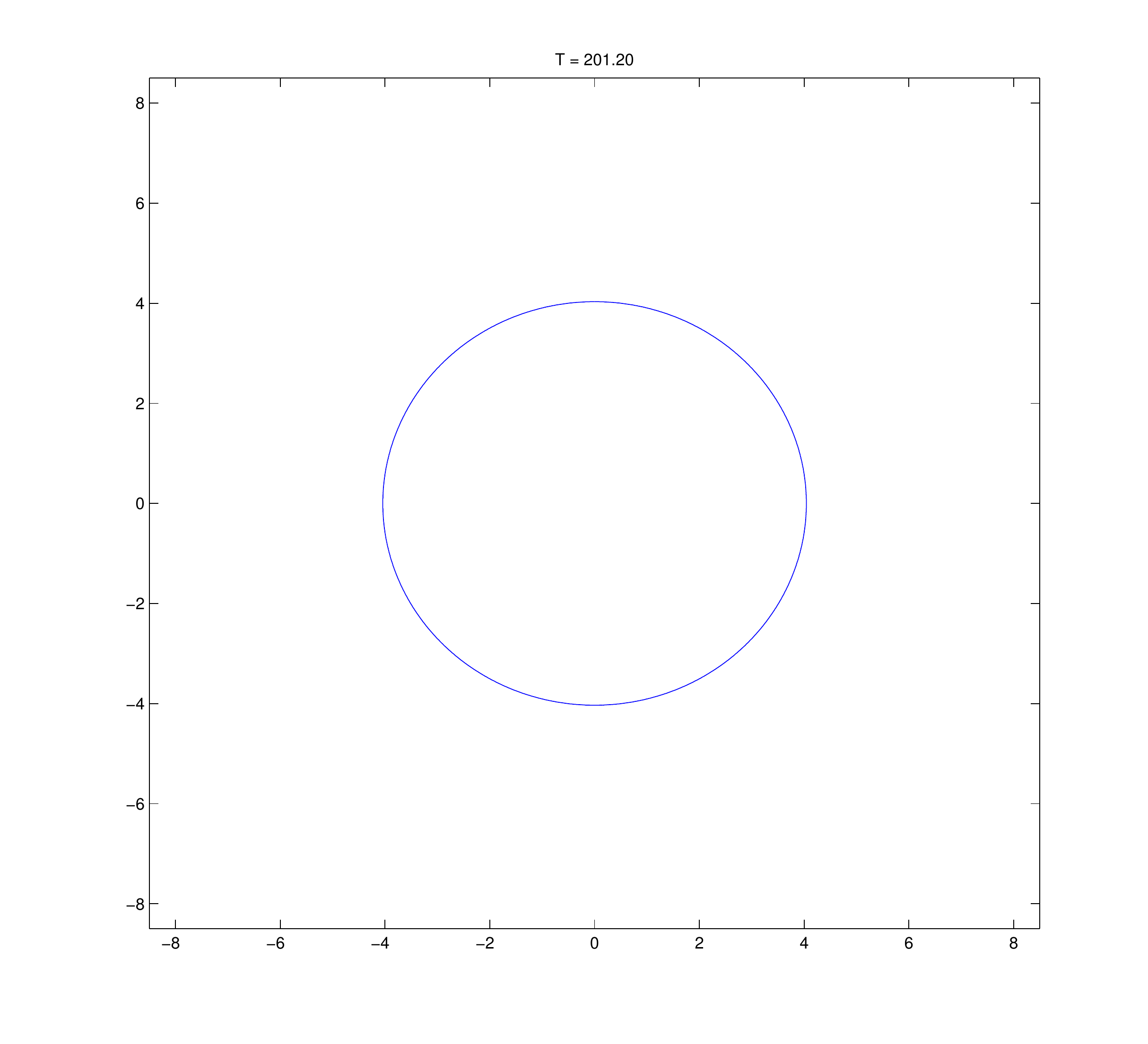}
}
\caption{(continue) A Single Cisternae without Barriers}
\end{figure}

\subsection*{Multiple Layers}
\begin{figure}[h!]
\centering
\includegraphics{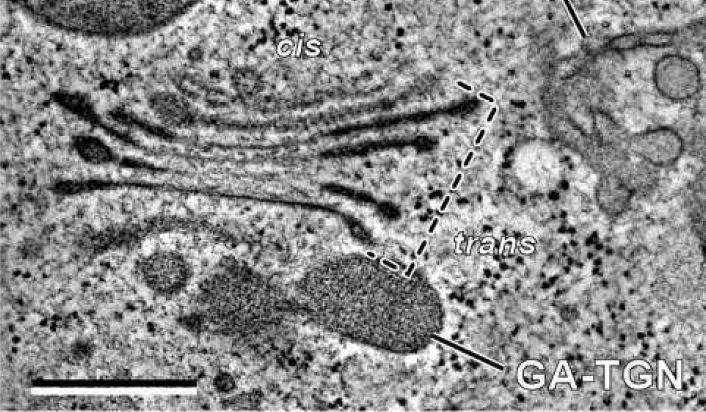}
\caption{Electron Tomographic Slice provided by Prof. Byung-Ho Kang.}
\label{fig:Golgi}
\end{figure}

The Golgi stacks consist of multiple layers of cisternaes (see Figure \ref{fig:Golgi}). This fact inspires us to form a model of multiple vesicles $\{\Gamma_i\}_{i=1}^M$.

\begin{figure}[h!]
\centering
\subfloat[]{
\includegraphics[width=0.33\textwidth]{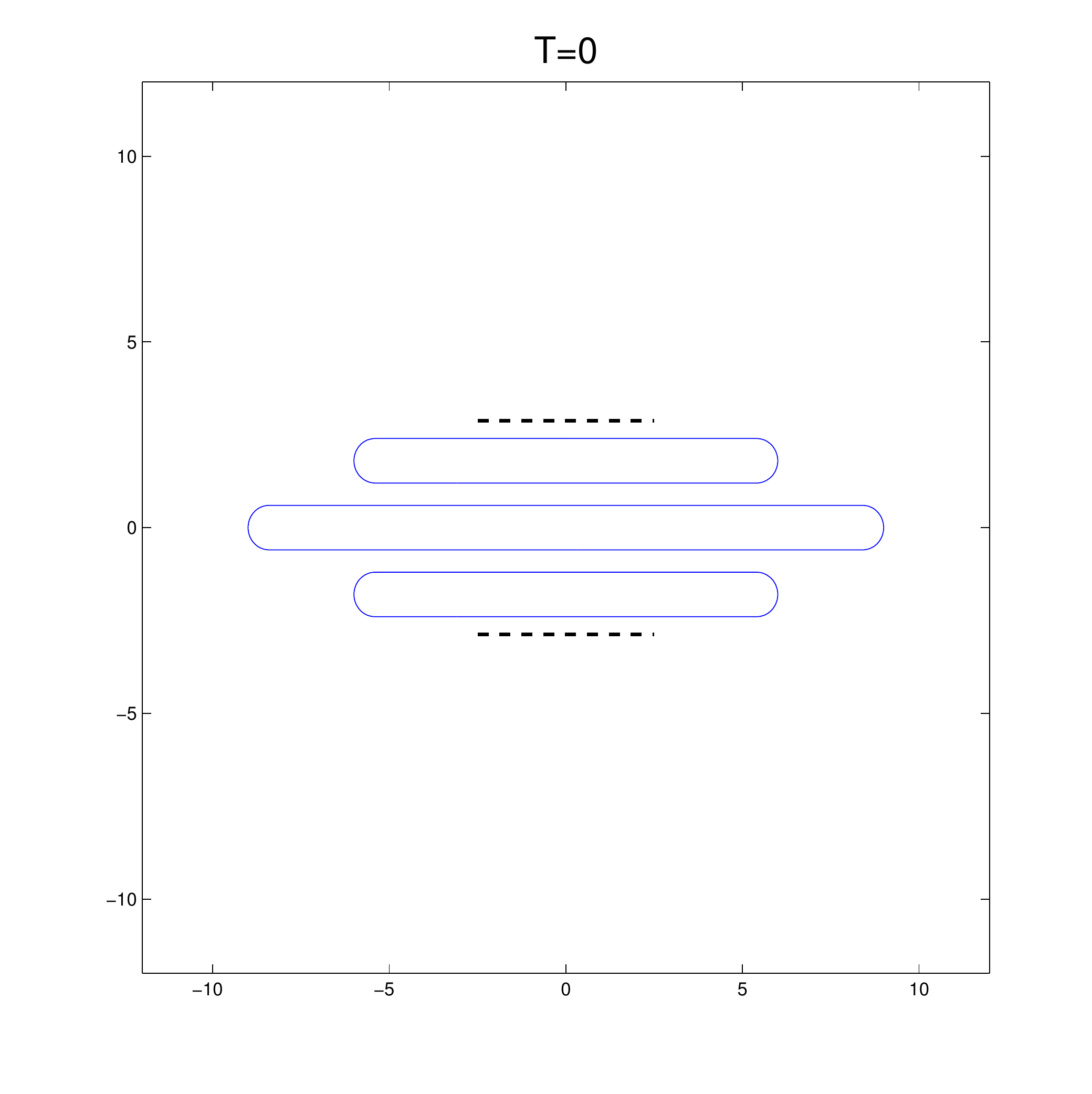}
}
\subfloat[]{
\includegraphics[width=0.33\textwidth]{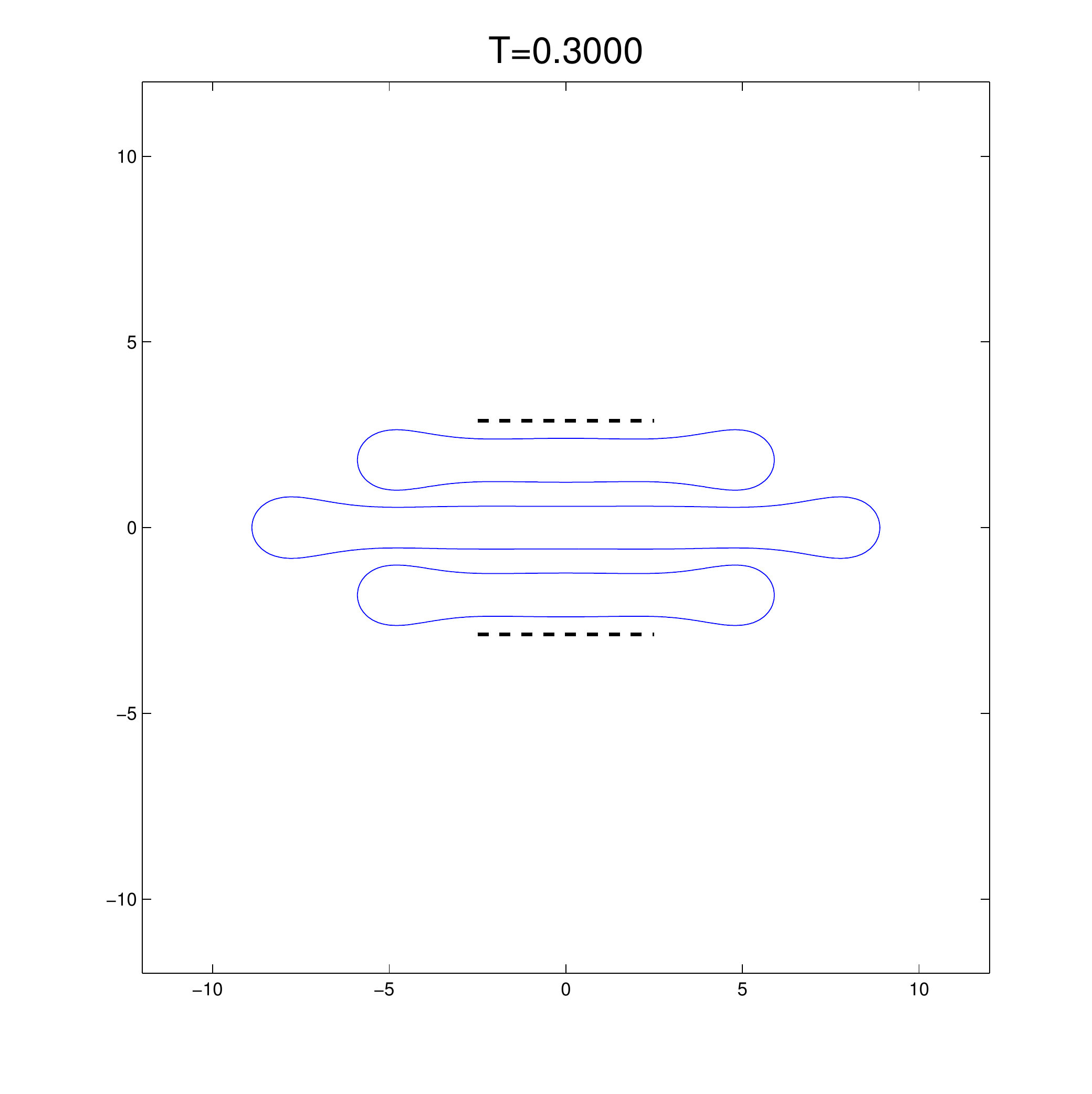}
}
\subfloat[]{
\includegraphics[width=0.33\textwidth]{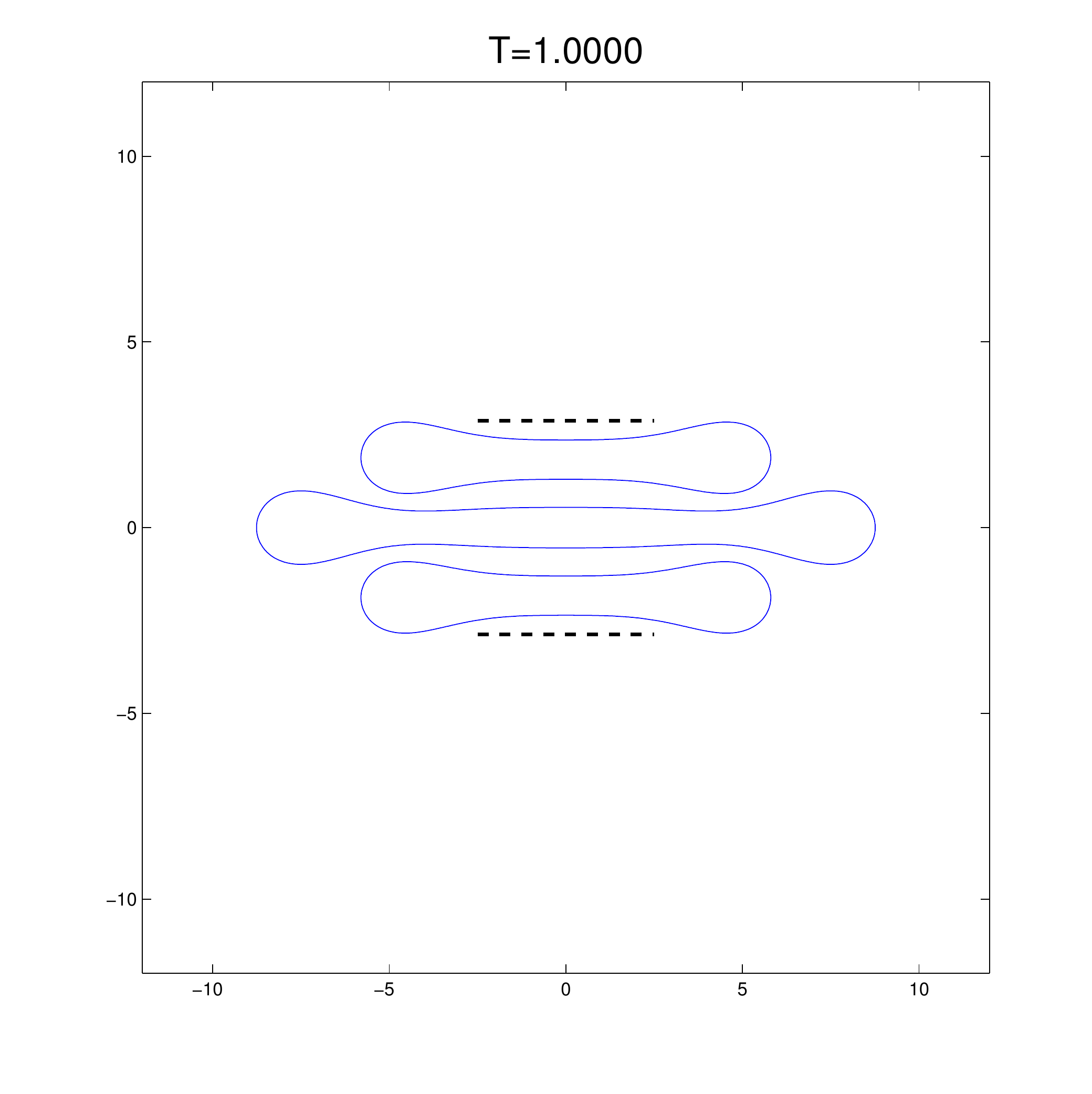}
}
\hspace{0mm}
\subfloat[]{
\includegraphics[width=0.33\textwidth]{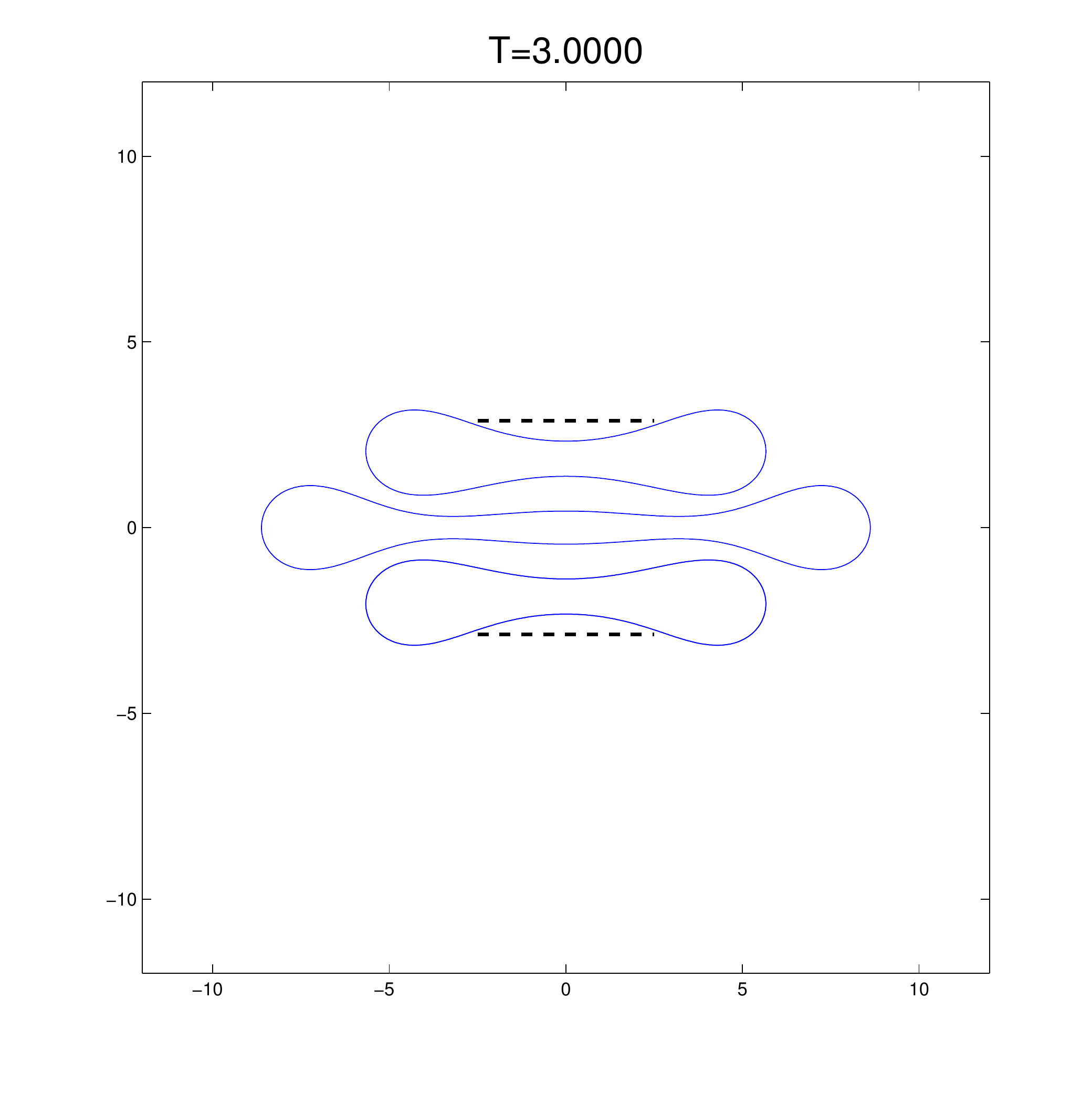}
}
\subfloat[]{
\includegraphics[width=0.33\textwidth]{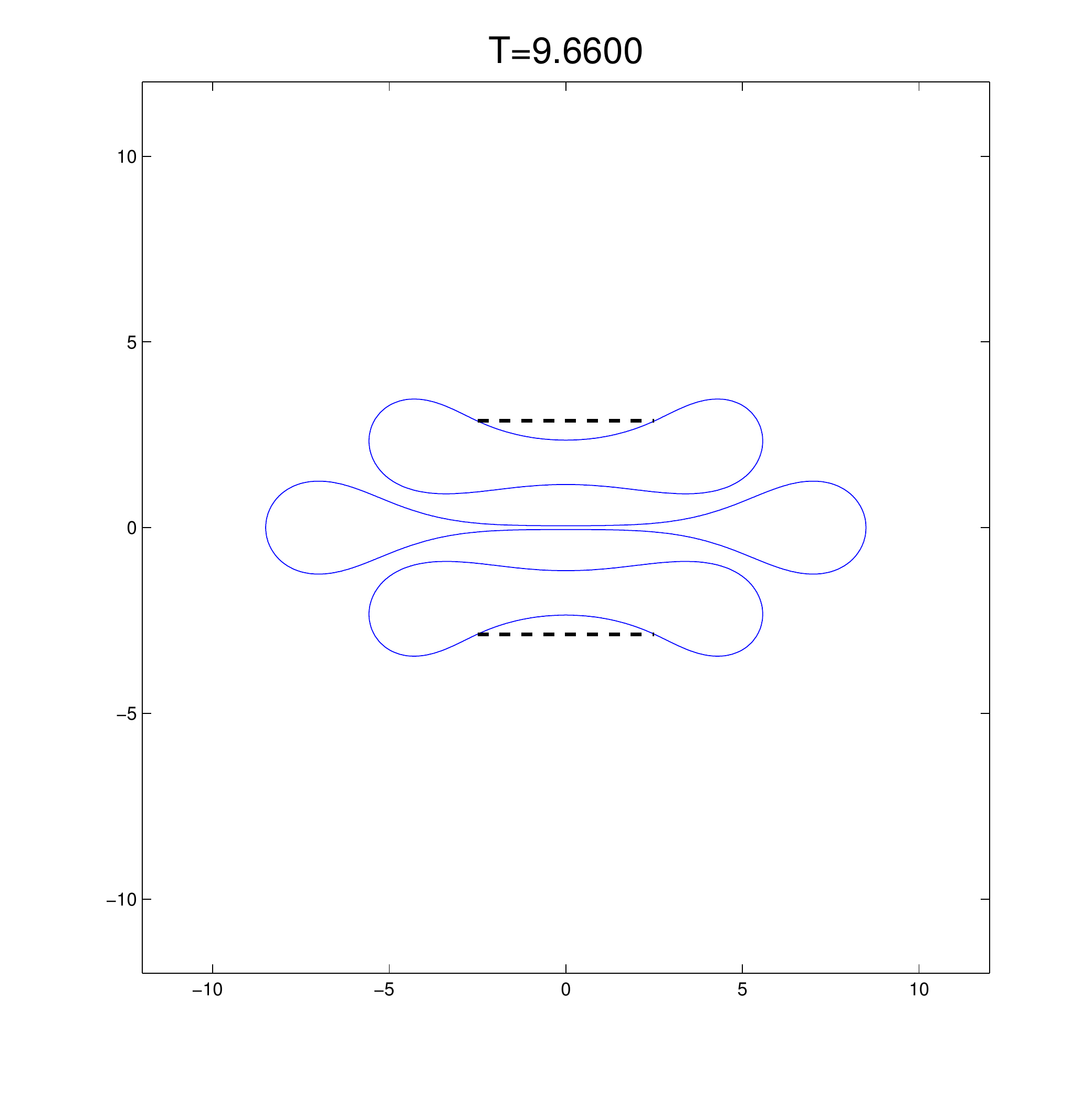}
}
\subfloat[From T=0 to T=9.66]{
\includegraphics[width=0.33\textwidth]{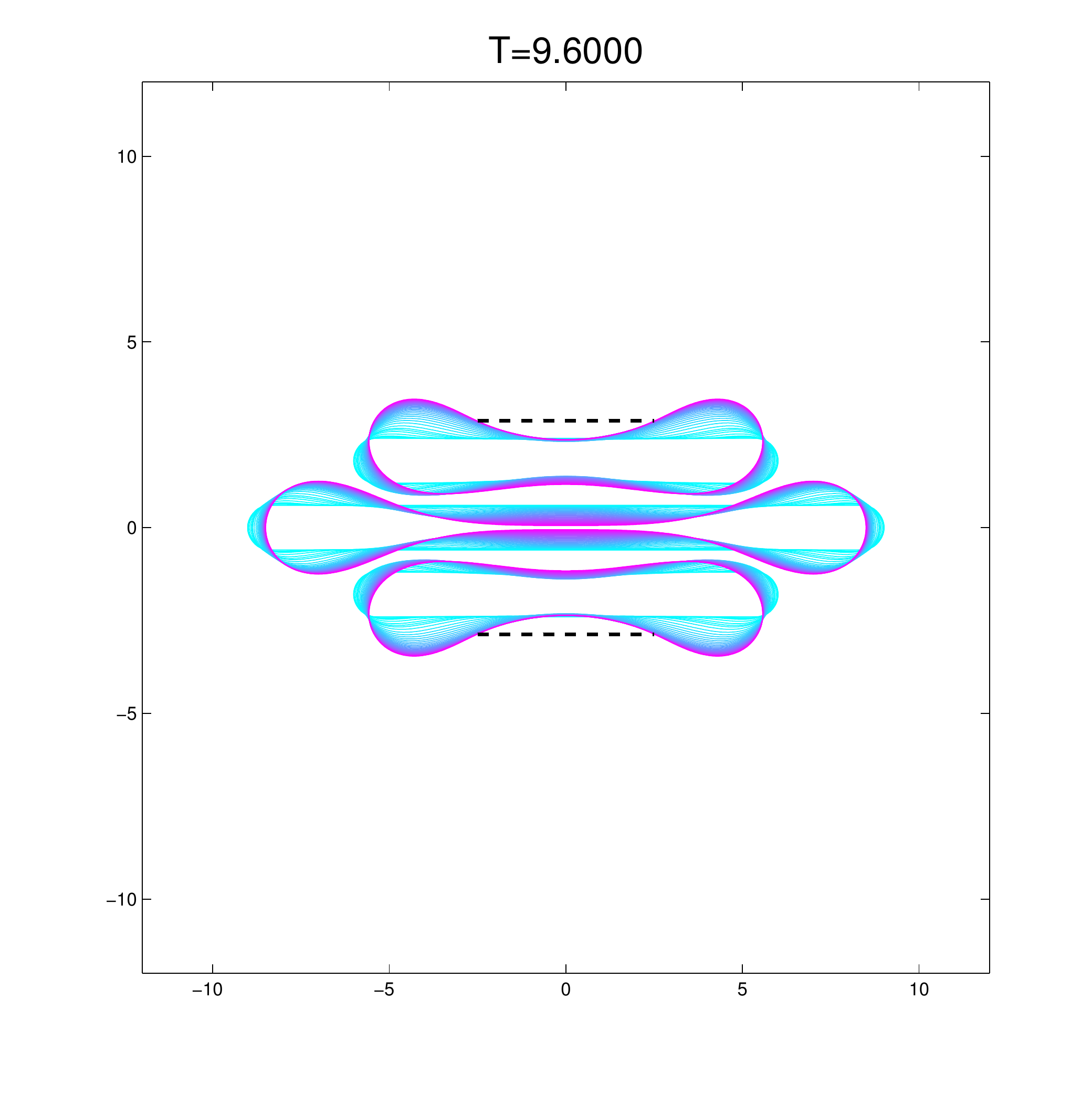}
}
\hspace{0mm}
\subfloat[Conservation of Length]{
\includegraphics[width=0.4\textwidth]{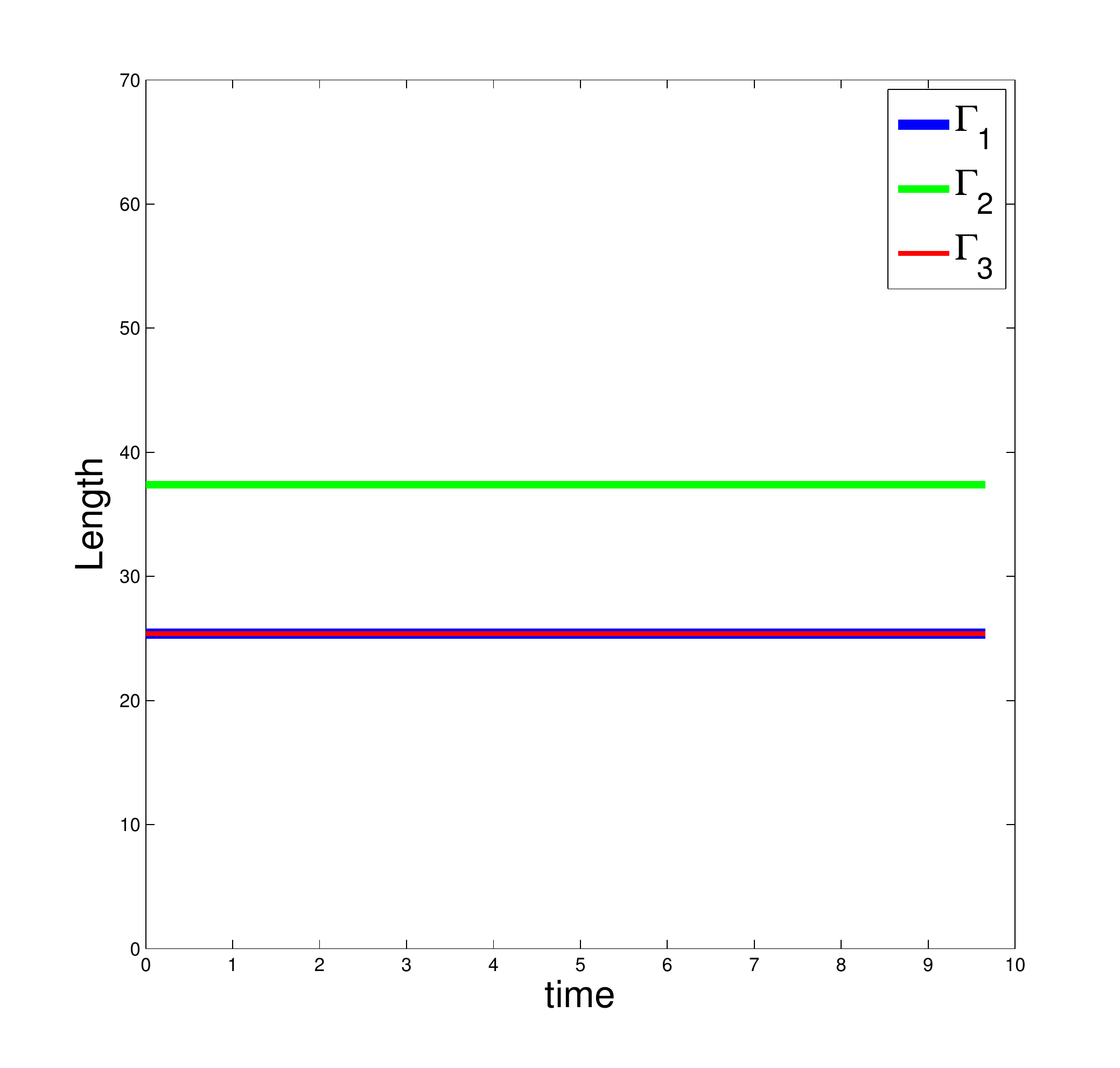}
}
\subfloat[Energy vs Time]{
\includegraphics[width=0.4\textwidth]{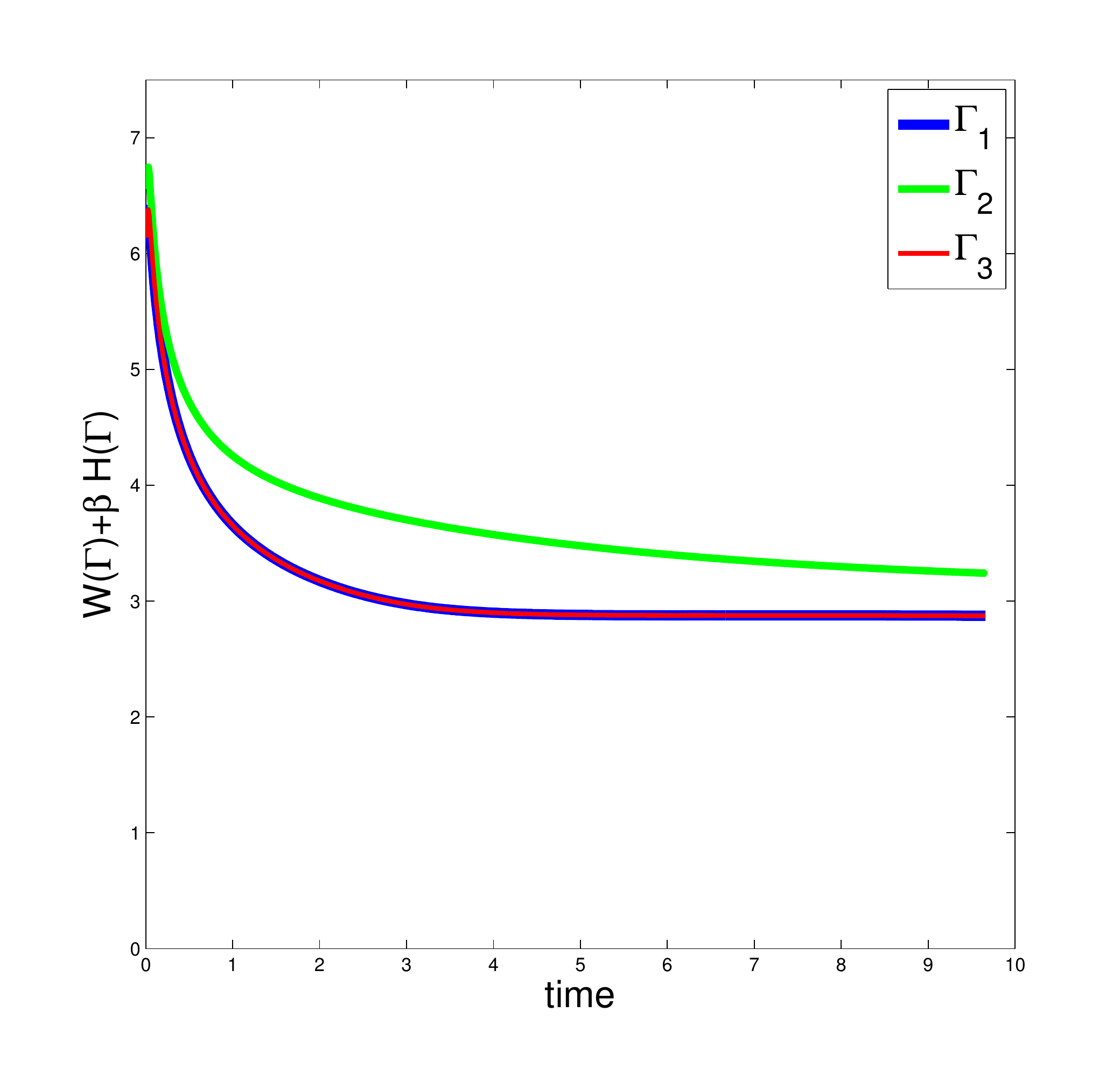}
}
\caption{Modeling of three cisternaes.}
\label{fig:L3}
\end{figure}

The details of the multiple vesicles case are illustrated in Section 3.3. The following example is an implementation of Model 3, with the initial shapes given as a set of parallelly ellipse-like shapes. In Figure \ref{fig:L3}, three cisternaes, named $\Gamma_1,\ \Gamma_2$ and $\Gamma_3$ from up to down are considered.

We do not draw conclusions about the Golgi stacks based on the modeling of multiple cisternaes, because the components of the Golgi stack is much more complicated. Not only the different biological properties of different stages of cisternaes make it complicated, but also the mechanism between them. Hence, we do these examples as a try only, though we can still find some similarity between our numerical results and the observed Golgi. For example, in Figure \ref{fig:Golgi}, the marginal parts of the upper half cisternaes go up, while those of the cisternaes below go down. Our numerical results also reveal this tendency. 

\begin{figure}[h!]
\centering
\includegraphics[width=0.45\textwidth]{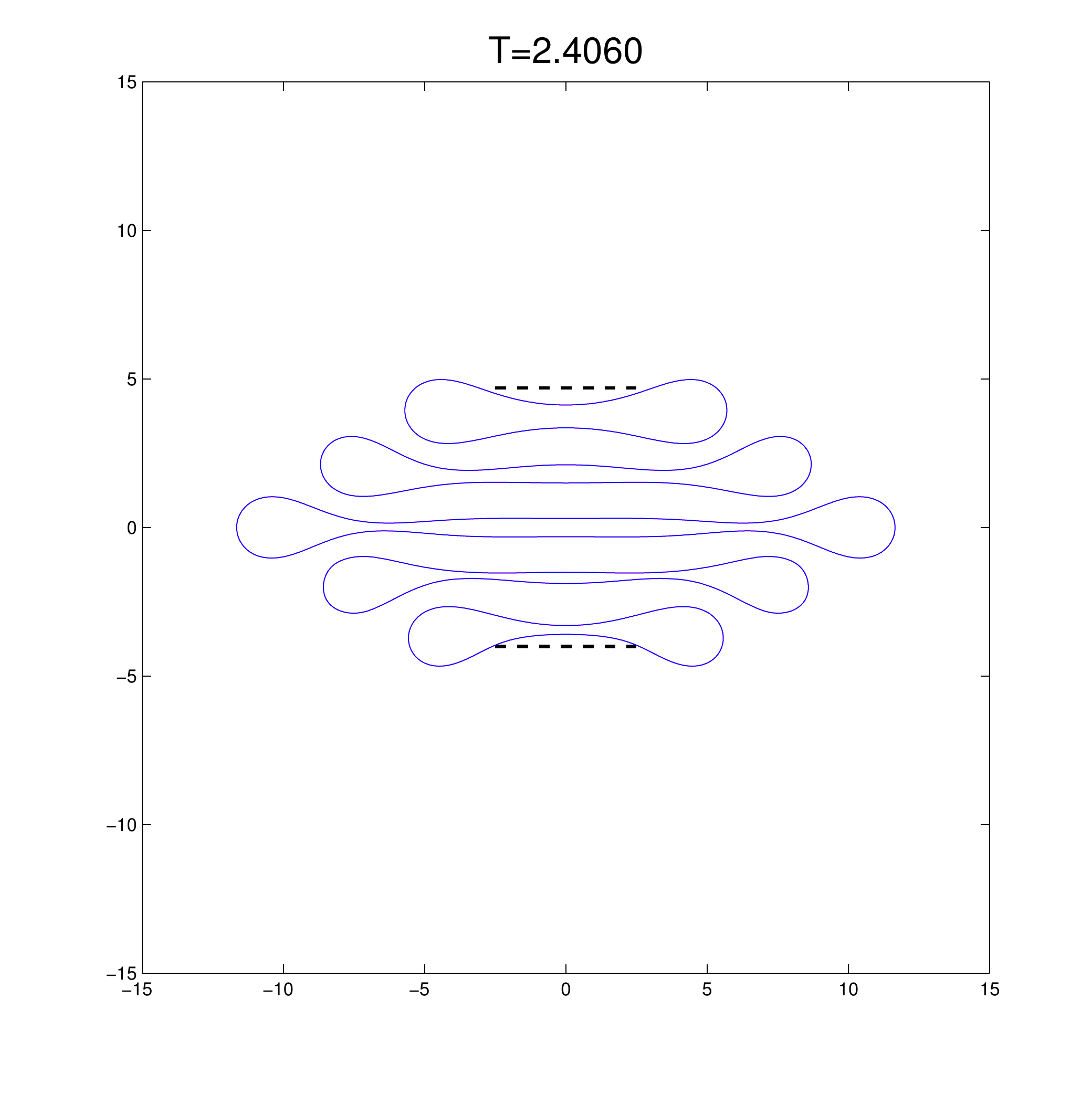}
\caption{Modeling of five cisternaes}	
\label{fig:L5}
\end{figure}

The simulation of more layers (5 cisternaes in Figure \ref{fig:L5}) gives similar results as those in Figure \ref{fig:L3}.

\end{document}